\input amstex
\documentstyle{amsppt}
\input bull-ppt
\keyedby{bull266e/jxs/pah}


  

\define\empha{\it}
\define\emphb{}
\define\emphc{\bf}

\define\refs#1:#2:#3:#4:#5:#6:#7:#8:{#1 #2, 
{\empha#3\/}, {\emphb#4} {\emphc#5} (#6), #7--#8.}
\define\refm#1:#2:#3:#4:#5:#6:#7:#8:#9:{#1 #2 #3, 
{\empha#4}, {\emphb#5} {\emphc#6} (#7), #8--#9.} 
\define\refmm#1:#2:#3:#4:#5:#6:#7:#8:#9:{#1 #2, #3, 
{\empha#4}, {\emphb#5} {\emphc#6} (#7), #8--#9.} 
\define\refbs#1:#2:#3:#4:#5:#6:{#1 #2, 
{\empha  #3}, #4,\ #5,  #6.} 
\define\refbm#1:#2:#3:#4:#5:#6:#7:{#1 #2 #3, 
{\empha  #4},\ #5,\ #6,\  #7.} 
\define\refps#1:#2:#3:{#1 #2, {\empha#3\/},\ preprint.}
\define\refpm #1:#2:#3:#4:{#1 #2  #3, {\empha#4\/},\ 
preprint.}
\define\refpmm #1:#2:#3:#4:{#1 #2,  #3, {\empha#4\/},\ 
preprint.}
\define\refprm #1:#2:#3:#4:{#1  #2, #3, {\empha#4\/},\ in 
preparation.}
\define\reftas #1:#2:#3:#4:{#1  #2, {\empha#3\/}, 
{\emphb#4}, to appear.}
\define\reftam #1:#2:#3:#4:#5:{#1 #2  #3, {\empha#4\/}, 
{\emphb#5}, to appear.}
\define\refts #1:#2:#3:#4:#5:{#1 #2, {\empha#3\/}, Thesis, 
#4, #5.}
\define\reflm#1:#2:#3:#4:{#1  #2  #3, {\empha#4\/}.}
\define\refres#1:#2:#3:#4:#5:{#1  #2, {\empha#4\/}, #4 
(#5).}
\define\reftesta#1:#2:#3:#4:#5:#6:{#1 #2, {\empha#3\/}, #4 
 (#5) (and to appear, #6).}


\define\amount{\hbox{\vrule height.1pt depth.2pt 
width.7truein}\kern.01truein}

\define\refsr#1:#2:#3:#4:#5:#6:#7:#8:{\amount, 
{\empha#3\/}, {\emphb#4} {\emphc#5}  (#6), #7--#8.}
\define\refmr#1:#2:#3:#4:#5:#6:#7:#8:#9:{
\amount, {\empha#4}, {\emphb#5} {\emphc#6}  (#7), #8--#9.} 
\define\refbsr#1:#2:#3:#4:#5:#6:{\amount, 
{\empha  #3}, #4,\ #5,  #6.} 
\define\refbmr#1:#2:#3:#4:#5:#6:#7:{\amount, 
{\empha  #4},\ #5,\ #6,\  #7.} 
\define\refpsr#1:#2:#3:{\amount, {\empha#3\/},\ preprint.}
\define\refpmr #1:#2:#3:#4:{\amount, {\empha#4},\ preprint.}
\define\refprmr #1:#2:#3:#4:{\amount, {\empha#4},\ in 
preparation.}
\define\reftasr #1:#2:#3:#4:{\amount,{\empha#3\/}, #4, to 
appear.}
\define\reftamr #1:#2:#3:#4:#5:{\amount, {\empha#4}, #5, 
to appear.}
\define\reftsr #1:#2:#3:#4:#5:{\amount, {\empha#3\/}, 
Thesis, #4, #5.}
\define\reflmr#1#2#3#4{\amount, {\empha#4}.}
\define\reftestar#1:#2:#3:#4:#5:#6:{\amount, {\empha#3\/}, 
#4  (#5) (and to appear, #6).}
\define\refresr#1:#2:#3:#4:#5:{\amount, {\empha#3\/}, #4 
(#5).}




\define\myProclaim#1. #2\par{\medbreak
\noindent{\bf#1.\enspace}{\it#2}\par\medbreak}

\define\m{{\cal S}(N)}
\define\secno{0} 
\define\trace{\text{ trace}} 
\define\hx{\hat x}
\define\hy{\hat y} 
\define\p{(x,y)} 
\define\sgl{_\gl}
\define\hp{\hat p}
\define\ph{(\hx,\hy)}
\define\oo#1{{1\over #1}}
\define\ot#1{{#1\over 2}} 
\define\usc#1{\roman{USC}(\R^{#1})}
\define\jtp{J^{2,+}} 
\define\jtm{J^{2,-}} 
\define\jtpo{J_\gO^{2,+}}
\define\jtmo{J_\gO^{2,-}} 
\define\jtpoc{\overline J_\gO^{2,+}}
\define\jtmoc{\overline J_\gO^{2,-}} 
\define\co{{\cal O}}
\define\jtpco{J_\co^{2,+}} 
\define\jtmco{J_\co^{2,-}}
\define\jtpcoc{\overline J _\co^{2,+}} 
\define\jtmcoc{\overline J_\co^{2,-}} 
\define\jtpc{\overline J ^{2,+}}
\define\jtmc{\overline J^{2,-}} 
\define\rsn{\RN\times\m}
 
\define\q#1{\ip{A#1,#1}}
\define\qp#1{\ip{A(#1),#1}}
\define\jt{{J}^2} 
\define\nos#1{|#1|^2}
\define\goc{\overline\gO} 
\define\bgo{\partial\gO}
\define\uscc{\roman{USC}(\goc)}
\define\lscc{\roman{LSC}(\goc)}
\define\gos{\gO\times\gO}
\define\xmys{\nos{x-y}}

\define\gocs{\goc\times\goc}
\define\ut{{}^t\kern -.08em} 
\define\md#1{{\cal S}(#1)}
\define\diag{\pmatrix X_1&\ldots&0\cr\vdots&\ddots&\vdots%
\cr0&\ldots&X_k\endpmatrix}
\define\jtpcoic{\overline J^{2,+}_{\scr O}}

\define\gS{\Sigma}
\define\sas{\hbox{$_*$}} 
\define\notes#1{\medskip \leftline{\bf Notes on Section 
#1}\smallskip} 
\define\set#1{\left\{#1\right\}} 
\define\dey{\delta y}
\define\dex{\delta x} 
\define\good{proper}
\define\Good{Proper}
\define\goodness{properness}
\define\hZ{\hat Z}

\define\pxij#1{#1_{x_i x_j}}
\define\sij{_{i,j}}
\define\spec#1{\gl_{1}(#1),\ldots,\gl_N(#1)}
\define\eig#1#2{\gl_{#1}(#2)}
\define\paren#1{\left(#1\right)}
\define\dt#1#2{{\partial^2#1\over\partial x_#2^2}}
\define\do#1#2{{\partial#1\over\partial x_#2}}
\define\dtm#1#2#3{{\partial^2#1\over\partial x_#2\partial 
x_#3}}
\define\cl{{\cal L}}

\define\ega{^\ga} 
\define\egab{^{\ga,\gb}} 
\define\ir#1{\ \hbox{\it #1}\ } 
\define\irq#1{\hbox{\it\quad #1\quad}}
\define\tx{\tilde x}
\define\sga{_{\ga}}
\define\sgan{_{\ga_n}}
\define\gan{\ga_n}
\define\junk{Z+\frac{1}{\|A\|}Z^2}
\define\ub{\underline u}
\define\ovp{\overline u}
\define\cf{{\cal F}}
\define\ws{w^*}

\define\us{u_*} 
\define\sgd{_{\gd,\gamma}}
\define\messf#1{\gl  M e^{-\gl d(#1)}Dd(#1)} 
\define\messt#1{\gl  M e^{-\gl 
d(#1)}D^2d(#1)-\gl^2Me^{-\gl d(#1)} Dd(#1)\otimes Dd(#1)} 
\define\messT#1{\gl   cDd(#1),\gl  cD^2d(#1)-\gl^2 c 
Dd(#1)\otimes Dd(#1)}
\define\psin{\psi_n}
\define\np#1{n^{#1}}
\define\hxn{\hat x_n}
\define\uc{\roman{UC}(\RN)}
\define\rtn{\RN\times\RN}
\define\gbr{\gb_R}
\define\bx{\overline x}
\define\ou{\overline u}
\define\uu{\underline u}
\define\hU{\overline U}
\define\uU{\underline U}
\define\ox{\overline  x}
\define\oy{\overline  y}
\define\bc{(BC)}
\define\jtpcgo{\overline J _{\goc}^{2,+}} 
\define\jtmcgo{\overline J _{\goc}^{2,-}}

\define\glp{\gl_+}
\define\glm{\gl_-}
\define\e#1{e^{#1}}
\define\sep{_\gep}
\define\tenp#1{#1\otimes #1} 
\define\cy{\co_T}


\define\dppc{\overline{\cal P}^{2,+}}

\define\dppco{{\cal P}_\co^{2,+}} \define\dpmco{{\cal 
P}_\co^{2,-}}
\define\dppcoc{\overline{\cal P}_\co^{2,+}}
\define\dpmcoc{\overline{\cal P}_\co^{2,-}}
\define\sa{_{\ga}}
\define\chatt{\hat t}
\define\tenpn#1{{#1\otimes #1\over\nos{ #1}}}
\define\uf{\underline F}
\define\of{\overline F}
\define\hw{\hat w}
\define\hu{\hat u}
\define\hv{\hat v}
\define\step#1{\noindent{\bf Step #1.\ }}
\define\sgf{\partial\gf}\define\dgf{D\gf}

\define\coi{{\cal O}_i}
\define\tpsi{\tilde\psi}
\define\cs{{\cal S}}

\define\ip#1{\langle #1\rangle}
\define\ti#1{\text{  #1  }}
\define\iq#1{\text{\quad #1\ }}
\define\R{\Bbb R}  
\define\C{\Bbb C}  
\define\vno#1{\|#1\|}

\define\RN{\R^N}

\define\ga{\alpha}
\define\gb{\beta}
\define\gG{\Gamma}
\define\gf{\varphi}
\define\gF{\Phi}
\define\gd{\delta}
\define\gD{\Delta}
\define\gep{\varepsilon}
\define\gm{\mu}

\define\gk{\kappa}
\define\gn{\nu}
\define\gl{\lambda}
\define\gL{\Lambda}

\define\gr{\rho}
\define\go{\omega}
\define\gO{\Omega}


\define\rtra{\rightarrow}

\define\downa{\downarrow}
\define\ua{\uparrow}

%
\newcount\equationno
\equationno=0
\newcount\resultno
\resultno=0
\newcount\sectionno
\sectionno=0
\newcount\bibno
\bibno=0




\define\marginal#1{{\vskip0pt\noindent\llap{\fiverm #1\ 
}\vskip-\baselineskip}}
\newif\ifshowlabel
\newif\ifremembermode

\define\eq{\global\advance\equationno by 1 
\leqno(\roman{\secno}.\the\equationno)}
\define\re#1{\ifshowlabel\marginal{\string#1}\fi
\ifremembermode\immediate\write\mem{\define\string#1
{(\secno.\the\equationno}}\fi 
\global\edef#1{(\secno.\the\equationno}}
\define\pre#1{{\advance\equationno by 1\re#1}}

\define\li{\global\advance\resultno by 1 
\secno.\the\resultno}
\define\rre#1{\ifshowlabel\marginal{\string#1}\fi
\ifremembermode\immediate\write\mem{\define\string#1
{\secno.\the\resultno}}\fi 
\global\edef#1{\secno.\the\resultno} }
\define\lref{\global\advance\resultno by 1 \the\resultno}
\define\rpre#1{{\advance\resultno by 1\rre#1}}

\define\sli{\global\advance\sectionno by 1 \the\sectionno}
\define\sre#1{\ifshowlabel\marginal{\string#1}\fi
\ifremembermode\immediate\write\mem{\define\string#1
{\the\sectionno}}\fi \global\edef#1{\the\sectionno} }
\define\spre#1{{\advance\sectionno by 1\sre#1}}

\define\bli{\global\advance\bibno by 1  [\the\bibno]}
\define\bre#1{\ifshowlabel\marginal{\string#1}\fi
\ifremembermode\immediate\write\mem{\define\string#1
{\the\bibno}}\fi \global\edef#1{\the\bibno} }
\define\bpre#1{{\advance\bibno by 1\bre#1}}



\define\equations{1}
\define\notions{2}
\define\maximum{3}
\define\perron{4}
\define\variations{5}
\define\rlimits{6}
\define\parabolic{8}
\define\boundary{7}
\define\singular{9}
\define\perspectives{10}
\define\tha{3.2}
\define\lej{A.3}
\define\thbcb{7.9}
\define\leat{A.4}
\define \alea {1}
\define \aleb {2}
\define \alec {3}
\define \alz {4}
\define \aiztom {5}
\define \baka {6}
\define \barbuc {7}
\define \barbu {8}
\define \barbua {9}
\define \barbdapa {10}
\define \bardi {11}
\define \bardiev {12}
\define \bardiper {13}
\define \barlesa {14}
\define \barlesb {15}
\define \barlesc {16}
\define \barlesgeo {17}
\define \barlese {18}
\define \barlesra {19}
\define \barlesf {20}
\define \barlesg {21}
\define \barlesbronsoug {22}
\define \barlesevsoug {23}
\define \barlespll {24}
\define \bapa {25}
\define \bapb {26}
\define \bapc {27}
\define \barlessoug {28}
\define \bara {29}
\define \barish {30}
\define \barjenpon {31}
\define \barjenm {32}
\define \barjena {33}
\define \barjenb {34}
\define \barjenfa {35}
\define \ben {36}
\define \bony {37}
\define \brezis {38}
\define \cafa {39}
\define \cansona {40}
\define \cansonb {41}
\define \cansoni {42}
\define \capev {43}
\define \capl {44}
\define \giga {45}
\define \cra {46}
\define \cel {47}
\define \craish {48}
\define \revisit {49}
\define \cla {50}
\define \clb {51}
\define \rclap {52}
\define \cln {53}
\define \cplld {54}
\define \clii {55}
\define \clinv {56}
\define \cls {57}
\define \cn {58}
\define \cnt {59}
\define \dupish {60}
\define \dupisha {61}
\define \dupishb {62}
\define \leneng {63}
\define \evco {64}
\define \evac {65}
\define \evpert {66}
\define \evgar {67}
\define \evish {68}
\define \evsonsoug {69}
\define \evsoug {70}
\define \evsougb {71}
\define \evsp {72}
\define \falcone {73}
\define \fed {74}
\define \fr {75}
\define \flsoa {76}
\define \flsougb {77}
\define \fra {78}
\define \gigis {79}
\define \gigi {80}
\define \gtr {81}
\define \ishiold {82}
\define \ish {83}
\define \ishm {84}
\define \ishiolda {85}
\define \isp {86}
\define \iseik {87}
\define \ishiirep {88}
\define \ishse {89}
\define \isbvp {90}
\define \ishiob {91}
\define \ishps {92}
\define \ishkoia {93}
\define \ishkob {94}
\define \ishiikoike {95}
\define \lii {96}
\define \ishyam {97}
\define \ivan {98}
\define \koia {99}
\define \koib {100}
\define \jen {101}
\define \jenb {102}
\define \jenr {103}
\define \jls {104}
\define \jensoug {105}
\define \kn {106}
\define \kru {107}
\define \nvka {108}
\define \nvkb {109}
\define \lali {110}
\define \lalib {111}
\define \lena {112}
\define \lenc {113}
\define \lenyama {114}
\define \lenyamb {115}
\define \plla {116}
\define \pllb {117}
\define \pllc {118}
\define \pllv {119}
\define \pllbony {120}
\define \pllneumann {121}
\define \pllreg {122}
\define \pll {123}
\define \pla {124}
\define \pllcid {125}
\define \plld {126}
\define \plle {127}
\define \pllpap {128}
\define \pllper {129}
\define \pln {130}
\define \pllro {131}
\define \pllso {132}
\define \pllsoa {133}
\define \pllsov {134}
\define \mignot {135}
\define \newc {136}
\define \nuna {137}
\define \nunb {138}
\define \olr {139}
\define \osse {140}
\define \pera {141}
\define \persan {142}
\define \pucci {143}
\define \roto {144}
\define \rudin {145}
\define \satoa {146}
\define \yam {147}
\define \sayaa {148}
\define \sona {149}
\define \sonb {150}
\define \sonc {151}
\define \sond {152}
\define \taka {153}
\define \takc {154}
\define \takd {155}
\define \takb {156}
\define \tata {157}
\define \trua {158}
\define \trub {159}
\define \truc {160}
\define \zara {161}
\define \zarb {162}

\topmatter
\cvol{27}
\cvolyear{1992}
\cmonth{July}
\cyear{1992}
\cvolno{1}
\cpgs{1-67}

  \title  user's guide to viscosity solutions\\ of
second order\\ partial differential equations\endtitle
\author Michael G. Crandall, Hitoshi Ishii, and 
Pierre-Louis Lions\endauthor
\shortauthor{M. G. Crandall, Hitoshi Ishii, and 
Pierre-Louis Lions}
\shorttitle{User's guide to viscosity solutions}
\address Department of Mathematics, University of 
California, 
Santa Barbara, California 93106\endaddress
\address Department of Mathematics, Chuo University, 
Bunkyo-ku, 
Tokyo 112, Japan\endaddress
\address Ceremade, Universit\'e Paris-Dauphine, Place de 
Lattre 
de Tassigny, 75775 Paris Cedex 16, France\endaddress 
\thanks First author supported in part by the
Army Research Office  DAAL03-87-K-0043 and 03-90-G-0102,
National Science Foundation  DMS-8505531 and 90-02331, and
Office of Naval Research   N00014-88-K-0134\endthanks 
\date November 16, 1990\enddate
\thanks This paper was given as a Progress in Mathematics
Lecture at the August 8--11, 1990 meeting of the American
Mathematical Society in Columbus, Ohio\endthanks
\abstract The notion of viscosity solutions of scalar
fully nonlinear partial differential equations of second
order provides a framework in which startling
comparison and uniqueness theorems, existence theorems,
and theorems about continuous dependence may now be
proved by very efficient and striking arguments. The
range of important applications of these results is
enormous. This article is a self-contained exposition
of the basic theory of viscosity solutions.\endabstract
\endtopmatter

\document

\heading Introduction\endheading

The
theory of viscosity solutions   applies to certain partial 
differential
equations of the form $F(x,u,Du,D^2u)=0$ where
$F\:\RN\times\R\times\RN\times\m\rtra\R$ and $\m$ is the 
set of
symmetric $N\times N$ matrices.    The primary virtues  of 
this
theory are that it allows merely continuous functions to be 
solutions of fully nonlinear equations of second order, 
that it
provides very general existence and uniqueness theorems 
and that it
yields precise formulations of general boundary conditions. 
Moreover, these features go hand-in-hand with a great 
flexibility
in passing to limits in various settings and relatively 
simple
proofs.  In the expression $F(x,u,Du,D^2u)$  
$u$ will  be a real-valued
function defined on some subset $\co$ of $\RN$,  $Du$ 
corresponds
to the gradient of $u$ and  $D^2u$ corresponds to the 
matrix of
second derivatives of $u$.  However, as explained below, 
$Du$ and
$D^2u$ will not have classical meanings and, in fact, the 
examples
exhibited in \S \equations\ will convince the reader that  
the
theory encompasses classes of equations that have no 
solutions
that are differentiable in the classical sense.

In order that the theory apply to a given equation $F=0$, 
we will
require  $F$ to  satisfy a fundamental  monotonicity 
condition 
\pre\eqaa
$$F(x,r,p,X)\le F(x,s,p,Y)\iq{whenever} r\le s\ti{and} 
Y\le X;\eq$$ 
where $r,s\in\R$, $x,p\in\RN$, $X,Y\in\m$ and $\m$ is
equipped with its usual order.   

Regarding \eqaa) as made up of the
two conditions
\pre\eqzzz
$$
F(x,r,p,X)\le F(x,s,p,X)\quad\text{whenever } r\le s, \eq$$
and 
\pre\eqwww
$$F(x,r,p,X)\le F(x,r,p,Y)\iq{whenever}  Y\le X,\eq$$ 
we will give the name ``degenerate ellipticity" to the 
second.  That is,
$F$ is said to be {\it degenerate elliptic} if \eqwww) 
holds.  When \eqzzz)
also holds (equivalently, \eqaa) holds), we will say that 
$F$ is {\it \good}.

The examples given in \S\equations\ will illustrate
the fact that the antimonotonicity in $X$ is indeed an
``ellipticity" condition.    The possibility of 
``degeneracies"  is
clearly exhibited by considering the case  in which 
$F(x,r,p,X)$ does
not depend on $X$---it is then degenerate elliptic.  The 
monotonicity in
$r$, while easier to understand, is a slightly subtle 
selection criterion
that, in particular,  excludes the use of the viscosity 
theory for first
order equations of the form $b(u)u_x=f(x)$ in $\R$ when 
$b$ is not a
constant function, since then $F(x,r,p)=b(r)p-f(x)$ is not 
nondecreasing
in $r$ for all choices of $p$ (scalar conservation laws 
are outside of the
scope of this theory).   

The presentation begins with \S\equations, which, as already
mentioned, provides a list of examples.   This rather long 
list is
offered to meet several objectives.  First, we seek  to 
bring the
reader to our conviction that the scope of the theory is 
quite
broad while providing a spectrum of meaningful 
applications and, at
the same time, generating some insight as regards the 
fundamental
structural assumption \eqaa).   Finally, in the 
presentation of
examples involving famous second order equations,  the 
very act of
writing the equations in a form compatible with the theory 
will
induce an interesting modification of the classical 
viewpoint concerning
 them.

In \S\notions\ we begin an introductory presentation of the
basic facts of the theory.  The style is initially 
leisurely and
expository and technicalities are minimized, although 
complete
discussions of various key points are given and some simple
arguments inconveniently scattered in the literature are
presented.  Results are illustrated with simple examples 
making
clear their general nature. Section \notions\  presents 
the basic 
notions of solution used in the theory, the analytical 
heart of
which lies in comparison results.  Accordingly, \S\maximum  
\ is devoted to explaining  comparison results in the 
simple setting
of the Dirichlet problem;  roughly speaking, they are  
proved by
  methods involving extensions of the maximum principle  to
semicontinuous functions.   Once  these comparison results 
are
established, existence assertions can be established by
 Perron's method, a rather striking tale that is told in \S
\perron.  With this background in hand, the reader will 
have an
almost complete (sub)story  and with some effort (but not 
too
much!) should be able to absorb in an efficient way some 
of the
 more technical features of the theory that are outlined
in the rest of the paper. 

 Other important ideas are to be found in  \S\rlimits, 
which  
 is concerned with the issue of taking limits of viscosity
solutions and applications of this and in \S\boundary, which
describes the adaptation of the theory  to accommodate 
problems
with other boundary conditions and problems in which the 
boundary
condition cannot be strictly satisfied.  In the later 
case, the
entire problem has a generalized interpretation for which 
there is
often existence and uniqueness.   While the description of 
these
results is deferred to \S\boundary, they are fundamental and
dramatic.  For example, if $G(p,X)$ is uniformly  
continuous,
degenerate elliptic, and independent of $x$ and 
$\gO\subset\RN$ is
open and bounded, $n(x)$ is the unit exterior normal at 
the point $x$ in 
its smooth boundary $\partial \gO$ and  $f\in C(\goc)$,  
then the Neumann
problem 
$$u+G(Du,D^2u)-f(x)=0 \quad\text{in }{\gO},\qquad
u_n=0\quad\text{on }\partial\gO$$ 
has a unique properly interpreted solution
(which may not satisfy $u_n=0$ on $\partial\gO$).  

Sections \variations, \parabolic, \singular\  discuss 
variations of
the basic material and need not be read in sequence.  
Section
\perspectives \  is devoted to a commentary about 
applications
(which are  not treated in the main text), references, and 
possible lines of 
future development of the subject.  We conclude with an 
appendix where
the reader will find  a self-contained presentation of the 
proof of the
analytical heart on the presentation we have chosen.

References are not given in the main text (with the 
exception of 
\S \perspectives),  but are to be  found at
the end of each section.   In particular, the reader 
should look to
the end of a section for further comments,  references that
contain details ommitted in the main text, and technical 
generality,
historical comments, etc.  The references, while numerous, 
are 
not intended to be complete, except that we have sought to 
represent 
all the major directions of research and areas of 
application. 
There are original aspects 
of the current presentation and the 
reader will note differences between the flavor and 
clarity of our
presentation and that of many of the citations.  However, 
equipped with the view presented here, we hope and expect 
that perusing the 
amazing literature that has so quickly matured will be a 
much more 
rewarding and efficient endeavor.

We are grateful to R. Dorroh, M. Kocan and A. Swiech for 
their kind help in 
reducing the number of errors herein.


\ac
\widestnumber\head{10}
\widestnumber\subhead{7.C\cprime}
\heading 1. Examples\endheading
\heading 2. The notion of viscosity solutions\endheading
\heading 3. The maximum principle  for semicontinuous 
functions and
comparison for the Dirichlet problem\endheading
\heading 4. Perron's method and existence\endheading
\heading 5. Comparison:  Variations on the theme\endheading
\subheading{5.A. Comparison with more regularity}
\subheading{5.B. Estimates from comparison}
\subheading{5.C. Comparison with strict inequalities and 
without
coercivity in $u$}
\subheading{5.D. Comparison and existence of unbounded 
solutions on 
unbounded domains}
\heading 6. Limits of viscosity solutions and an 
application\endheading
\heading 7. General and generalized boundary 
conditions\endheading
\subheading{7.A. Boundary conditions in the viscosity sense}
\subheading{7.B. Existence and uniqueness for the Neumann 
problem}
\subheading{7.C. The generalized Dirichlet problem }
\subheading{7.C\cprime. The state constraints problem}
\subheading{7.D. A remark on (BC) in the classical sense}
\subheading{7.E. Fully nonlinear boundary conditions}
\heading 8. Parabolic problems\endheading
\heading 9. Singular equations: An example from 
geometry\endheading
\heading 10. Applications and perspectives\endheading
\heading {} APPENDIX The proof of Theorem 3.2\endheading
\endac

\resultno=0 \equationno=0 \spre\equations 
\heading \sli. Examples \endheading
\def\secno{\the\sectionno} 

We  will record here many examples of degenerate elliptic 
equations
mentioning, when appropriate, areas in which they arise.  
The
reader is invited to scan the list and pause where 
interested---it is 
possible to proceed to \S \notions\ at any stage.  Below
we will say either that a function $F(x,r,p,X)$ is 
degenerate
elliptic  or that the equation $F(x,u,Du,D^2u)=0$ is 
degenerate
elliptic or that the ``operator" or expression 
$F(x,u,Du,D^2u)$ is
degenerate elliptic and always mean the same thing, i.e.,  
\eqwww)
holds; the term ``\good" is used in a similar fashion.  
 
\rpre\exta
\ex{Example  \li}  Laplace's equation. 
We revisit
an old friend, the equation 
$$-\gD u+c(x)u=f(x)\eq$$ 
\re\eqexa---note
the sign in front of the Laplacian.    The corresponding 
$F$ is naturally
given by $F(x,r,p,X)=-\trace(X)+c(x)r-f(x)$, which is \good 
\ if $c\ge0$.  
\endex
 
\rpre\extb 
\ex{Example\ \li}  Degenerate
elliptic linear equations. 
Example \exta\  immediately extends to the more general 
linear
equation 
\pre\eqexb
$$-\sum_{i,j=1}^N a_{i,j}(x)\dtm u i
j+\sum_{i=1}^Nb_i(x)\do  u i+c(x)u(x)=f(x)\eq$$
  where the
matrix  $A(x)=\{a\sij(x)\}$ is  symmetric; the 
corresponding $F$ is
\pre\eqexc
$$F(x,r,p,X)=-\trace(A(x)X)+\sum_{i=1}^N b_i(x)p_i+
c(x)r-f(x).\eq$$

 In this case, $F$ is degenerate elliptic if and only if  
$A(x)\ge 0$ and it is \good\ if also $ c(x)\ge0$.   In the 
event that
there are constants $\gl,\gL>0$ such that $\gl I\le 
A(x)\le \gL I$ for all
$x$ where $I$ is the identity matrix, $F$ is said to be 
{\it uniformly
elliptic}. 

Of course, the linear equation in divergence form 
$$-\sum_{i,j=1}^N
{\partial\over\partial x_i}\left(a_{i,j}(x)\do u j\right)
+\sum_{j=1}^Nb_j(x)\do  u i+c(x)u(x)=f(x)$$ 
can be written as above
with  
$$F(x,r,p,X)=-\trace(A(x)X)+
\sum_{j=1}^N\left(b_j(x)-\sum_{i=1}^N\do
{a_{i,j}}i(x)\right)p_j+c(x)r-f(x)$$ 
provided that the indicated
derivatives of the $a_{i,j}$ exist.
\endex

We leave the interesting class of linear equations to turn 
to the
totally degenerate case of first order equations.

\rpre\extc 
\ex{Example\  \li}  First order
equations. 
The main point is  that a first order operator
$F(x,u,Du)$ is always degenerate elliptic and thus it is 
\good\ if and  only
if $F(x,r,p)$ is nondecreasing in $r\in\R$.  \Good\  
equations of
the form $F(x,u,Du)=0$ play a  fundamental role in the 
classical
Calculus of Variations and in Optimal Control Theory of 
ordinary
differential equations; in this context they are often 
called
Bellman or Hamilton-Jacobi equations and then $F(x,r,p)$ 
is convex
in $(r,p)$.  These equations, in the full generality of 
nearly
arbitrary \good\ functions $F$, are also crucial in 
Differential Games
Theory where they are known as Isaacs's equations. 
\endex

\rpre\extd 
\ex{Example\  \li}  Quasilinear
elliptic equations in divergence form. 
The usual notion of ellipticity for equations of the form 
$$-\sum_{i=1}^N  {\partial\ \ \over\partial
x_i}\left(a_i(x,Du)\right)+b(x,u,Du)=0\eq$$ 
\re\eqexd 
is the
monotonicity of the vector field $a(x,p)$ in  $p$ as a 
mapping from
$\RN$ to $\RN$.    If enough regularity is available to 
carry out
the differentiation, we write \eqexd) as 
$$-\sum_{i,j=1}^N{\partial a_i \over\partial 
p_j}(x,Du)\dtm u i j
+b(x,u,Du)-\sum_{i=1}^N{\partial a_i \over\partial
x_i}(x,Du)=0\eq$$ 
\re\eqexe 
and correspondingly set 
$$F(x,r,p,X)=-\trace((D_p a(x,p))X)+b(x,r,p)- 
\sum_{i=1}^N{\partial
a_i \over\partial x_i}(x,p).\eq$$ 
The monotonicity of $a$ in $p$ is
precisely the condition that guarantees that $F$ is 
degenerate
elliptic, and then it is \good\ provided we ask that $b$ 
be nondecreasing
in $r$.

Two well-known instances are provided by the equation of 
minimal
surfaces
 in nonparametric form and the ``$m$-Laplace's" equation 
that are
given, respectively,  by $a(x,p)=(1+|p|^2)^{-1/2}p,\ 
b=b(x,u)$ and
$a(x,p)=|p|^{m-2}p,\ b=b(x,r)$ where $m\in(1,\infty)$. 
Computations show that the corresponding $F$'s are, 
respectively, 
$$F(x,r,p,X)=-(1+|p|^2)^{-{1\over 2}}\trace (X)+(1+
|p|^2)^{-{3\over
2}}\trace ((p\otimes p)X)+b(x,r)$$ 
and
$$F(x,r,p,X)=-|p|^{m-2}\trace (X)-(m-2)|p|^{m-4}\trace 
((p\otimes
p)X)+b(x,r).$$ 
\endex

\rpre\extnd 
\ex{Example\  \li}  Quasilinear
elliptic equations in nondivergence form. 
The equation
$$-\sum_{i,j=1}^N a\sij(x,p)\dtm u ij+ b(x,u,Du)=0,\eq$$ 
\re\eqexex
where $A(x,p)=\{a\sij(x,p)\}\in\m$, contains all of the 
above as
special cases and  corresponds to 
$$F(x,r,p,X)=-\trace(A(x,p)X)+b(x,r,p),$$ 
which is \good\  if $A\ge 0$ and
$b$ is nondecreasing with respect to $r$.  Two relevant 
special cases are
$$-\gn\gD u+f(x,u,Du)=0$$ 
with $\gn>0$ and $f$ nondecreasing in $u$,
which may be regarded as a first-order Hamilton-Jacobi 
equation
perturbed by an additional ``viscosity" term $-\gn\gD u$ 
(equations of
this type arise in optimal stochastic control),   and the 
L\'evi's equation 
$$
\align
&-\paren{\dt u 1+\dt u 2}\paren{1+\paren{\do u 3}^2}-\dt
u 3\paren{\paren{\do u 1}^2+\paren{\do u 2}^2}\\
&\qquad+ 2\dtm u 1
3\paren{\do u 3\do u 1-\do u 2}+2\dtm u 2 3\paren{\do u 
3\do u
2+\do u 1}=0,\endalign$$ 
which is the nonparametric formulation for a
hypersurface in $\C^2$ with vanishing L\'evi's form.  Note 
that in
this example  $F=-\trace(A(p)X)$ where
$$A(p)=\pmatrix 1+p_3^2&
0&p_3p_1-p_2\\ 0&1+p_3^2&p_3p_2+p_1 \\ 
p_3p_1-p_2&p_3p_2+p_1&p_1^2+p_2^2\endpmatrix$$ 
so that $A\ge 0$ but
det($A(p))=0$ for all $p$.
\endex

 \rpre\exte 
\ex{Example\ \li}  
Hamilton-Jacobi-Bellman and Isaacs equations. 
Hamilton-Jacobi-Bellman and Isaacs equations are, 
respectively, 
the fundamental  partial differential equations for 
stochastic
control and stochastic differential games.  The natural 
setting
involves a collection of elliptic operators of second-order
depending either on one parameter $\ga$ (in the
Hamilton-Jacobi-Bellman case) or two parameters $\ga,\gb$ 
(in the
case of Isaacs's equations).   These parameters lie in 
index sets we
will not display in the discussion.   Thus we take as 
ingredients
\good\  expressions of the form  
$$\cl\ega u=-\sum_{i,j=1}^N a_{i,j}\ega(x)\dtm u i j+
\sum_{i=1}^Nb_i\ega(x)\do  u
i+c\ega(x)u(x)-f\,\ega(x)\eq$$  
\re\eqexf 
or 
$$\cl\egab u=-\sum_{i,j=1}^N
a_{i,j}\egab(x)\dtm u i j+\sum_{i=1}^Nb_i\egab(x)\do  u
i+c\egab(x)u(x)-f\,\egab(x)\eq$$ 
 \re\eqexg
\noindent where 
 all the coefficients are bounded with
respect to $\ga$ or $\ga,\gb$.  Hamilton-Jacobi-Bellman 
equations
include those of the form 
$$\sup_\ga\left\{\cl\ega u\right\}=0\eq$$
\re\eqxh 
while  
$$\sup_\ga\inf_\gb\left\{\cl\egab u\right\}=0\eq$$
\re\eqxi 
is a typical Isaacs's equation.  The corresponding
nonlinearities $F$ have the form
$$F(x,r,p,X)=\sup_\ga\left[-\trace(A\ega(x)X)+
\ip{b\ega(x),p}+c\ega(x)r-
f\,\ega(x)\right]$$ 
and
$$
\align
F(x,r,p,X)=\sup_\ga\inf_\gb[ & -\trace(A\egab(x)X)\\
& +\ip{b\egab(x),p}+c\egab(x)r
-f\,\egab(x)],\endalign$$ 
each of which is clearly also \good.
Notice that in the first case $F$ is convex in $(r,p,X)$ 
while in the
second case this is not so; indeed, if one allows for 
``unbounded
envelopes" (i.e., coefficients that are unbounded in the
parameters), one can show that essentially any \good
\ $F$ satisfying minor regularity assumptions can be 
represented as a
``sup inf"  of linear expressions as above.
\endex

Indeed, the above process is quite general and does not 
require
linear ingredients. Suppose $F\egab$ is \good\    for each
$\ga,\gb$.  Then  $F(x,r,p,X)=\sup_\ga\inf_\gb 
F\egab(x,r,p,X)$ and
$F(x,r,p,X)= \inf_\ga\sup_\gb F\egab(x,r,p,X)$ are also 
evidently
\good\  (for the moment, we set aside considerations of
finiteness and continuity).   

\rpre\extf 
\ex{Example\ \li}  Obstacle and gradient constraint
problems. 
A special case of the last remarks above is met in the
consideration of ``obstacle problems."  Very general forms 
of such 
problems may be written  
$$\max\{F(x,u,Du,D^2u), u-f(x)\}=0\eq$$
\re\eqexj 
or  
$$\min\{F(x,u,Du,D^2u), u-f(x)\}=0\eq$$ 
\re\eqexk 
or even 
$$\max\{\min\{F(x,u,Du,D^2u), u-f(x)\},u-g(x)\}=0.\eq$$
\re\eqexl 
In accordance with remarks  made in the previous example,
if $F$ is \good\  then so are \eqexj)--\eqexl).

Likewise, ``gradient constraints" may be imposed in this 
way.   A
typical example corresponds to 
$\max\{F(x,u,Du,D^2u),|Du|-g(x)\}=0$.
\endex

\rpre\extg 
\ex{Example\ \li}  Functions of the
eigenvalues. 
For $X\in\m$ we let 
$\lambda_1(X),\dotsc,
\lambda_N(X)$
be its
eigenvalues arranged in increasing order, $\eig i X\le\eig
{i+1}X$.  If
$g(x,r,p,s_1,\ldots,s_N)$ is defined on
$\R^{3N+1}$ and is nondecreasing in $r$ and each
$s_i$,
then $F(x,r,p,X)=g(x,r,p,-\eig 1 X,\ldots,-\eig N X)$ is 
\good.   
 For
instance, $F(X)=-\max\{\spec X\}=-\eig N X$,
$F(X)=
-\min\{\spec X\}=-\eig 1 X$ and $F(X)=-(\eig 2 X)^3$ are
degenerate elliptic.  Another example is 
$$F(x,r,p,X)=-|\trace(X)|^{m-1}\trace (X)+|p|^q+
c(x)r-f(x)$$ 
where
$c\ge0$ and $m,q>0$.  The corresponding equation is 
$$-|\gD u|^{m-1}\gD u+|Du|^q+c(x)u=f(x),$$ 
which provides another example of
the generality we are dealing with, even if we have no
interpretation of this equation in mind.
\endex

\rpre\exth 
\ex{Example\ \li}   Sums and
increasing  functions of proper  functions. 
If $F_i$ is \good\   for $i=1,\ldots, M$, then so is
$F_1+\dotsb+F_M$.  More generally,  if $g(s_1,\dotsc,s_M)$ 
is nondecreasing
in each variable, then $g(F_1,\dotsc,F_M)$ is \good.  
One may build very complex
examples using the cases discussed above and these remarks.
\endex

\rpre\exti 
\ex{Example\ \li} Parabolic
problems. 
We just observe that if $(x,r,p,X)\rtra F(t,x,r$, $p,X)$ is 
\good\  for fixed
$t\in[0,T]$, then so is the associated ``parabolic" problem
$$u_t+F(t,x,u,Du,D^2u)=0\eq$$ 
\re\eqexm 
\noindent when considered as an
equation in the $N+1$ independent variables $(t,x)$. We 
mention only one
example (there are, of course, infinitely many) that has 
some
geometrical interest since it describes the evolution of a 
surface (given
by a level set of the initial condition) with a motion 
along its normal with
a speed proportional to the mean curvature
$$u_t-|Du|\text{
div}\left({Du\over |Du|}\right)=0.\eq$$ 
\re\eqexn 
Carrying out the
differentiations yields 
$$u_t-\gD u+\sum_{i,j=1}^N\dtm u i j\do  u
i\do u j|Du|^{-2}=0.\eq$$ 
This may be written in the form \eqexm)
with  
$$F(x,p,X)=-\trace\left(\left(I-{p\otimes p\over
|p|^2}\right) X\right).$$  
\endex

\rpre\exti 
\ex{Example\ \li} Monge Amp\`ere
equations. 
The Monge-Amp\`ere equation may be written as 
$$u\iq{is convex,} \qquad
\det(D^2u)=f(x,u,Du)\eq$$ 
\re\eqexo 
where $f(x,r,p)\ge 0$.  We
are dealing here with the real Monge-Amp\`ere equation, but
everything that will be said adapts to the complex case 
and to
other curvature equations.  Allowing $F$ to be 
discontinuous (even
more, to become infinite), we may write  \eqexo) in our 
form by
putting 
$$F(x,r,p,X)=\cases -\det (X)+f(x,r,p)&\text{if $X\ge
0,$}\cr+\infty&\text{otherwise};\endcases$$ 
$F$ is then degenerate elliptic.
This follows from the fact that 
$$g(s_1,\cdots,s_N)=\cases  \prod_{i=1}^N s_i&\text{if 
$s_i\ge 0,  
\ i=1,\ldots,N$,}\cr -\infty & \text{otherwise,}\endcases$$
is nondecreasing in each of its arguments and Example \extg.
\endex

\rpre\extj 
\ex{Example\ \li}  Uniformly
elliptic functions. 
This ``example" is really a definition.  One says that 
$F(x,r,p,X)$
is {\it uniformly elliptic} if there are constants $\gl,
\gL>0$ for which 
$$\gl\trace (P)\le F(x,r,p,X-P)-F(x,r,p,X)\le\gL\trace 
(P)\iq{for}  P\ge 0$$ 
and all $x,r,p,X$ 
and then  calls the constants $\gl,\
\gL$  ellipticity constants.   We note that $L$ in Example
\extb\ is uniformly elliptic with constants $\gl,\ \gL$ 
exactly when
$\gl\le\eig 1{A(x)}$ and $\eig N{A(x)}\le\gL$.  One notes 
that sums
of uniformly elliptic functions are again uniformly 
elliptic and
that the sup inf process over a family of uniformly elliptic
functions with  common ellipticity constants produces 
another such
function.
\endex

\subheading{Notes on \S 1}
We will disappoint the reader in the following sections by 
not 
applying the theory developed therein to the many examples 
given above.
The goal of this section was to  exhibit clearly the 
breadth and importance 
of the class of \good\ equations.  We simply do not have 
enough space here  to
develop applications of the theory of
these equations beyond that which follows immediately from 
the general results 
presented.

Most of  the examples listed have been considered 
via classical approaches.  We give some references 
containing 
classical presentations:  D. Gilbarg and N. S. Trudinger 
[\gtr] is a 
basic source concerning linear and quasi-linear uniformly 
elliptic equations;
O. A. Ole\u\i nik and E. V. Radkevic [\olr],  J. J. Kohn 
and L. Nirenberg [\kn], 
and A. V. Ivanov [\ivan] treat degenerate 
elliptic equations;  W. H.  Fleming and R. Rishel [\fr], 
P. L. Lions [\pllb], and 
N. V. Krylov [\nvka,
\nvkb] are sources for Hamilton-Jacobi-Bellman equations; 
S. Benton [\ben] and P. L. Lions [\plla] 
discuss first-order Hamilton-Jacobi equations.  Most of 
these references present 
some of the ways these equations arise.

\resultno=0 \equationno=0 \spre\notions  
\heading \sli. The notion of viscosity solutions 
\endheading 
\def\secno{\the\sectionno} 

It is always assumed that $F$  satisfies  \eqaa) (i.e., 
$F$ is  
\good ) and, {\it unless otherwise said}, is   
continuous.   To motivate the notions, we begin by  
supposing that
$u$ is  $C^2$ (i.e., twice continuously differentiable) on 
$\RN$
and 
$$
F(x,u(x),Du(x),D^2u(x))\le 0$$ 
holds for all $x$ (that is, $u$
is a classical subsolution of $F=0$ or, equivalently, a 
classical
solution of $F\le 0$ in $\RN$).  Suppose that $\gf$ is 
also $C^2$
and $\hx$ is a local maximum of $u-\gf$.    Then calculus 
implies
$Du(\hx)=D\gf(\hx)$ and $D^2u(\hx)\le D^2\gf(\hx)$ and so, 
by
degenerate ellipticity, 
$$F(\hx,u(\hx),D\gf(\hx),D^2\gf(\hx))\le
F(\hx,u(\hx),Du(\hx),D^2u(\hx))\le 0.\eq$$ 
\re\eqaax 
\noindent The extremes
of this inequality do not depend on the derivatives of $u$ 
and so
we may consider defining an arbitrary function $u$ to be 
(some kind
of generalized)  subsolution  of $F=0$  if 
$$F(\hx,u(\hx),D\gf(\hx),D^2\gf(\hx))\le 0\eq$$ 
 \re\eqac 
\noindent whenever
$\gf\ti{is} C^2 $ and $\hx$ is a local maximum of $ 
u-\gf$. Before
making any formal definitions, let us also note that 
$u(x)\le
u(\hx)-\gf(\hx)+\gf(x)$ for $x$ near $\hx,  \gf\in C^2$ and
Taylor approximation imply  
$$u(x)\le u(\hx)+\ip{p,x-\hx}+\tfrac 12
\ip{X(x-\hx),x-\hx}+o(|x-\hx|^2)\iq{as}x\rtra\hx\eq$$ 
\re\eqab 
\noindent where
$p=D\gf(\hx)$ and $X=D^2\gf(\hx)$.  Moreover, if \eqab) 
holds for
some $(p,X)\in\RN\times\m$ and $u$ is twice differentiable 
at
$\hx$, then $p=Du(\hx)$ and $D^2u(\hx)\le X$.  Thus if $u$ 
is a
classical solution of $F\le 0$ it follows that
$F(\hx,u(\hx),p,X)\le 0$ whenever \eqab) holds; we may also
consider basing a definition of nondifferentiable 
solutions $u$ of
$F\le 0$ on this fact.    Roughly speaking, pursuing 
\eqac) leads
to notions based upon test functions $\gf$ but does not 
immediately
lead us, as will pursuing \eqab), to define ``$(Du,D^2u)$" 
for
nondifferentiable  functions $u$, which will turn out to 
be a good
idea.  For this reason, we begin by developing  the line 
suggested
by \eqab).  Next we introduce a set $\co\subset\RN$ on which
$F\le0$ is to hold and the appropriate notation to deal with
inequalities like \eqab) ``on $\co$."  At the moment, 
$\co$ is
arbitrary; later  we  require it to be locally compact.  
Taking off
from \eqab), if $u\:\co\rtra\R$, $\hx\in\co$, and \eqab) 
holds as $
\co\ni x\rtra\hx$, we say $(p,X)\in\jtpco u(\hx)$ (the 
second-order
``superjet" of $u$ at $\hx$).   This defines a mapping 
$\jtpco u$
from $\co$ to the subsets of $\rsn$. 

\rpre\exax 
\ex{Example \ \li}    By way of
illustration, if $u$ is defined on $\R$ by  
$$u(x)=\cases 0&\text{for $x\le 0$,}\cr ax+{b\over2}x^2
&\text{for $x\ge 0$,}\endcases$$
 then
$\jtp_{[-1,0]}u(0)=((-\infty,0)\times\R)\cup(\{0\}%
\times[0,\infty))$,
while 
$$\jtp_{\R}u(0)=\cases \emptyset&\text{if 
$a>0,$}\cr\{0\}\times[\max\{0,b\},\infty)&\text{if 
$a=0,$}\cr
((a,0)\times\R)\cup(\{0\}\times[0,\infty))\cup(\{a\}\times 
[b,\infty))&\text{if  $a<0.$}\endcases$$ 
\endex

Having thought through this example, the reader will see 
that
$\jtpco u(x)$ depends on $\co$ but realize it is the same 
for all
sets $\co$ for which $x$ is an interior point; we let 
$\jtp u(x)$
denote this common value.  If we repeat the above 
discussion after
switching the inequality sign in \eqab), we arrive at the
definitions of the second-order ``subjets" $\jtmco u$, 
$\jtm u$;
equivalently, $\jtmco u(x)=-\jtpco(- u)(x)$, etc. 

We are ready to define the notions of viscosity
subsolutions, supersolutions, and solutions. It will be 
useful to
have the notations 
$$\eqalign{&\roman{USC}(\co)=\{\ti{upper semicontinuous
functions}u\:\co\rtra\R\},\cr &
\roman{LSC}(\co)=\{\ti{lower semicontinuous 
functions}u\:\co\rtra\R\}.}$$

\rpre\dea 
\dfn{Definition \li} Let $F$ satisfy \eqaa) and
$\co\subset\RN$.  A {\it viscosity subsolution} of $F=0$
(equivalently, a viscosity solution of $F\le 0$) on $\co$ 
is a 
function $u\in \roman{USC}(\co)$ such that 
 $$F(x,u(x),p,X)\le0\quad\text{for
all}\quad x\in\co\quad\text{and}\quad(p,X)\in\jtpco 
u(x).\eq$$
 \re\eqad  
\noindent Similarly,
a {\it viscosity supersolution} of $F=0$ on $\co$ is a   
function
$u\in  LSC(\co)$ such that 
$$F(x,u(x),p,X)\ge0\quad\text{for
all}\quad x\in\co\quad\text{and}\quad(p,X)\in\jtmco 
u(x).\eq$$
 \re\eqae 
Finally, $u$
is a {\it viscosity  solution} of $F=0$ in $\co$ if it is 
both a
viscosity subsolution and a viscosity supersolution of 
$F=0$ in
$\co$.\enddfn

\rpre\rea 
\rem{Remarks\  {\rm \li}}  Since these ``viscosity
notions" are the primary ones for the current discussion, we
immediately agree (at least, we hope you agree) to drop 
the term
``viscosity" and hereafter simply refer to subsolutions,
supersolutions, and solutions.  This is a happy idea, as 
the term
``viscosity," which lacks elegance, is an artifact of the 
origin of
this theory in the study of first-order equations and the 
name was
then motivated by the consistency of the notion with the 
method of
``vanishing viscosity," which is irrelevant for many 
second-order
equations.   It follows from the discussion preceding the
definition that, for example, if $u$ is a solution of 
$F\le 0$ in
$\co$, $\gf$ is $C^2$ in a  neighborhood of $\co$, and 
$u-\gf$ has a
local maximum (relative to $\co$) at $\hx\in\co$, then  
\eqac)
holds.  Analogous remarks hold for supersolutions.  These 
remarks
motivate  the requirement that a subsolution be
upper semicontinuous, etc., in the sense that producing 
maxima of
upper semicontinous functions is straightforward.  
Solutions,
being both upper semicontinuous and lower semicontinuous, 
are
continuous.  One might ask  if the validity of \eqac) for 
all
$\gf\in C^2$ (with the maxima relative to $\co$) for an
upper semicontinuous function $u$ is equivalent to  $u$ 
being a
subsolution.  This is so. In fact, if $\hx\in\co$ then
  $$\jtpco u(\hx)=\{(D\gf(\hx),D^2\gf(\hx)):\gf\ti{is}C^2%
\ti{and}u-\gf\ti{has a
local maximum at}\hx\};$$ 
we leave the proof as an interesting
exercise.
\endrem

We next  record the definitions of the  closures of the 
set-valued
mappings  needed in the next section. 

With the above notation,  for $x\in\co$, we set 
\pre\eqzzp
$$\eqalign{\jtpcoc u(x)=\{&(p,X)\in\rsn:
\exists(x_n,p_n,X_n)\in\co\times\rsn\ni\cr & 
(p_n,X_n)\in\jtpco
u(x_n)\ti{and}(x_n,u(x_n),p_n,X_n)\!\rtra\!(x,u(x)\!,p,X)%
\}\cr}\eq$$
\pre\eqzzm
 and  
$$\eqalign{\jtmcoc u(x)=\{&(p,X)\in\rsn:
\exists(x_n,p_n,X_n)\in\co\times\rsn\ni\cr& 
(p_n,X_n)\in\jtmco
u(x_n)\ti{and}(x_n,u(x_n)\!,p_n,X_n)\!\rtra\!(x,u(x)%
\!,p,X)\}\!;\cr}\eq$$
 we are abusing standard practice as regards defining 
closures of
set-valued mappings a bit in that we put the extra condition
$u(x_n)\rtra u(x)$ in the definitions while the graphs of 
the
multifunctions $\jtpco u,\ \jtmco u$ do not themselves 
record the
values of $u$.  The reader may note the use of expressions 
like
``$x_n\rtra x$" as an abbreviation for ``the sequence 
$x_n$ satisfies
$\lim_{n\rtra\infty}x_n=x$," etc.   If 
$x\in\text{interior}(\co)$, we
define $\jtpc u(x)$, $\jtmc u(x)$ in the obvious way.

\rpre\reb 
\rem{Remark\ \ {\rm \li}}  If $u$ is a solution of 
$F\le0$ on $\co$, then $F(x,u(x),p,X)\le 0$ for $x\in\co$ 
and
$(p,X)\in \jtpco u(x)$.  This remains true, for reasons of 
continuity,  if $(p,X)\in\jtpcoc u(x)$ and $F$ is continuous
(or even lower semicontinuous).  Similar remarks apply to
supersolutions and solutions.
\endrem

\rem{Advice}  We advise the reader to either skip the
following material or to  scan it lightly at the present 
time and
proceed directly to the next section. The comments 
collected below
can be referred to as needed.\endrem

\rpre\reax 
\rem{Remark \ {\rm \li}}  While the definitions
above may seem reasonable, they contain subtleties.  In 
particular, 
they do not define ``operators" on domains in a familiar 
way.
 For example, let $N=1$ and consider the equation $(u')^2
-1 = 0$.  The functions $|x|$ and $-|x|$ both satisfy this 
equation
classically on $\R\backslash\{0\}$.  However, the semijets 
of the
function $u(x)=-|x|$  are given by   
$$\jtp u(x)=\cases \{(-1,X):
X\ge 0\}&\text{if $ 
x>0,$}\cr((-1,1)\times\R)\cup(\{-1,1\}\times[0,\infty))&
\text{if $ x=0$},\cr \{(1,X): X\ge 0\}&\text{if 
 $x<0$}\endcases$$
 and 
$$\jtm u(x)=\cases \{(-1,X): X\le 0\}&\text{if $
x>0,$}\cr\emptyset&\text{if $ x=0$},\cr \{(1,X): X\le 
0\}&\text{if 
 $x<0$,}\endcases$$ 
so $-|x|$ is a solution of the equation, while if
$v(x)=|x|$ then $(0,0)\in\jtm v(0)$ and $v$ is not a 
solution. 
Similarly, one sees that $v$ does solve $-(v')^2+1=0$.  
Hence a
solution of $F=0$ need not be a solution of $-F=0$ (of 
course, if
$-F$ is not \good, as is typical when  $F$ depends on $u, 
X$, this
is  meaningless).   It is also unknown if  the information 
 $H, u,
f, g\in C(\RN)$ and $u$ is a solution of $H(Du)-f(x)=0$ as 
well as
$H(Du)-g(x)=0$ in $\RN$ implies $f=g$ (although it is 
known that
$f=g$ if either $H$ is uniformly continuous or $N=1$).
\endrem

We do not mention these things to disturb the reader; 
indeed, they 
will not explicitly appear in the theory.  Our purpose is 
to emphasize 
that the theory to follow is not a variant of more 
classical developments.

\rpre\reaxx 
\rem{Remark \ {\rm \li}}  If we consider one of
the simplest examples of a discontinuous  upper 
semicontinuous
function, namely,  the function given by  $u(x)=0$ for 
$x\not=0$,
$u(0)=1$, we have $\jtp u(0)=\R\times\R$.  The relation
$F(0,1,p,X)\le 0$ for $(p,X)\in\jtp u(0)$ is then very
restrictive.   In fact,   a primary role played by allowing
semicontinuous subsolutions and supersolutons in the basic 
theory
is only as a technical device to produce continuous 
solutions via
Perron's method as described in \S \perron.  If solutions 
are 
produced by other methods (e.g., via numerical or other
approximations)  there may be no need to invoke results on
semicontinuous functions.  However, one of the simplest 
ways to
show uniform convergence of approximation schemes,
 involves comparison of semicontinuous subsolutions and
supersolutions.  Thus one often deals with semicontinuity 
mainly 
as a device in proofs---as such, it is a powerful labor 
saving
device. 
\endrem

\rpre\rebx 
\rem{Remarks\ {\rm \li}}  We point out some facts
concerning second-order semijets.   

(i)  It is clear that $\jtpco u(x)$ is a convex
subset of $\RN\times\m$;  Example \exax\  shows that it is 
not
necessarily closed.   However,  a little thought shows 
that if
$p\in\RN$, then $\{X:(p,X)\in\jtpco u(x)\}$ is closed. 

(ii)  If $\gf\in C^2$ in some neighborhood of
$\co$, then $\jtpco (u-\gf)(x)=\{(p-D\gf(x), 
X-D^2\gf(x)):(p,X)\in
\jtpco u(x)\}$.   As a consequence, the same statement 
holds if
$\jtpco$ is replaced everywhere by $\jtmco,\ \jtpcoc,$ or $
\jtmcoc$.  Indeed, if $(q,Y)\in\jtpco (u-\gf)(\hx)$, then
 $$u(x)-\gf(x)\le u(\hx)-\gf(\hx)+\ip{q,x-\hx}+\tfrac 1
2\ip{Y(x-\hx),x-\hx}+o(|x-\hx|^2)\quad\ti{as}x\rtra\hx\eq$$ 
\re\eqabx 
\noindent and 
$$\gf(x)-\gf(\hx)=\ip{D\gf(\hx),(x-\hx)}+\tfrac 
12\ip{D^2\gf(\hx)(x-\hx),x-\hx}+
o(|x-\hx|^2)$$ 
imply that $(q+D\gf(\hx),Y+D^2\gf(\hx))\in\jtpco
u(\hx)$ and so $\jtpco (u-\gf)(\hx)\subset\{(p-D\gf(\hx),
X-D^2\gf(x)):(p,X)\in \jtpco u(x)\}$.   The other 
inclusion follows
from this as well, since $\jtpco u(\hx)= \jtpco ((u-\gf)+
\gf)(\hx)$.
It is also clear that one always has 
$$\jtpco (u+v)(x)\supset\jtpco u(x)+\jtpco v(x).$$

(iii)  We consider $\jtpco \gf(\hx)$ when
$\gf\in C^2(\RN)$ and will end up with a general statement
corresponding to  Example \exax.   In view of (ii), we may 
as well
assume that $\gf\equiv 0$, and we will write ``Zero" for 
the zero
function.  We know that if $\hx\in{\text{interior}}(\co)$, 
then
$\jtpco {\text{Zero}}(\hx)=\jtp 
{\text{Zero}}(\hx)=\{(0,X)\}: X\ge
0\}$.   In general,
 $(p,X)\in\jtpco {\text{Zero}}(\hx)$ if 
$$\eqalign{0\le\ip{p,x-\hx}+\tfrac 1
2\ip{X(x-\hx),x-\hx}+o(|x-\hx|^2)\quad\ti{as}\co\ni 
x\rtra\hx.\cr}\eq$$
\re\eqzb 
\noindent Assuming that $x_n\in\co, \ 0<|x_n-\hx|\rtra 0$,
and $(x_n-\hx)/|x_n-\hx|\rtra q$, we may put $x=x_n$ in 
\eqzb),
divide the result by $|x_n-\hx|$ and pass to the limit to 
conclude
that  
$$0\le\ip{p,q}\iq{for}q\in \roman{UT}_\co(\hx)\eq$$ 
\re\eqzc 
\noindent where
$\roman{UT}_\co(\hx)$ is the set of ``generalized unit 
tangents" to $\co$
at $\hx$; it is given by 
$$\roman{UT}_\co(\hx)=\left\{q:\exists
x_n\in\co\backslash \{\hx\},\
x_n\rtra\hx,\iq{and}{x_n-\hx\over|x_n-\hx|}\rtra 
q\right\}.\eq$$ 
\re\eqzd  
\noindent If $\co$ is a smooth $N$-submanifold  of $\RN$ 
with boundary and
$\hx\in\partial\co$, then the generalized tangent cone 
$$T_{\co}(\hx)=\text{
convex\ hull}(\roman{UT}_\co(\hx))$$ 
is a halfspace and $\co$ has an
exterior normal $\vec n$ at $\hx$.  In this event, \eqzc) 
says that
$p=-\gl \vec n$ for some $\gl\ge0$.  Moreover, if  $p=0$, 
\eqzb)
then implies that $0\le \ip{Xq,q}$ for $q\in 
\roman{UT}_{\co}(\hx)$, and we
conclude that   
$$(0,X)\in\jtpco {\text{Zero}}(\hx)\iq{if and only
if}0\le X\eq$$  
\re\eqze  
provided $T_\co(\hx)$ is a halfspace.
\endrem

Life is more complex if $p=-\gl\vec n$ and $\gl>0$.  In 
this case,
\eqzb) reads 
$$\gl\ip{\vec n,x-\hx}-\tfrac 12\ip{X(x-\hx),x-\hx}\le
o(|x-\hx|^2) \quad\ti{as}\co\ni x\rtra\hx.\eq$$ 
\re\eqzf
\noindent We study \eqzf) when $\hx=0$ and we can 
represent $\co$ near $0$
in the form
$$\left\{(\tilde x,x_N):x_N\le 
g(\tx)\right\}\eq$$ 
\re\eqzg 
\noindent where
$\tx=(x_1,\ldots,x_{N-1})$, $g(\tx)=\oo 2\ip{Z\tx,\tx}+
o(|\tx|^2)\}$,
and $Z\in{{\cal S}(N-1)}$.  That is, we
assume that the boundary of $\co$ is twice differentiable 
at $0$
and rotate so that the normal is $\vec 
n=e_N=(0,\ldots,0,1)$.  With
these normalizations, we put $(\tx,g(\tx))$ into \eqzf) to 
find $\gl
\ip{Z\tx,\tx}-\ip{X(\tx,0),(\tx,0)}\le o(\nos{\tx})$ or $ 
\gl \widetilde Z\le
P_{N-1}XP_{N-1}$ where 
$\widetilde Z(\tilde x,x_n)=(Z\tilde x,0)$ and
$P_{N-1}$ is the  projection on the first
$N-1$ coordinates.  It is not hard to  see that this is also
sufficient.  

When unraveled, the above considerations lead to the 
following
conclusion. Let $\co$ be an $N$-submanifold of $\RN$ with
boundary, $\hx\in\partial\co$, $\partial \co$ be twice
differentiable at $\hx$, $\vec n$ be the outward normal,
$T_{\partial\co}(\hx)$ be the tangent plane to 
$\partial\co$ at
$\hx$, and $P\in\m$  be orthogonal projection on
$T_{\partial\co}(\hx)$ (which we regard as a subspace of 
$\RN$). 
Finally, let $S$ be the symmetric operator in
$T_{\partial\co}(\hx)$ corresponding to the second 
fundamental form
of $\partial\co$ at $\hx$ (oriented with the exterior 
normal to
$\co$) extended to $\Bbb R^N$ by $S\vec n=0$. Then
 $$
\cases (p,X)\in \jtpco \gf(\hx) & \text{if and only if 
either } 
p=D\gf(\hx)\ti{and} D^2\gf(\hx)\le X,\text{ or}\\
p=D\gf(\hx)-\gl \vec n, & \gl>0 \ti{and}PD^2\gf(\hx)P\le
PXP-\gl S.\endcases\eq$$ 
\re\eqzfxx 
Noting that $P\vec n\otimes\vec
n P=0$, we see that if  $(p,X)\in\jtpco u(\hx)$, $S\le 0$, 
and $\gl>0$, 
then  $(p-\gl \vec n,X+\gm \vec n\otimes\vec n)\in\jtpco 
u(\hx)$  for $\gm\in\R$.  For
later use, we denote by $X_\gl$ the particular choice of 
$X$ such that $(D\gf(\hx)-\gl
\vec n, X)\in\jtpco \gf(\hx)$ given by 
$$X_\gl(\xi+\gb \vec
n)=D^2\gf(\hx)(\xi+\gb\vec n) +\gl S\xi\iq{for} \xi\in
T_{\partial\co}(\hx), \gb\in\R.\eq$$ 
\re\eqzfx

Before leaving this topic, we note that the set of exterior
``normals" 
$$N(\hx)=(\vec n\in\RN:\ip{\vec n,q}\le
0\iq{for}q\in T_{\co}(\hx))\eq$$ 
\re\eqzh 
\noindent can also be described by
$$N(\hx)=\paren{\vec n\in\RN:\ip{\vec n,x-\hx}\le
o(|x-\hx|)\ti{as}\co\ni x\rtra \hx}\eq$$ 
\re\eqzi 
\noindent no matter how
regular $\partial\co$ is at $\hx$ and that we showed above 
that
$$\jtpco \gf(\hx)\subset
\left(D\gf(\hx)-N(\hx)\right)\times\m.\eq$$ 
\re\eqzj

\subheading{Notes on \S 2}
Viscosity solutions were introduced by M. G. Crandall and 
P. L.
Lions [\clb]  (following an announcement in [\cla]);  
analogies
with S. N.  Kruzkov's theory of  scalar conservation laws 
[\kru]
provided guidance for the notion  and its  presentation. The
presentation of [\clb] emphasized the first-order case,  
owing to
the fact  that uniqueness  results were only obtainable 
for this
case at the time.   It was also pointed out in [\clb] that 
there
were several  equivalent ways to formulate  the notion of 
viscosity
solutions; one of these  was intimately connected with 
earlier work
by L. C. Evans [\evco, \evac] concerning ``weak passages  
to the
limit" in equations satisfying the maximum principle. In 
these
developments, directions were  indicated by aspects of  
nonlinear
functional analysis   and nonlinear  semigroup  theory.

It became apparent that working  with these  alternative 
formulations was superior to the approach taken in [\cla], 
and
perhaps  M. G. Crandall, L. C. Evans, and P. L. Lions 
[\cel] is the
first readable account of the early theory.  (We note that 
the 
definition most emphasized in [\cel] is not the one put 
foremost 
here.)

The necessity of defining the semijets on closed sets became
apparent  with the investigation of boundary value 
problems.  We
return to  this issue in \S \boundary, but a  couple of
significant papers here are  P. L. Lions [\pllneumann] and 
H. M.
Soner [\sona]  (see also
M. G. Crandall and R. Newcomb [\cn]). The
closure of the semijets was introduced in  H. Ishii and P. 
L. Lions
[\lii] and was fully exploited in M. G. Crandall  [\cla] 
and M. G.
Crandall and H. Ishii [\craish].   Analogous operations  
have long
been employed by  ``nonsmooth analysts."

The fact that the theory required only appropriate 
semicontinuity
was  realized earlier on, but the first striking 
applications of
the fact  (to existence, stability, optimal control, 
$\ldots)$
were given in H. Ishii [\isp, \isbvp], G. Barles and B. 
Perthame
[\bapa, \bapb].

Let us finally mention that some structural properties of 
viscosity 
solutions were and will be ommitted from our presentation 
and, in 
addition to the above references, some can be found in P. 
L. Lions
[\pla],  R. Jensen and P. E. Souganidis [\jensoug], and H.
Frankowska [\fra]. \medskip

\resultno=0 \equationno=0 \spre\maximum 
\heading \sli. The maximum principle for
semicontinuous functions \\
    and comparison for the Dirichlet problem\endheading

\def\secno{\the\sectionno}  
Let $\gO$ be a bounded open subset of
$\RN$, $\goc$ be its closure, and $\partial\gO$ be its 
boundary. 
Suppose $u\in\roman{USC}(\overline\Omega)$ 
is a solution of $F\le 0$ in $\gO$, $v\in
\roman{LSC}(\overline\Omega)$
is a solution of $F\ge0$, and $u\le v$ on $\partial\gO$.  
We seek to
show that $u\le v$ on $\goc$.  

In the event $u$ and $v$ are classical sub and 
supersolutions, we
could employ the classical ``maximum principle."   Let us 
recall
some elementary facts in this regard.   Let  $w$ be 
defined in a
neighborhood of $\hx\in\RN$.  If there exists  
$(p,X)\in\rsn$ such
that
$$w(x)=w(\hx)+\ip{p,x-\hx}+\tfrac12 
\ip{X(x-\hx),x-\hx}+o(|x-\hx|^2)\iq{as}x\rtra\hx,\eq$$
\re\eqba 
\noindent we say that $w$ is twice differentiable at $\hx$ 
and
$Dw(\hx)=p$, $D^2w(\hx)=X$ (it being obvious that $p$ and 
$X$ are
unique if they exist).  It is clear that $w$ is twice
differentiable at $\hx$ if and only if $\jt w(\hx)\equiv\jtp
w(\hx)\cap\jtm w(\hx)$ is nonempty (in which case $\jt
w(\hx)=\{(Dw(\hx),D^2w(\hx))\}$).  Classical 
implementations of the
maximum principle are based on the fact, already used 
above, that
if $w$ is twice differentiable at a local maximum $\hx$, 
then
$Dw(\hx)=0$ and $D^2w(\hx)\le 0$.  Thus, if $u$ and $v$ 
are twice
differentiable everywhere and $w=u-v$ has a local maximum
$\hx\in\gO$, we would have $Du(\hx)=Dv(\hx)$ and 
$D^2u(\hx)\le
D^2v(\hx)$ and then, in view of \eqaa),
$$
\align
F(\hx,u(\hx),Du(\hx),D^2u(\hx)) & \le0\le
F(\hx,v(\hx),Dv(\hx),D^2v(\hx))\\
&  \le
F(\hx,v(\hx),Du(\hx),D^2u(\hx)).\endalign$$ 
In the event that $F(x,r,p,X)$
is strictly nondecreasing in $r$ (a simple but 
illustrative case),
it follows that $u-v$ is nonpositive at an interior 
maximum and so
$u\le v$ in $\gO$ since $u\le v$ holds on $\bgo$.

We seek to extend this argument to the case $u\in\uscc,\
v\in\lscc$.  We are unable to simply plug 
$(Du(\hx),D^2u(\hx))$ and
$(Dv(\hx),D^2v(\hx))$ into $F$ since these expressions 
must be
replaced by the set-valued functions $\jtp u$ and $\jtm v$ 
(and
their values may well be empty at many points, including 
maximum
points of $u(x)-v(x)$). To use $\jtp u$ and $\jtm v$, we 
employ a 
device that doubles the number of variables and then 
penalizes
this doubling.  More precisely, we  maximize the function
$u(x)-v(y)-(\ga/2)\nos{x-y}$ over $\goc\times\goc$; here 
$\ga>0$ is
a parameter.   As $\ga\rtra\infty$, we  closely approximate
maximizing $u(x)-v(x)$ over $\goc$.  More precisely, we have
\rpre\lea \myProclaim Lemma \li. Let $\co$ be a   
subset of $\RN$, $u\in \roman{USC}(\co),\ 
v\in\roman{LSC}(\co)$ and  
$$M_\ga
=\sup_{\co\times\co}(u(x)-v(y)-\tfrac\ga 2\nos{x-y}) \eq$$  
\re\eqbb
for $\ga>0$.  Let $M_\ga<\infty$ for large $\ga$ and
$(x_\ga,y_\ga)$ be such that 
$$\lim_{\ga\rtra\infty}(M_\ga -
(u(x_\ga)-v(y_\ga)-\tfrac\ga 2\nos{x_\ga-y_\ga}))=0.\eq$$ 
\re\eqbc 
\noindent Then
the following holds{\rm :} 
$$\left\{\eqalign{&{\roman{(i)}\ 
}\lim\nolimits_{\ga\rtra\infty}\ga\nos{x_{\ga}-y_{\ga}}=0%
\ti{and}
\cr&\eqalign{{\roman{(ii)}\  }
&\lim\nolimits_{\ga\rtra\infty}M_\ga=u(\hx)-v(\hx)=
\sup\nolimits_{\co}(u(x)-v(x))\cr&\ti{whenever}
\hx\in\co\ti{is a limit point
of}x_\ga\ti{as}\ga\rtra\infty.\cr}\cr}\right.\eq$$ 
\re\eqbd

  Deferring the elementary proof of the lemma and 
returning to our sub and supersolutions
$u\in\uscc,\ v\in\lscc$ satisfying $u\le v$ on $\bgo$, we 
note that
$M_\ga=\sup_{\gocs}(u(x)-v(y)-(\ga/2)\nos{x-y})$ is finite 
since
$u(x)-v(y)$ is upper semicontinuous and $\goc$ is compact. 
 Since
we seek to prove $u\le v$, we assume to the contrary that
$u(z)>v(z)$ for some $z\in\gO$; it follows that  
$$M_\ga\ge
u(z)-v(z)=\gd>0\iq{for}\ga>0.\eq$$ 
\re\eqbe 
\noindent Choosing
$(x\sga,y\sga)$ so that
$M\sga=u(x\sga)-v(y\sga)-(\ga/2)\nos{x\sga-y\sga}$ (the 
maximum is
achieved in view of upper semicontinuity and compactness), 
it
follows from \eqbd)(i), (ii) and $u\le v$ on $\bgo$ that
$(x\sga,y\sga)\in\gos$ for  $\ga$ large.  The next step is 
to use
the equations to estimate $M_\ga$ and contradict \eqbe) 
for large
$\ga$.   This requires producing suitable values of $\jtp 
u$ and
$\jtm v$, and we turn to this question.

To know what to look for, let us proceed more generally 
and assume
that $u$, $v$ are defined in  neighborhoods of 
$\hx,\hy\in\RN$ and
twice differentiable at $\hx,\ \hy$ respectively.  Assume,
moreover, that $\gf$ is $C^2$ near $\ph$ in $\RN\times\RN$ 
and
$\ph$ is a local maximum of $u(x)-v(y)-\gf(x,y)$.  
Applying the
classical maximum principle to this situation (in the 2N 
variables
$\p$), we learn that $Du(\hx)=D_x\gf(\hx,\hy)$, 
$Dv(\hy)=-D_y\gf(\hx,\hy)$, and  
$$\pmatrix X&\ 0\\0&-Y\endpmatrix\le
D^2\gf\ph\eq$$ 
\re\eqbf 
\noindent where $X=D^2u(\hx),\ Y=D^2v(\hy)$. Note
that with the choice $\gf\p=(\ga/2)\xmys$ the above reads
$$\pmatrix X&\ 0\cr0&-Y\endpmatrix\le\ga\pmatrix \ 
I&-I\cr-I&\ I\endpmatrix\eq$$
\re\eqbg 
\noindent where $I$ will stand for the identity matrix in 
any
dimension and, since the right-hand side annihilates  
vectors of
the form $(\smallmatrix \xi\\ \xi\endsmallmatrix)$ 
(also written $\ut(\xi,\xi)$ where
$\ut Z$ denotes the transpose of a matrix $Z$), \eqbg) 
implies $X\le
Y$, making further contact with the maximum principle.

It is a remarkable fact that perturbations of the above 
results may
be obtained in the class of semicontinuous functions.   
The main
result we use in this direction is the following theorem, 
which is
formulated in a useful but distracting generality.  (See 
\eqzzp)
and \eqzzm) regarding  notation below.)     
 \rpre\tha  \myProclaim Theorem \li. Let $\co_i$ be a 
locally compact
subset of $R^{N_i}$ for  $i=1,\ldots,k$, 
$$\co=\co_1\times\cdots\times\co_k,$$ 
$u_i\in \roman{USC}(\co_i),$ and
$\gf$  be twice continuously differentiable in a 
neighborhood of
$\co$.    Set   
  $$w(x)=u_1(x_1)+\cdots+u_k(x_k)\quad\text{for 
}x=(x_1,\cdots,x_ k)\in\co,$$
and suppose $\hx=(\hx_1,\dotsc,\hx_k)\in\co$ is a local 
maximum of
$w-\gf$ relative to $\co$. Then for each $\gep>0$ there 
exists
$X_i\in\md{N_i}$  such that  
$$
(D_{x_i}\gf(\hx),X_i)\in\overline J^{2,+}_{\scr 
O_{i}}u_i(\hx_i)  \quad\text{for }i=1,\ldots,k,
$$ 
and the block diagonal matrix
with entries $X_i$ satisfies  
$$-\left(\oo\gep +\vno A\right)I\le 
\diag\le A+\gep A^2\eq$$ 
where $A=D^2\gf(\hx)\in\m$,
$N=N_1+\cdots+N_k$. 
\re\eqbh

The norm of the symmetric matrix $A$ used in \eqbh) is
$$\vno{A}=\sup\{|\gl|: \gl \ti{is an eigenvalue
of}A\}=\sup\{|\ip{A\xi,\xi}|:|\xi|\le 1\}.$$ 
We caution the reader
that this result is an efficient summarization of the  
analytical
heart of the theory suitable for the presentation we have 
chosen
and its   proof, which is outlined in the appendix, is 
deeper and
more difficult than  its applications that we give in the
main text.  See also the notes to this  section.

In order  to apply  Theorem \tha\ in the above situation,  
we put
$k=2,$ 
 $\co_1=\co_2=\gO$, $u_1=u$, $u_2=-v$, 
$\gf(x,y)=(\ga/2)\nos{x-y}$,
and recall that $\jtmoc v=-\jtpoc (-v)$.    In this case  
$$
\gather
D_x\gf(\hx,\hy)=-D_y\gf(\hx,\hy)=\ga(\hx-\hy), \qquad
A=\ga\pmatrix\ I&-I\\ -I&\ I\endpmatrix,\\
A^2=2\ga A,\qquad \vno A=2\ga,\endgather$$ 
and  we conclude that for every $\gep>0$ there exists $X, 
\ Y\in\m$
such that 
$$(\ga(\hx-\hy),X)\in\jtpoc u(\hx),\qquad
(\ga(\hx-\hy),Y)\in\jtmoc v(\hy)\eq$$ 
\re\eqbhx 
and
$$-\left(\oo\gep+2\ga\right)\pmatrix I&0\cr 0& 
I\endpmatrix\le
\pmatrix X&0\cr 0&-Y\endpmatrix\le  
\ga(1+2\gep\ga)\pmatrix\ I&-I\cr -I&\ I\endpmatrix.$$ 
Choosing 
$\gep=1/\ga$ yields the  elegant relations 
$$-3\ga\pmatrix I&0\cr 0& I\endpmatrix\le
\pmatrix X&\ 0\cr 0&-Y\endpmatrix\le 3\ga \pmatrix \ 
I&-I\cr -I&\ I\endpmatrix;\eq$$ 
\re\eqbi 
\noindent in all, we
conclude that at a maximum $\ph$ of 
$u(x)-v(y)-(\ga/2)\xmys$ there exists $X,
-Y\in\m$ such that \eqbhx) and \eqbi) hold.  Observe that, 
in
comparison with \eqbg), we  have worsened the upper  bound 
by a
factor of 3 (but can still conclude from \eqbi) that $X\le 
Y$) and
obtained a lower bound.  We will not actually use the 
lower bound
in what follows, but its presence corresponds to an 
essential
compactness issue in the proof of the theorem.

Let us finish this chapter of the story.  Writing $\ph$ in 
place of
$(x\sga,y\sga)$ for simplicity,  we only need to assume 
that 
$$F(\hx,u(\hx),\ga(\hx-\hy),X)\le 0\le
F(\hy,v(\hy),\ga(\hx-\hy),Y)\eq$$ 
\re\eqbj
\noindent is inconsistent with the other information we 
have in hand, which
includes  
$$\left\{\eqalign{&0<\gd\le
M_\ga=u(\hx)-v(\hy)-(\ga/2)|\hx-\hy|^2\le u(\hx)-v(\hy) 
\text{ and}\cr
&\ga|\hx-\hy|^2\rtra 0\quad\text{as 
}\ga\rtra\infty\cr}\right.\eq$$
\re\eqbjx 
\noindent
(the second assertion follows from \eqbd)(i)) and
\eqbi).  For example, assume that there exists $\gamma>0$ 
such that
\pre\eqbk  
$$\gamma(r-s)\le F(x,\!r,\!p,\!X)-F(x,\!s,\!p,\!X)\iq{for} 
r\ge s,\
(x,\!p,\!X)\in\goc\times\RN\times\m\eq$$ 
and there is a
function $\go\:[0,\infty]\rtra[0,\infty]$ that satisfies
$\go(0+)=0$ such that
\pre\eqbl
$$
\cases
F(y,r,\ga(x-y),Y)-F(x,r,\ga(x-y),X)\le
\go(\ga\xmys+|x-y|)\\
\text{whenever } x,y\in\gO,\ r\in\R,\
X,Y\in\m,\ti{and} \eqbi) \text{ holds}.\endcases\eq$$   

Proceding,  we deduce from  \eqbjx), \eqbk), \eqbj), and 
\eqbl)
that 
\pre\eqblx
$$
\aligned
\gamma\gd & \le\gamma(u(\hx)-v(\hy))\le
F(\hx,u(\hx),\ga(\hx-\hy),X)- F(\hx,v(\hy),\ga(\hx-\hy),X)\\
&=
(F(\hx,u(\hx),\ga(\hx-\hy),X)- 
F(\hy,v(\hy),\ga(\hx-\hy),Y))\\
&\hphantom{=}+
(F(\hy,v(\hy),\ga(\hx-\hy),Y)-
F(\hx,v(\hy),\ga(\hx-\hy),X))\\
&\le \go(\ga|\hx-\hy|^2+|\hx -\hy|);\endaligned \eq$$  
here we
used \eqbj) to estimate the first term on the right by 0 
and 
\eqbl) on the second term.  Since $\go(\ga|\hx-\hy|^2+|\hx 
-\hy|)\rtra0$ as
$\ga\rtra\infty$ by \eqbjx), we have a contradiction.  We 
have proved

\rpre\thb \myProclaim Theorem \li.  Let $\gO$ be a bounded 
open
subset of $\RN$, $F\in C(\gO\times\R\times\RN\times\m)$ be
\good\  and satisfy {\rm \eqbk), \eqbl)}.  Let $u\in
\roman{USC}(\goc)$ {\rm (}respectively, $v\in 
\roman{LSC}(\goc))$ be a   subsolution
{\rm (}respectively, supersolution{\rm )} of $F=0$ in 
$\gO$ and $u\le v$ on
$\bgo$.  Then $u\le v$ in $\goc$.

\rpre\rec 
\rem{Remark\  {\rm \li}}  We motivated the 
structure condition \eqbl) by the natural way it arose in 
the
course of analysis and \eqbl) was used in the proof and not
degenerate ellipticity.  However, \eqbl) implies degenerate
ellipticity.  To see this,  suppose that \eqbl) holds, 
$X,Y\in\m$
and  $X\le Y$.  Observe that then for $\xi,\eta\in\RN$, and
$\gep>0$ 
$$\eqalign{\ip{X\xi,\xi}-\ip{Y\eta,\eta}&\le
\ip{Y\xi,\xi}-\ip{Y\eta,\eta}=
2\ip{Y\eta,\xi-\eta}+\ip{Y(\xi-\eta),\xi-\eta}\cr &
\le\gep\nos\eta+\left(1+{\vno Y\over \gep}\right)\vno 
Y\nos{\xi-\eta}.\cr}$$
 The last inequality may be written 
$$\pmatrix X&\ \ \
0\cr0&-(Y+\gep I)\endpmatrix\le\left(1+{\vno Y\over 
\gep}\right)\vno Y
\pmatrix \ I&-I\cr-I&\ I\endpmatrix.$$ 
Thus if 
$\ga>(1/3)\max\{\vno X,\vno Y,(1+\vno Y/\gep)\vno Y\}$ and 
$\gep$ is
sufficiently small,  the pair $X,Y+\gep I$ satisfy \eqbi). 
 Next we
fix $x, r$, and $p$, and put $y=x-p/\ga$ in \eqbl) so that 
by \eqbl)
$$F\left(x-{p\over\ga},r,p,Y+\gep I\right)-F(x,r,p,X)\le
\go\left(\frac 1\alpha(|p|^2+|p|)\right)$$ 
and then we let $\ga\rtra\infty$
and  $\gep\downa 0$ and conclude that $F$ is degenerate
elliptic.
\endrem

\rpre\recxx  
\rem{Remark\ {\rm \li}}   The proof of 
Theorem \thb\ adapts to provide modulus of continuity 
estimates on solutions of $F=0$.
\endrem

\rpre\recx 
\ex{Example\ \li} We turn to some examples
in which \eqbl) is satisfied. First notice  that \eqbl)  
evidently
holds (with $\go$ the modulus of continuity for $f$) 
for $F(x,r,p,X)=G(r,p,X)-f(x)$ if $G$ is
degenerate elliptic and $f$ is continuous; this is because 
\eqbi)
implies $X\le Y$.    

Secondly, the linear first order expression $\ip{b(x),p}$ 
satisfies \eqbl) if 
$$\ga\ip{b(y)-b(x),x-y}\le\go(\ga\nos{x-y}+|x-y|)$$ 
for some $\go$
and this holds with $\go(r)=cr$ if there is a constant 
$c>0$ such that
$\ip{b(x)-b(y),x-y}\ge-c\nos{x-y}$, i.e., the vector field 
$b+cI$
is ``monotone."  In fact, it is not hard to see that this 
is a
necessary condition as well.

Next, the linear expression  
$$G(x,X)=-\trace(\ut \gS(x)\gS(x)
X),\eq$$ 
\re\eqbm 
\noindent where $\gS$ maps $\goc$ into the $N\times N$ real
matrices, is degenerate elliptic  (it is a special case of 
Example
\extb) with $A(x)=\ut \gS(x)\gS(x))$. We seek to estimate
$G(y,Y)-G(x,X)$ when \eqbi) holds.  Multiplying the  
rightmost
inequality in \eqbi) by the nonnegative symmetric matrix 
$$\pmatrix \ut\gS(x)\gS(x)&\ut\gS(y)\gS(x))
\cr\ut\gS(x)\gS(y)&\ut\gS(y)\gS(y)\endpmatrix$$ 
and taking traces preserves the
inequality and yields 
$$\eqalign{G(y,Y)-G(x,X)&=\trace(\ut\gS(x)\gS(x)
X-\ut\gS(y)\gS(y)Y)\cr&\le3\ga\trace((\ut\gS(x)-\ut%
\gS(y))(\gS(x)-\gS(y)),\cr}\eq$$
\re\eqbn 
so if $\gS$ is Lipschitz continuous with constant $L$ then
$$G(y,Y)-G(x,X)\le3L^2\ga\xmys$$ 
and we may choose $\go(r)=3L^2r$.   
We note that if $\gG(x)\:\goc\rtra\m$ is Lipschitz
continuous and $\gG(x)>0$ in $\goc$, then
$$-\trace(\gG(x)X)=-\trace(\gG(x)^{1/2}\gG(x)^{1/2}X)$$
is an example
of the sort just discussed since $\gG(x)^{1/2}$ is also 
Lipschitz
continuous.  Finally, it is known that it $\gG(x)\ge 0$ and 
$\gG\in W^{2,\infty}$, then $\gG^{1/2}$ is Lipschitz 
continuous.

Now we note that the  processes of forming sums or sups or 
sup
infs, etc.,  above produce examples obeying \eqbl) if the
ingredients obey \eqbl) with  a common $\go$.  

Finally, we may produce examples satisfying \eqbk) and 
\eqbl) by
putting 
$$F(x,r,p,X)=\gamma r+G(x,r,p,X)$$
whenever $G$ is degenerate
elliptic and satisfies \eqbl).
\endex


  We deferred discussion of the elementary Lemma \lea
\ since it is better presented in a  generality not 
required above.
 Indeed, Lemma \lea\  may be obtained from  the following
proposition via the correspondences $2N\rtra M,\ 
\co\times\co\rtra\co,\ (x,y)\rtra x,\ 
u(x)-v(y)\rtra\gF(x), \
(1/2)|x-y|^2\rtra\Psi(x)$. \rpre\pra  \myProclaim 
Proposition \li. Let
$\co$ be a subset of $\R^M$, $\gF\in \roman{USC}(\co),\
\Psi\in \roman{LSC}(\co)$, $\Psi\ge0$, and  
$$M_\ga
=\sup_{\co}(\gF(x)-\ga\Psi(x)) \eq$$  
\re\eqbbx 
\noindent for $\ga>0$.  Let
$-\infty<\lim_{\ga\rtra\infty}M_{\ga}<\infty$ and 
$x_\ga\in\co$ be
chosen so  that  
$$\lim_{\ga\rtra\infty}(M_\ga -
(\gF(x_\ga)-\ga\Psi(x_\ga))=0.\eq$$ 
\re\eqbcx 
\noindent Then the following
hold{\rm :}  
$$\left\{\eqalign{&{\roman{(i)}\ 
}\lim\nolimits_{\ga\rtra\infty}\ga\Psi(x_\ga)=0, 
\cr&\eqalign{{\roman{(ii)}}\
&\Psi(\hx)=0\ti{and} 
\lim\nolimits_{\ga\rtra\infty}M_\ga=\gF(\hx)=
\sup\nolimits_{\{\Psi(x)=0\}}\gF(x)\cr&\qquad\qquad\quad%
\ti{whenever}
\hx\in\co\ti{is a limit point
of}x_\ga\ti{as}\ga\rtra\infty.\cr}\cr}\right.\eq$$ 
\re\eqbdx

\demo{Proof} Put 
$$\gd\sga=M_\ga - (\gF(x_\ga)-\ga\Psi(x_\ga))$$ 
so that
$\gd\sga\rtra0$ as $\ga\rtra\infty$.   Since  $\Psi\ge 0$, 
$M_\ga$
decreases as $\ga$ increases and 
$\lim_{\ga\rtra\infty}M\sga$ exists;
it is finite by assumption.   Moreover,  
$$ M_{\ga/2}\ge
\gF(x_\ga)-\ot\ga\Psi(x_\ga)\ge\gF(x_\ga)-\ga\Psi(x_\ga)
+\ot\ga\Psi(x_\ga)\ge M\sga -\gd\sga+\ot\ga\Psi(x\sga),$$ 
so
$2(M_{\ga/2}-M\sga+\gd_\ga)\ge\ga\Psi(x\sga)$, which shows 
that
 $\ga\Psi(x\sga)\rtra0$ as
$\ga\rtra\infty$.  Suppose now that $\ga_n\rtra\infty$ and
$x_{\ga_n}\rtra\hx\in\co$.  Then $\Psi(x\sgan)\rtra0$ and 
by the
lower semicontinuity $\Psi(\hx)=0$.   Since
$$\gF(x\sgan)-\gan\Psi(x\sgan)\ge
M\sgan-\gd\sgan\ge\sup_{\{\Psi=0\}}\gF(x)-\gd\sgan$$ 
and $\gF$ is
upper semicontinuous, \eqbdx) holds.\enddemo

\rpre\reba 
\rem{Remark\ {\rm \li}}  We record, for later use, 
some observations concerning maximum points $\ph$ of
$u(x)-v(y)-\gf(x-y)$ over $\co\times\co$ for other choices 
of $\gf$
besides $(\ga/2)|x|^2$ and the implications of Theorem 
\tha.  We are
assuming that $u$ is upper semicontinuous, $v$ is
lower semicontinuous, and $\gf\in C^2$.   In this 
application, the
matrix $A$ of Theorem \tha\ has the form 
$$A=\pmatrix  \ Z&-Z\cr-Z&\ Z\endpmatrix,\ti{where} 
Z=(D^2\gf)(\hx-\hy),\eq$$ 
\re\eqbmx 
so  
$$A^2=\pmatrix\ 2Z^2&-2Z^2\cr-2Z^2&\ 2Z^2\endpmatrix\eq$$ 
\re\eqbnx 
\noindent and $\vno A=2\vno Z$. 
\endrem

 Choosing $\gep=\oo{\vno A}$ in \eqbh), we conclude that 
there are
$X,\ Y\in\m$ such that  
$$(D\gf(\hx-\hy),X)\in\jtpco u(\hx),\qquad
(D\gf(\hx-\hy),Y)\in\jtmo v(\hy),\eq$$ 
\re\eqbo
\noindent and 
$$-2\vno
A\pmatrix I& 0\cr0&I\endpmatrix\le \pmatrix \ X&\ 0\cr\ 
0&-Y\endpmatrix \le
\pmatrix \junk&-(\junk)\cr-(\junk)&\junk\endpmatrix;\eq$$ 
in particular,
$$\vno
X,\vno Y\le 2\vno A\quad\text{and}\quad X\le Y.\eq$$

\subheading{Notes on \S 3}
This section is self-contained except for the proof of 
Theorem
3.2, which is due  to M. G. Crandall and H. Ishii 
[\craish].  This
result distills and sharpens  the essence  of a line of 
development
that runs through   several of the  references listed 
below; the
proof is explained in the appendix for  the reader's 
convenience. 
We want to emphasize at the outset that our  presentation 
of the
proof of the comparison theorem, based as it is on Theorem 
\tha, 
is but one  among several possibilities.  Moreover, other
approaches may  be useful in other situations.  For 
example, one
may use elements of the proof  of Theorem \tha\ directly  
and look
at corresponding manipulations  on solutions---this 
approach seems
to be the leading presentation for the infinite-dimensional
extensions of the theory and allows for adaptations on the 
 test
functions and the notion itself and avoids the ``semijets" 
or
generalized  derivatives.  Another presentation selects 
another
building block, namely, regularization by  supconvolutions 
(which
occur in the proof of Theorem \tha) and its effect on 
viscosity 
solutions; this procedure is especially helpful in some 
particular 
problems (integrodifferential  equations, regularity 
issues) and 
stresses the regularization procedure.  This procedure has 
some
analogues with the use of mollification in the study of 
linear
partial differential  equations.  With apologies for 
confusing the
reader, we are saying that a full grasp    of all the 
elements of
the proof of Theorem \tha\  and exposure to other  
presentations 
(e.g., R. Jensen [\jenb] and  H. Ishii and P. L. Lions 
[\lii]) may
prove valuable.

The first uniqueness proofs for viscosity solutions were 
given for
first-order  equations in M. G. Crandall and P. L. Lions 
[\clb]
and then  M. G. Crandall, L. C. Evans, and P. L. Lions 
[\cel].  The
second-order case remained  open for quite a while during 
which the
only evidence that a general theory could be developed was 
in 
results for Hamilton-Jacobi-Bellman equations  obtained by 
P. L.
Lions [\pllv, \pllc].   The proof in these works involved 
ad hoc
stochastic control verification arguments.

A breakthrough was achieved in the second-order theory by R.
Jensen [\jen] with the  introduction of several key 
arguments; some
of these were simplified in R. Jensen, P. L.  Lions, and 
P. E.
Souganidis [\jls], P. L. Lions and P. E.  Souganidis 
[\pllso], and
R. Jensen [\jenb].   In particular, the use of the 
``supconvolution" 
regularization (see the appendix), which  is a
standard tool in convex and nonsmooth analysis, was 
somewhat 
``remise au go\^ut dujour"  by J. M. Lasry and P. L.
Lions [\lali].  

Progress in understanding these proofs so as to be able to 
 handle 
more examples was made  by H. Ishii [\ishse] who 
introduced matrix
inequalities of the  general form \eqbg).  This work 
contains an 
example (following [\ishse, Theorem 3.3]) showing the  
optimality of condition \eqbl) in the sense that 
we cannot replace the right-hand side of \eqbl) by, for 
example, 
$\go(\ga|x-y|^\theta+|x-y|)$ with $\theta<2$ (but see \S 
\variations.A).
Estimating the left-hand side of \eqbl) differently, for 
example
in the form  $\ga g(|x-y|)+|x-y|$, it is possible to prove 
uniqueness in cases in which $g(r)/r^2$ is unbounded near 
$r=0$ 
by refining the arguments of [\revisit].  A more  complete
understanding together with the sharpest structure 
conditions (up
to the present) were achieved  in [\lii], which also 
contains many
examples. We mention that the structure condition \eqbl) 
exhibited in this
 section pays no attention to whether the equation is of 
first order or possess 
stronger ellipticity properties.  In the case of uniformly 
ellipitic $F$, much
 more can
be done [\lii]. See also Jensen [\jenb] as  regards 
structure conditions. 

A  useful sharpening  and improved organization of the 
analytical
essence of the theory were contributed by [\cra]. The  
presentation 
here was based on Theorem \tha; as remarked above, it 
follows from 
[\craish] and the generality and proofs presented there 
represent 
the current state of the  art (see the appendix).

\medskip \resultno=0 \equationno=0 
\spre\perron 
\heading \sli.  Perron's method and
existence \endheading
 \def\secno{\the\sectionno}  
 Let $\gO$ be an arbitrary open
subset of $\RN$. By a  solution (respectively, 
subsolution, etc.) of the Dirichlet
problem  
$$
\left\{\matrix\format\r&\quad\l\\
F(x,u,Du,D^2u)=0 & \text{in }\gO,\\
u=0 & \text{on }\bgo\endmatrix\right.\leqno(\roman{DP})$$ 
we
mean a function $u\in C(\goc)$ (respectively, $u\in 
\roman{USC}(\goc)$,
etc.) that is a (viscosity) solution (respectively, 
subsolution,
etc.) of $F=0$ in $\gO$ and satisfies $u(x)=0$ 
(respectively,
$u(x)\le 0$, etc.) for $x\in\bgo$.   (We note  that this 
formulation imposes  the boundary condition in a strict 
sense
that  will be  relaxed in  \S\boundary.)

 Recall that we always assume that $F$ is  \good\  and,
{\it unless otherwise said\/}, continuous .  To discuss 
Perron's
method, we will use the following notations:  if $u\: 
\co\rtra
[-\infty,\infty]$ where $\co\subset\RN$, then
$$\cases
u^*(x)=\lim\nolimits_{r\downa
0}\sup\{u(y):y\in\co\ti{and}|y-x|\le r\},\\
u_*(x)=\lim\nolimits_{r\downa
0}\inf\{u(y):y\in\co\ti{and}|y-x|\le r\}.\endcases\eq$$ 
One calls $u^*$  the upper semicontinuous envelope of $u$; 
it is the
smallest upper semicontinuous function (with values in
$[-\infty,\infty]$) satisfying $u\le u^*$.   Similarly, 
$u_*$ is
the lower semicontinuous envelope of $u$.
 \rpre\thp 
\myProclaim
Theorem \li.  {\rm (Perron's Method).}  Let
comparison hold for {\rm (DP);} i.e., if $w$ is a 
subsolution of
{\rm (DP)}
and $v$ is a supersolution of {\rm (DP)}, then $w\le v$.  
Suppose also
that there is a subsolution $\ub$ and a supersolution 
$\ovp$ of {\rm (DP)}
that satisfy the boundary condition $\ub_*(x)=\ovp^*(x)=0$ 
for
$x\in\bgo$.  Then  
$$
W(x)=\sup\{w(x):\ub\le w\le\ovp
\text{ and $w$ is a subsolution of {\rm (DP)}}\}\eq$$
is a solution of {\rm (DP)}.

Theorem \thb\  provides conditions under which comparison 
holds for
(DP). \re\eqca  The proof consists of
two steps.  The first one is \rpre\lepa \myProclaim Lemma 
\li.  Let
$\co\subset\RN$ be locally compact, $F\in
\roman{LSC}(\co\times\R\times\RN\times\m)$,  and $\cf$ be 
a family of
solutions of $F\le 0$ in $\co$.  Let  
$w(x)=\sup\{u(x):u\in\cf\}$
and assume that 
  $w^*(x)<\infty$ for $x\in\co$.  Then $w^*$ is a solution 
of $F\le
0$ in $\co$.

Allowing $F$ to be merely lower semicontinuous does not 
affect the
proof and is used later.   The information produced while 
proving
this lemma is needed later as well, so we isolate the 
essential
point in another result.  \rpre\prcc   \myProclaim
Proposition \li.  Let $\co\subset\RN$ be locally compact, 
$v\in
\roman{USC}(\co)$, $z\in\co$, and $(p,X)\in\jtpco v(z)$.  
Suppose also that
$u_n$ is a sequence of upper semicontinuous functions on 
$\co$
such that
 $$\left\{
\matrix\format\r&\quad\l\\
\roman{(i)} & \text{there exists } x_n\in\co\text{ such
that } (x_n,u_n(x_n))\rtra (z,v(z)),\\
\roman{(ii)} & \text{if }
z_n\in\co\ti{and}\ z_n\rtra x\in\co,\ir{then} 
\limsup_{n\rtra\infty}
u_n(z_n)\le v(x).\endmatrix\right.\eq$$  \re\eqccx Then 
$$
\cases
\text{there exists $\hx_n\in \co,\ (p_n,X_n)\in
\jtpco u_{n}(\hx_n)$}\\
\text{such that $(\hx_n,u_n(\hx_n),p_n,X_n)\rtra
(z,v(z),p,X)$}.
\endcases\eq$$ \re\eqdx

\demo{Proof}Without loss of generality we put $z=0$.  By 
the assumptions, 
 for every $\gd>0$ there is an $r>0$ such that
$N_r=\{x\in\co:|x|\le r\}$ is compact and   $$v(x)\le
v(0)+\ip{p,x}+\tfrac 12\ip{Xx,x}+\gd  |x|^2\quad\text{for 
}x\in N_r.\eq$$
\re\eqcb By the assumption \eqccx)(i), there exists 
$x_n\in\co$
such that $(x_n,u_n(x_n))\rtra(0,v(0))$.  Let $\hx_n\in 
N_r$ be a
maximum point of the function  $u_n(x)-(\ip{p,x}+\oo
2\ip{Xx,x}+2\gd |x|^2)$ over $N_r$ so that 
$$
\eqalign{u_n(x)\le
u_n(\hx_n)+\ip{p,x-\hx_n}+
\tfrac12(\ip{Xx,x}-\ip{X\hx_n,\hx_n})+2\gd
(|x|^2-|\hx_n|^2)\hfill\cr
\hfill\text{for } x\in N_r.}\eq$$ \re\eqcc 
\noindent Suppose that
(passing to a subsequence if necessary) $\hx_n\rtra y$ as
$n\rtra\infty.$  Putting $x=x_n$ in \eqcc) and taking the 
limit
inferior as $n\rtra\infty$, we find 
$$v(0)\le\liminf_{n\rtra\infty}u_n(\hx_n)-\ip{p,y}-\tfrac 12
\ip{Xy,y}-2\gd |y|^2;$$ on the other hand,  by \eqccx)(ii)
 $\liminf u_n(\hx_n)\le v(y)$ while  \eqcb) implies
$v(y)-\ip{p,y}-\oo 2\ip{Xy,y}-2\gd |y|^2\le 
v(0)-\gd|y|^2$.  We
conclude that
$$\eqalign{v(0)&\le\liminf_{n\rtra\infty} 
u_n(\hx_n)-\ip{p,y}-\tfrac 12
\ip{Xy,y}-2\gd |y|^2\cr&\le v(y)-\ip{p,y}-\tfrac 
12\ip{Xy,y}-2\gd
|y|^2\le v(0)-\gd|y|^2.\cr}$$ From the extreme 
inequalities we
learn $y=0$, so $\hx_n\rtra0$ (without passing to a 
subsequence), and
then from the first inequality and \eqccx)(ii) one sees 
that 
$v(0)=\lim_{n\rtra\infty}u_n(\hx_n)$.  Since we have $(p +
4\gd
\hx_n+X\hx_n,X+4\gd I)\in \jtpco u_n(\hx_n)$ for large 
$n$, we are nearly done. 
To conclude, we merely note that the set  of 
$(q,Y)\in\RN\times\m$
such that there exists $z_n\in \co,\ (p_n,X_n)\in \jtpco
u_{n}(z_n)$ such that $(z_n,u_n(z_n),p_n,X_n)\rtra 
(0,v(0),q,Y)$ is
closed and contains $(p,X+4\gd I)$ for $\gd>0$ by the 
above.\enddemo

\demo{Proof of Lemma \rm \lepa}  With the notation of the
lemma, suppose that $z\in\co$ and $(p,X)\in \jtpco 
\ws(z)$.  We
seek to show that $F(z,\ws(z),p,X)\le 0$.   It is clear 
that we may
choose a sequence $(x_n,u_n)\in\co\times\cf$ such that
$(x_n,u_n(x_n))\rtra (z,\ws(z))$ and that \eqccx) then 
holds with
$v=\ws$.  Hence, by the existence of data satisfying  
\eqdx) and the
fact that each $u_n$ is a subsolution, we may pass to the 
limit in
the relation $F(\hx_n,u_n(\hx_n),p_n,X_n)\le 0$ to find
$F(z,\ws(z),p,X)
\le 0$ as desired.

  The second step in the proof  of
Theorem \thp\ is a simple construction that we now 
describe. 
Suppose that $\gO$ is open, $u$ is a solution of $F\le 0$, 
and
$u_*$ is not a solution of $F\ge 0$; in particular, assume
$0\in\gO$ and we have  
$$
F(0,\us(0),p,X)<0\quad\text{for
some }(p,X)\in\jtmo \us(0).\eq
$$ 
Then, by continuity,
$u\sgd(x)=\us(0)+\gd+\ip{p,x} +
\tfrac12\ip{Xx,x}-\gamma|x|^2$ is a classical
solution of $F\le 0$ in $B_r=\{x:|x|< r\}$ for all small 
$r,\gd,
\gamma>0$.  Since  $$u(x)\ge\us(x)\ge
\us(0)+\ip{p,x}+\tfrac12\ip{Xx,x}+o(|x|^2),$$ if we choose
$\gd=(r^2/8)\gamma$  then $u(x)>u\sgd(x)$ for $r/2\le 
|x|\le r$  if $r$
is sufficiently small and then, by Lemma \lepa,  the 
function
$$U(x)=\cases \max\{u(x),u\sgd(x)\} & \text{if} |x|<r,\\
u(x)&\text{otherwise},\endcases$$ is a solution of $F\le 
0$ in $\gO$.  The
last observation is that in every neighborhood of 0 there 
are
points such that $U(x)>u(x)$; indeed, by definition, there 
is a
sequence $(x_n,u(x_n))$ convergent to $(0,\us(0))$ and then 
$$\lim_{n\rtra\infty}(U(x_n)-u(x_n))=u\sgd(0)-\us(0)=%
\us(0)+\gd-\us(0)>0.$$
We summarize what this ``bump" construction provides in the
following lemma, the proof of which consists only of 
choosing $r,\
\gamma$ sufficiently small.  
\enddemo

\rpre\lepb \myProclaim Lemma \li. Let $\gO$
be open  and  $u$ be solution of $F\le 0$ in $\gO$.  If  
$\us$
fails  to be a supersolution at some point $\hx$,  i.e., 
there
exists $(p,X)\in\jtmo \us(\hx)$ for which 
$F(\hx,\us(\hx),p,X)<0$,
then for any small $\gk>0$ there is a subsolution $U_\gk$ 
of $F\le
0$ in $\gO$ satisfying  
$$
\left\{\eqalign{&U_\gk(x)\ge u(x)\ti{and}
\sup\nolimits_{\gO}(U_\gk-u)>0,\cr
&U_\gk(x)=u(x)\quad\text{for 
}x\in\gO,|x-\hx|\ge\gk.\cr}\right.\eq
$$
\re\eqcd
 
\demo{Proof \,of\ \,Theorem \rm \thp}  With the notation 
of the
theorem observe that $\ub_*\le W_*\le W\le W^*\le\ovp^*$ 
and, in particular,
$W_*=W=W^*=0$ on $\bgo$.  
By Lemma \lepa\    $W^*$ is a
subsolution of (DP) and hence, by comparison, $W^*\le 
\ovp$.  
  It then follows from  the
definition of $W$ that $W=W^*$ (so $W$ is a subsolution).  
 If
$W_*$ fails to be a supersolution at some point 
$\hx\in\gO$, let 
$W_\gk$ be provided by Lemma \lepb.  Clearly $\ub\le 
W_\gk$ and $W_\gk=0$ 
on $\bgo$ for sufficiently small $\gk$.  By comparison, 
$W_\gk\le\ovp$ and since $W$
is the maximal subsolution between $\ub$ and $\ovp$,  we 
arrive at the contradiction
$W_\gk\le W$.  Hence $W_*$ is a supersolution of (DP) and 
then, by comparison
for (DP), $W^*=W\le W_*$, showing that $W$ is continuous 
and is a
solution.
\enddemo
\rpre\reb 
\rem{Remarks\  {\rm \li}}  The assumption that
$\gO$ was open in Lemma \lepb\ was used only to know that 
classical
subsolutions in small relative neighborhoods of points of 
$\gO$
were subsolutions in our sense.  In order to generalize 
this and to
formulate the version of Theorem \thp\ we will need in  \S
\boundary\   (which we did not do above for pedagogical 
reasons),
we now make some remarks and invite the reader to ignore 
them until
they are called for later.   Suppose $\co$ is locally 
compact, 
$G^+$, $G_-$ are defined on $\co\times\R\times\RN\times\m$ 
and have
the following properties:
 $G^+$ is upper semicontinous, $G_-$ is lower 
semicontinuous, and
classical solutions (twice continuously differentiable  
solutions
in the pointwise sense) of $G_+\le0$ on relatively open 
subsets of
$\co$ are  solutions of $G_-\le 0$.  Suppose, moreover, that
whenever $u$ is a solution of $G_-\le 0$ on $\co$ and $v$ 
is a
solution of $G^+\ge 0$ on $\co$ we have $u\le v$ on $\co$. 
 Then we
conclude that the existence of such a subsolution and 
supersolution
guarantees that there is a unique function $u$, obtained 
by the
Perron construction, that is a solution of both $G_+\ge 0$ 
and
$G_-\le0$ on $\co$.  The proof is unchanged except in 
trivial ways. 
\endrem

\rpre\exp 
\ex{Example \li}   Theorem \thp\ leaves open
the question of when a subsolution $\ub$ and a 
supersolution  $\ovp$
of (DP) that vanish on $\bgo$ can be found.  Let us 
consider this
problem for the  equation  \pre\eqpd 
$$
\left\{\matrix\format\r&\quad\l\\
u+G(x,Du,D^2u)=0 & \text{in }\gO,\\
u=0 & \text{on }\bgo\endmatrix\right.\eq$$ 
where $G$ is degenerate elliptic.  
Sometimes there may be an obvious choice for $\ub$ or 
$\ovp$; e.g.,
if  \pre\eqpe $$G(x,0,0)\le 0\iq{for} x\in\gO,\eq$$ then 
$\ub\equiv
0$ is a subsolution.   We assume \eqpe) and seek a 
supersolution 
``near $\bgo$" in the form \pre\eqpf 
$$u_1(x)=M(1-e^{-\gl
d(x)})\eq$$ where $M, \gl>0$ are parameters to be chosen
later  and  $$d(x)=\inf_{y\in\bgo}|y-x|$$ is the distance to
$\bgo$.  For $\gamma>0$ set $$\gO_\gamma=\set{x\in\gO: 
d(x)<\oo\gamma};$$ if
$\bgo$ is of class $C^2$ and $\gamma$ is large, then $d\in
C^2(\goc_\gamma)$.  We have     
$$\align
&u_1+G(x,Du_1(x),D^2u_1(x))\\
&\qquad\geq\!G(x\!,\lambda Me^{-\lambda 
d(x)}Dd(x)\!,\!\lambda Me^{-\lambda
d(x)}D^2d(x)-\lambda^2Me^{-\lambda d(x)}Dd(x)\otimes 
Dd(x)\!)\endalign
$$ 
in $\gO_\gamma$, where $p\otimes q$ is the matrix with
entries $p_iq_j$.   \pre\eqpg If $x\in\gO\sgl$, then $\gl 
d(x)\le
1$ and so  \pre\eqph $$ e^{-1}\le e^{-\gl d(x)}\le 1\iq{in}
\gO\sgl.\eq$$ Choose  $M>0$ so that  
\pre\eqpi
$$M(1-e^{-1})+G(x,0,0)>1\iq{in}\gO.\eq$$ Then we assume 
that there
is a large $\gl>0$ so that  
\pre\eqpj 
$$
\eqalign{&G(x,\messT x)\ge 0\cr
&\qquad\qquad\qquad\qquad\text{for }x\in\gO\sgl, \ 
Me^{-1}\le c\le M.\cr}
\eq$$ 
Then, by
assumption, $u_1$ is a classical solution of $u+
G(x,Du,D^2u)\ge 0$
in  $\gO\sgl$ and $u_1=0$ holds (continuously) on $\bgo$.

To complete the construction, choose $C>0$ so that \pre\eqpk
$$C+G(x,0,0)\ge 0\iq{in}\gO\eq$$ and \pre\eqpl
$$0<C<M(1- e^{-1});\eq$$ this is possible in view of 
\eqpi).  We
claim that  if $\ovp$ is defined by  
$$\ovp=\cases u_1\wedge C&\text{in
$\goc\sgl$,}\cr C&\text{in 
$\gO\setminus\goc\sgl$,}\endcases$$ 
then it is a
supersolution as desired.  Indeed, $\ovp=0$ on $\bgo$ is 
evident. 
The constant  $C$ is a supersolution of the equation in 
$\gO$ by
\eqpk) and the subset of $\goc\sgl$  on which $u_1=C$ does 
not
meet  $\partial(\gO\sgl)$ (because of \eqpl) and $u_1=0$ 
on $\bgo$
and  $u_1=M(1-e^{-1})$ when $\gl d(x)=1$).  Since the 
property of
being a supersolution is local and closed under finite 
minimums, we
are done. 
\endex

We have left the sufficient condition for our 
construction, \eqpj),
in a ``raw form." It is not useful to refine it further 
without
specializing $G$ further.  For example,  if $G$ has the 
form $G(x,p,X)=G_0(p,X)-f(x)$ where $f\in C(\goc)$,  in 
which the
$x$-dependence is  ``separated," the condition  that
$$\lim_{\gl\rtra\infty}G_0(\messT x)=\infty$$ uniformly in 
$c$ in
each interval $\gep\le c\le 1/\gep$, $\gep>0$, and $x$ in 
some 
neighborhood of $\bgo$ in $\goc$ will certainly allow us 
to satisfy
\eqpj).  In particular,  if $G_0$ depends only on $p$, 
$$\lim_{|p|\rtra\infty}G_0(p)=\infty$$ and $G_0(0)-f(x)\le 
0$, we can
uniquely solve the Dirichlet problem for $u+
G_0(Du)-f(x)=0$. Here we
rely on the fact that $|Dd(x)|=1$ at points of 
differentiability of
$d$; near $\bgo$  (which is all we need) this follows from 
 the
obvious fact that $Dd(x)=-n(x)$ on $\bgo$ where  $n(x)$ is 
the 
exterior normal to $\gO$ at $x\in\bgo$.   

In other cases, second-order terms dominate.  If
$$G(x,p,X)=-\trace(A(x)X)+\ip{b(x),p}-f(x)$$ where $A, b, 
f$ are
continuous, and  $$\ip{A(x)n(x),n(x)}> 0 \iq{on}\bgo,$$ 
then the
left-hand side of \eqpj) has the form 
$$c\gl^2\ip{A(x)Dd(x),Dd(x)}
+O(\gl) \iq{as}\gl\rtra\infty.$$ Since $Dd(x)=-n(x)$ on 
$\bgo$, the
coefficient of $\gl^2$ is positive near  $\bgo$ and it is 
easy to
achieve \eqpj) by taking $\gl$ large.    Note that  we are 
thus
able to assert the unique existence of a solution to  
\eqpd) in
this case (provided comparison holds) even though $A$ may be
completely degenerate inside $\gO$.

The reader can invent an unlimited array of examples. It 
should be
noted that  one may often  produce both a subsolution and a
supersolution by constructions  like the above  and then 
one need
not assume $u\equiv 0$ is a subsolution. One may also 
usually take
maximums and minimums of operators for which \eqpj)  can be
verified and stay within this class.  Thus, for example, 
if  
$$G(x,p,X)=\max\{-|p|^{2\theta-\gep}-|\trace\,
X|^{\theta-1}{\trace\, X} -f(x),|p|^\ga-g(x)\}$$ and 
$0<\gep\le
2\theta$, $\ga>0$, $f,g\in C(\goc)$, and $f,g\ge0$, then 
\eqpd)  has
a unique solution.

\subheading{Notes on \S 4} The combination of Perron's 
method and viscosity
solutions was introduced by H. Ishii [\isp].  The 
definition of 
a viscosity subsolution (respectively, supersolution) $u$ 
in 
[\isp] was that $u^*$ is a subsolution (respectively, 
$u_*$ is a 
supersolution) in the current  sense.  Solutions are then 
functions that are both a
subsolution and   a supersolution  and continuity   is not 
required.  (Note that then 
the characteristic function of the rationals is a solution 
of $u'=0$).  
  With this 
setup,  Perron's method does not require the comparison 
assumption and the
statements become more elegant (see [\isp]).

 We mention some
other approaches  to existence, for even if they are in 
general
much more complicated and of a  more limited scope, they 
can be
useful in some delicate situations.  For example,  one can 
use
formulas from control and differential games to write 
explicit
solutions   for approximate equations and  then use limiting
arguments; this approach  is used in M. G. Crandall and P. 
L. Lions
[\clii, Parts III and V] and D. Tataru [\tata].  

Two other approximations that have been used are 
discretization
and elliptic  regularization (for first-order equations; 
this is
the method of ``vanishing viscosity" and its relation to 
the theory
accounts for the term ``viscosity solutions").  Having 
solved  an
approximate problem, one then needs to pass to the limit 
(with some
a priori estimates---but see \S \rlimits!).  Existence
schemes of this sort have been used in P. L. Lions 
[\plla],  M. G.
Crandall and P. L. Lions [\cln], G. Barles [\barlesa, 
\barlesb], H. Ishii 
[\ishiold, \ishiolda],
P. E. Souganidis [\takc], and I. Capuzzo-Dolcetta and  P. 
L. Lions
[\capl].

\medskip

\resultno=0 \equationno=0 \spre\variations  
\heading \sli.  Comparison  and Perron's method:
Variations\endheading
\def\secno{\the\sectionno} 
\subheading{\num{\secno.A.}   Comparison with more
regularity} 
Suppose that the supersolution $u$ in Theorem \thb\  is more
regular; for example, $u\in C^\gamma(\gO)$ where 
$C^\gamma(\gO)$
denotes the set of functions that are H\"older continuous
with exponent $\gamma\in(0,1]$ on each compact subset of 
$\gO$.
Then the assertions of Theorem \thb\  remain valid if \eqbl)
is weakened to 
$$
\cases
F(y,r,\ga(x-y),Y)-F(x,r,\ga(x-y),X)\le
\go(\ga|x-y|^\theta+|x-y|) \\
\ti{whenever} x,y\in\goc,\
r\in\R,\ X,Y\in\m,\ti{and} \eqbi) \ti{holds}\endcases\eq$$ 
\re\eqdl 
\noindent for some $\theta>2-\gamma$.  Indeed, 
in the notation of the proof of Theorem \thb,
$$u(x)-v(y)-\ot\ga\xmys\le u(\hx)-v(\hy)-\ot\ga|\hx -\hy|^2
$$ implies, upon putting $x=y=\hy$ and recalling that $\hx$
will remain in a compact  subset of $\gO$ as 
$\ga\rtra\infty$,
$\ga|\hx -\hy|^2\le 2(u(\hx)-u(\hy))\le C|\hx-\hy|^\gamma$ 
and
then $\ga|\hx-\hy|^\theta\rtra 0$ as $\ga\rtra\infty$ 
provided
only $\theta >2-\gamma$, so the rest of the proof is 
unchanged.

More generally and precisely, if $u$ is uniformly continuous
with the modulus $\gr/2$, i.e., 
$|u(x)-u(y)|\le\gr(|x-y|)/2$, we
have $\ga|\hx -\hy|^2-\gr(|\hx-\hy|)\le 0$.  Defining 
$h(\ga)=\sup\{r:\ga r^2-\gr(r)\le 0\}$,  we will succeed if
$\go(\ga|x-y|^\theta+|x-y|)$ in \eqdl) is replaced by 
$\go(\ga k(|x-y|)+|x-y|)$ where $\ga k(h(\ga))\rtra0 $ as
$\ga\rtra\infty$.  Since the best possible modulus for
nonconstant functions behaves like $\gr(r)=cr$, in which
case $h(\ga)=c/\ga$ and $k(r)=o(r)$ succeeds, the Lipschitz
continuous  case is limiting with respect to these
arguments (however, see below for the $C^1$ case).  

 If we assume still more regularity of $u$,
e.g., $Du\in C^\gamma(\gO)$ where $0\le\gamma\le1$ (with 
$\gamma$=0 
meaning continuity and  $\gamma$=1 meaning  Lipschitz
continuity of $Du$), we need to adapt our strategy to 
obtain 
further results.  

Let $K\subset\gO$ be a compact
neighborhood of the set  of maximum points of $u-v$ in
$\gO$. If $Du\in C^\gamma$  there exists 
a sequence $\psin$
 in $C(\goc)\cap C^2(K)$ with the  properties  
\pre\eqdaa $$u-\psin\rtra
0\iq{as}n\rtra\infty\iq{uniformly on}\gO\eq$$ and  
\pre\eqdab
$$
\cases
|Du-D\psin|\le
Cn^{-\gamma} & \text{on }K,\\
|D^2\psin|\le Cn^{1-\gamma} & \text{on }
K,\\
|D^2\psin(x)-D^2\psin(y)|\le
Cn^{2-\gamma}|x-y| & \text{for }x,y\in K\endcases\eq$$ as
$n\rtra\infty$ where $C$ is a constant.  
The $\psin$ may then be
constructed by mollification of $u$ near $K$ and then 
extended.

Let $\ph$ be a maximum point of  \pre\eqdac
$$(u(x)-\psin(x))-(v(y)-\psin(y))-{\ga\over 2}|x-y|^2;\eq$$
as before, we may assume that $\hx,\ \hy\in K$ for large
$\ga$ and $n$.  Since  $Du(\hx)-D\psin(\hx)=\ga(\hx-\hy)$
and \eqdab) holds, 
we have \pre\eqdag $$\ga|\hx-\hy|\le 
Cn^{-\gamma}\iq{as}n\rtra\infty.\eq$$ Modifying the 
arguments
in the proof of Theorem \thb\ in an obvious way, we will 
still be able to establish comparison if we can show that
\pre\eqdah $$
\aligned
&(F(\hy,r,\ga(\hx-\hy)+D\psin(\hy),Y+D^2\psin(\hy))\\
&\qquad- 
F(\hx,r,\ga(\hx-\hy)+D\psin(\hx),X+D^2\psin(\hx)))^+\rtra 0
\endaligned\eq$$
as $n,\ga\rtra \infty$ in some appropriate fashion when $r$
remains bounded  and \eqbi), \eqdab), \eqdag) hold.

It is not very informative to try and analyze this condition
in  general.  Instead, let us note that in the linear case
(see Remark \recx) $$F(x,r,p,X)=-\trace(\ut
\gS(x)\gS(x)X)+\ip{b(x),p}+r-f(x)$$ with $f\in C(\gO)$, the
first-order terms are harmless from the point of view  of
verifying  \eqdah) if $b$ is continuous.  Via Remark \recx, 
the second-order terms contribute  at most  $$3\ga\trace(\ut
(\gS(\hx)-\gS(\hy))(\gS(\hx)-\gS(\hy))+\trace(-A(\hy) 
D^2\psin(\hy)+
A(\hx)D^2\psin(\hx))$$ when estimating \eqdah) above.  Here
we have put $A(x)=\ut\gS(x)\gS(x)$. Assuming  $\gS\in
C^\gl(\gO)$ (so $A\in C^\gl(\gO))$, we invoke \eqdab) and
estimate the  above expression by a constant times
$$\ga |\hx-\hy|^{2\gl}+n^{1-\gamma}
|\hx-\hy|^\gl+n^{2-\gamma}|\hx-\hy|.$$ By \eqdag),
this expression may, in turn,  be  estimated above in the 
form
 \pre\eqdam
$$
\left(\ga\left(\frac
{n^{-\gamma}}{\alpha}\right)^{2\lambda}+n^{1-\gamma}
\left(\frac
{n^{-\gamma}}{\alpha}\right)^\lambda+n^{2-\gamma}
\left(\frac
{n^{-\gamma}}{\alpha}\right)\right)
\iq{as}n\rtra\infty.\eq$$
Putting $\ga=n^\gb$, this becomes 
$$\left( \np{\gb-2\gl(\gamma+\gb)}+\np{-\gl(\gamma+\gb)+
1-\gamma}+\np{2(1-\gamma)-\gb} \right)\iq{as}n\rtra\infty.$$
Regarding $\gl$ as fixed so as to see what we require of
$\gamma$, we note that we may  make all of the exponents 
above
vanish by the choices $\gamma=(1-2\gl)/(1-\gl)$,
$\gb=2\gl/(1-\gl)$.  (The case $\gl>\tfrac12$ was treated 
above, so
here we are concerned with $\gl\le \tfrac12$.) If we 
increase $\gamma$ from this
value, so that $\gamma>(1-2\gl)/(1-\gl)$, it follows that 
comparison holds.  

Comparison is
thus assured if $\gS\in C^\gl$ and $Du\in C^\gamma$ under 
the relation 
$\gamma>(1-2\gl)/(1-\gl)$.  The  limit cases 
$\gl=\tfrac12$  and 
$\gl=0$ need to be discussed separately, for  $\gamma=0$ 
suffices if $\gl=
\tfrac12$ and 
 $\gamma=1$ suffices if $\gl=0$.  In the event that 
$\gamma=0$, the first relation
in \eqdab) should be replaced by $|u-\psi_n|\le o(1)$ and 
\eqdag) by
$\ga|\hx-\hy|\le o(1)$ and then the proof runs as before 
for $\gamma=0$ but with
the constant $o(1)$ appearing everywhere. To treat the 
other limiting
case, $\gl=0,\ \gamma=1$, note that if $\gr$ is a modulus 
of continuity of
$\gS$ and $A$,   \eqdam) should be replaced by  $$\ga
\gr(|\hx-\hy|)^2+\gr(|\hx-\hy|)+n|\hx-\hy|$$  and \eqdag) 
should read $\ga|\hx-\hy|\le
C/n$.  We  conclude upon letting $n,\ga\rtra\infty$ in a 
manner so that
$\ga\gr(C/(n\ga))^2\rtra0$.

It is worthwhile to compare the uniqueness results above
with those that may be obtained when solutions are more
nearly  classical in the sense that they  possess second
derivatives in a strong enough sense.   To simplify the
exposition, we will always assume that \eqbk) holds.   Then,
of course, if we assume that one of $u$ or $v$ is $C^2$ (or
even everywhere twice differentiable), the comparison result
holds without other structure conditions on $F$ beyond
\goodness, since we can work directly with
maxima of $u-v$.  Unfortunately, this regularity is rarely
available, even for uniformly elliptic fully nonlinear
equations.  If we require less regularity of $u$ and 
$v$---but 
still much more than assumed in Theorem \thb---different 
structure 
conditions suffice to guarantee  the
comparison result via finer considerations about the
pointwise twice differentiability of functions.    

We present below a strategy that establishes comparison
without further assumptions on $F$ if both $u+c|x|^2$ and $
-v+c|x|^2$ are convex for some $c>0$.  It also shows that
comparison holds if $u,\ v\in W^{2,N}_{\text{loc}}(\gO)$ and
one of $u$ or $v$ lies in $W^{2,p}_{\text{loc}}(\gO)$ for 
some
$p>N$ and the structure condition
$$\cases
F(x,t,p,X)\text{ is uniformly continuous in $p$}\\
\text{uniformly for }
x\in\goc,\ X\in\m, \text{ and }r,p \text{ 
bounded}\endcases$$ 
holds.

 In both of these cases, the argument runs
as follows:  we may assume, without loss of generality, that
all maxima of $u-v$ lie in $\gO$ and then general
optimization results imply that there exists a sequence
$p_n\rtra 0$ such that $ u(x)-v(x)-\ip{p_n,x}$ has a strict
maximum at $\hxn$ and $\hxn\rtra \hx$ where $\hx$ is a 
maximum
point of $u-v$ over $\goc$.  Then one knows (cf. Lemma 
\lej) 
that if $r>0$ is
small enough, the set of maximum points of
$u(x)-v(x)-\ip{p_n,x}-\ip{q,x}$ in $|x-\hxn|<r$ as $q$
ranges over the ball $B_\gd$ contains a set of positive 
measure if $\gd<r^{-1}\gep_n$ where
$$\gep_n=\inf_{|y-\hxn|=r}\left[(u(\hxn)-v(\hxn))-
\ip{p_n,\hxn-y}-(u(y)-v(y))
\right].$$ Furthermore, either $u$ or $v$ is twice
differentiable a.e. (cf. Theorem A.2)  and therefore we 
can find maxima $z$
with $|\hxn-z|<r$ for some $q\in B_\gd$ that is a  point 
where, in the first case, $u$ and $v$ are  twice
differentiable and in the second case $v$ or $u$ is twice
differentiable.  In the first case we have 
$$F(z,u(z),Dv(z)+p_n+q,D^2v(z))\le 0\le 
F(z,v(z),Dv(z),D^2v(z)).$$  In the second case, if $v$ is 
twice 
differentiable the same inequality holds and otherwise we 
have
$$F(z,u(z),Du(z),D^2u(z))\le 0\le 
F(z,v(z),Du(z)-p_n-q,D^2u(z)).$$  
The conclusions are easily
reached upon sending $\gd,r$ to $0+$ and then $n\rtra
\infty$.  \medskip

\subheading{\num{\secno.B.}   Estimates from
comparison } 
We make some simple remarks that hold as soon as one
has comparison by any means; there are many variants of 
these
(including parabolic ones).
 For example, suppose $K>0$ and
$u-v\le K$ on $\bgo$ instead of $u\le v$ on $\bgo$.  Then
$u-K$ is also a subsolution since $\jtpo(u-K)=\jtpo u$
and $F(x,u-K,p,X)\le F(x,u,p,X)$ by \goodness. 
Thus $u-v\le K$ in $\goc$. In particular, the variant of
Theorem \thb\ in which $u\le v$ on $\bgo$ is dropped and the
conclusion is changed to $u-v\le\sup_{\bgo}(u-v)^+$ holds.

In a similar spirit, suppose we have solutions $u$ and $v$
of $F(x,u,Du,D^2u)
\le 0$,  and $\hat F(x,u,Du,D^2u)\ge 0$ in
$\gO$ where $F$ is \good,  satisfies \eqbl)
and \eqbk), $K>0$, and $ F(x,r,p,X)+K\ge\hat F(x,r,p,X).$ 
Then
 $w=v+\max(\sup_{\bgo}(u-v)^+,K/\gamma)$ is a solution of
$F(x,w,Dw,D^2w)\ge 0$ in $\gO$. Since we have comparison
for  $F$, we then conclude $u-w\le\sup_{\bgo}(u-w)^+=0$ or 
$u-v\le\max(\sup_{\bgo}(u-v)^+,K/\gamma)$.  In particular, 
if
$$u+G(x,u,Du,D^2u)-f(x)\le 0\quad\text{and }
v+G(x,v,Dv,D^2v)-g(x)\ge 0\quad\text{in }\gO,$$ 
\noindent $G$ is \good
\ and satisfies \eqbl), and $f,g\in C(\goc)$, then
$$u-v\le\max(\sup_{\bgo}(u-v)^{+}
,\sup_\gO(f-g)^+).  $$

\subheading{\num{\secno.C.} 
Comparison with strict inequalities and without coercivity 
in $u$} 
The condition \eqbk) was used in the proof of Theorem \thb
 \  in order to have \eqblx).  If we simply assume there 
is a
$\gd>0$ such that either $F(x,u,Du,D^2u)\le -\gd$ or
$F(x,v,Dv,D^2v)$ $\ge\gd$ we are in the same situation and
do not need \eqbk) to hold.  Moreover, if only
$F(x,u,Du,D^2u)\le 0$ and $F(x,v,Dv,D^2v)$ $\ge 0$ but for
$\gep>0$ we can find $\psi_\gep\in C^2$, $\gd_\gep>0$ such
that $|\psi_\gep|\le \gep$ and
$F(x,u+\psi_\gep,D(u+\psi_\gep),D^2(u+\psi_\gep))\le
-\gd_\gep$, we conclude that $(u+\psi_\gep)-v\le
\sup_{\bgo}(u+\psi_\gep-v)$ and then
$u-v\le\sup_{\bgo}(u-v)+2\gep$ and we recover the result as
$\gep\downa 0$.  This construction can be carried out in 
some
cases.  

\subheading{\num{\secno.D.} Comparison and existence of 
unbounded solutions on
unbounded  domains}
We first
illustrate a  method to establish comparison of unbounded
solutions in unbounded domains by showing that if $u$ and
$v$ grow at most linearly and solve  $$u+F(Du,D^2u)-f(x)\le
0\iq{and}v+F(Dv,D^2v)-f(x)\ge 0\iq{in}\RN\eq$$ \re\eqdm
\noindent where $f$ is uniformly continuous on $\RN$ 
(i.e., $f\in
\roman{UC}(\RN))$, then $u\le v$.  After this, we will 
prove that if
$u$ is a solution of $u+F(Du,D^2u)=0$ of linear growth, then
$u\in \roman{UC}(\RN)$ and finally that Perron's method 
supplies
existence.  Thus we will prove  \rpre\thda
\myProclaim  Theorem \li.  If $f\in \roman{UC}(\RN)$, then
$u+F(Du,D^2u)-f(x)=0$ has a unique solution $u$ that grows
at most linearly and $u\in\uc$.

\demo{Proof of comparison}  The proof proceeds in
two steps.  First we note that $f\in \roman{UC}(\RN)$ 
implies that
there is a constant $K$ such that 
$$\sup_{\rtn}(f(x)-f(y)-K|x-y|)<\infty\eq$$ \re\eqdn 
\noindent and 
then we show that
$$\sup_{\rtn}(u(x)-v(y)-2K|x-y|)<\infty.\eq$$ \re\eqdo 
\noindent By
the assumed linear growth, we have  $$u(x)-v(y)\le
L(1+|x|+|y|)\iq{on}\rtn\eq$$ \re\eqdp 
\noindent for some $L>0$.   We
choose a family $\gb_R$ of $C^2$ functions  on $\RN$
parameterized by $R\ge 1$ with the properties
 $$\left\{\eqalign{&{\text{(i)}}\ \gb_R\ge 0,\cr 
&{\text{(ii)}}\
\liminf\nolimits_{|x|\rtra\infty}\gb_R(x)/|x|\ge 2L, \cr& 
{\text{(iii)}}
\ |D\gb_R(x)|+\vno{D^2\gb_R(x)}\le C\ti{for} R\ge 1,\
x\in\RN,\cr &{\text{(iv)}}\
\lim\nolimits_{R\rtra\infty}\gb_R(x)=0\ti{for}x\in\RN,\cr}%
\right.\eq$$
\re\eqdq 
\noindent where $C$ is some constant.  In view of \eqdp) and
\eqdq)(ii),  the function  
$$\gF(x,y)=u(x)-v(y)-2K(1+\nos{x-y})^{1/2}-(\gb_R(x)+
\gbr(y))$$
  attains its maximum at some point
$(\hx,\hy)$.  Now either \eqdo) holds or for large $R$ we
have $\gF(\hx,\hy)>0$, which implies $$2K|\hx-\hy|\le
u(\hx)-v(\hy).\eq$$ 
\re\eqdr 
\noindent Noting that
$$(\hat p + D\gb_R(\hx),\hZ + D^2\gb_R(\hx)) \in \jtp 
u(\hx),$$
$$(\hat p - D\gb_R(\hy),- \hZ - D^2\gb_R(\hy)) \in \jtm 
v(\hy),$$
where
$$\hat p = (2KD_z(1 + |z|^2)^{1/2})|_{z = \hx - \hy}, \
\hZ = (2KD_z^2(1 + |z|^2)^{1/2})|_{z = \hx - \hy}, $$
we see that
$$u(\hx) + F(\hat p + D\gb_R(\hx),\hZ + D^2\gb_R(\hx)) 
\leq f(\hx), $$
$$v(\hy) + F(\hat p - D\gb_R(\hy),- \hZ - D^2\gb_R(\hy)) 
\geq f(\hy). $$
From this, using \eqdr) and observing that $\hat p$ and  
$\hZ$ are
bounded independently of $R \geq 1$, we have
\pre\eqds
$$\aligned
u(\hx) - v(\hy) &\leq f(\hx) - f(\hy)
+ F(\hat p - D\gb_R(\hy),- \hZ - D^2\gb_R(\hy)) \\
&\hphantom{\leq}- F(\hat p + D\gb_R(\hx),\hZ + 
D^2\gb_R(\hx))\\
& \leq K|\hx - \hy| + C \leq \tfrac12 (u(\hx) - v(\hy)) + 
C\endaligned\eq$$
where $C$ is a constant independent of $R \geq 1$.  
Therefore,
$u(\hx) - v(\hy)$ is bounded independently of $R \geq 1$.

 Since $\gF(x,y)\le\gF(\hx,\hy)\le u(\hx)-v(\hy)$,
we may then send $R\rtra \infty$ to conclude that
$u(x)-v(y)-2K(1+\nos{x-y})^{1/2}$ is bounded and thus that
\eqdo) holds.    With the information
\eqdo) in hand, we repeat the above line of argument, this
time assuming $u(\bx)-v(\bx)=2\gd>0$ holds for some $\bx$
and putting $$\gF(x,y)=u(x)-v(y)-\ot\ga\nos{x-y}-\gep(\nos
x+\nos y)\eq$$ \re\eqdu 
\noindent where $\gep,\ \ga$ are positive
parameters.  For small $\gep$ we have $\gF(\bx,\bx)\ge\gd$
and by \eqdo) $\gF$ will have a maximum $(\hx,\hy)$ at which
$$
\ot\ga\nos{\hx-\hy}+\gep(\nos\hx+\nos\hy)\le
u(\hx)-v(\hy)\le 2K|\hx-\hy|+C 
\leq {\ga\over4} |\hx -\hy|^2 + 4K^2/\ga + 
C\vphantom{\sum^M}\eq
$$
\re\eqdv
\noindent for some $C$. 
Moreover, there are $X,\ Y\in\m$ such that
$$(\ga(\hx-\hy)+2\gep\hx,X+2\gep I)\in\jtpc u(\hx),\ 
(\ga(\hx-\hy)-2\gep\hy,Y-2\gep I)\in\jtmc v(\hy)$$ and 
$$-3\ga\pmatrix I&0\cr0&I\endpmatrix\le\pmatrix \  X&\
0\cr
\ 0&-Y\endpmatrix\le3\ga\pmatrix \ I&-I\cr-I&\ 
I\endpmatrix.\eq$$ \re\eqdw 

As usual, we have
$$
\aligned
u(\hx)-v(\hy)\le&
f(\hx)-f(\hy)+F(\ga(\hx-\hy)-2\gep\hy,Y-2\gep I)\\
&-F((\ga(\hx-\hy)+2\gep\hx,X+2\gep
I).\endaligned\eq$$ \re\eqdz
Since $\gd\le\gF(\bx,\bx)\le u(\hx)-v(\hy)$, and since 
$X\le Y$ by
\eqdw), 
$$\gd\leq\go_f(|\hx-\hy|) + 
F(\ga(\hx-\hy)-2\gep\hy,X-2\gep I)
- F(\ga(\hx-\hy)+2\gep\hx,X+2\gep I),$$
where $\go_f$ is the modulus of continuity of $f$,   
and we will obtain a contradiction if we show that  the 
iterated limit inferior
of the right-hand side as $\gep\downa0$ and then
$\ga\rtra\infty$ is nonpositive.  This information is 
deduced
from \eqdv), first using the extreme inequalities, which
imply that $\ga|\hx-\hy|^2$ and $\gep(\nos\hx+\nos\hy)$
are bounded independently of $\ga\ge 1$ and $0<\gep\le1$. 
Thus $\gep\hx, \gep\hy\rtra0$ 
and $\alpha(\hat x-\hat y)$ remains bounded 
as $\gep\downa0$, while $|\hat x-\hat y|\to 0$
as $\alpha\to\infty$ uniformly in $\varepsilon>0$. 
Therefore, using 
(5.17) and the 
uniform continuity of $F$ on bounded sets, 
$$\gd\le\liminf_{\ga\rtra\infty}\liminf_{\gep\downa0}%
\go_f(|\hx-\hy|)=0,$$
which is a contradiction.
\enddemo

\demo{Proof that solutions lie in $\uc$}  We take
the easy way here, using that $F$ is independent of $x$, and
note that if  $y\in\RN$ and $u$ solves
$u+F(Du,D^2u)-f(x)=0$, then $w(x)=u(x+y)$ solves
$w+F(Dw,D^2w)-f(x+y)=0$, so that
$|u(x)-u(x+y)|\le\sup_{z\in\RN}|f(z+y)-f(z)|$ by \secno.B;
thus $u\in\uc$. 
\enddemo

\demo{Proof of existence}  We need only to produce
a linearly growing subsolution and supersolution and invoke
Perron's method.  Letting $\ou(x)=A+$$B(1+\nos x)^{1/2}$, we
note that $D\ou$ and $D^2\ou$ are bounded so long as $B$ is
bounded.  By \eqdn), $|f(x)|\le C+K|x|$ for some $C$.  If we
put $B=K$ and  $A=\sup_{\RN} |F(D\ou,D^2\ou)|+C$, $\ou$ is a
supersolution and $\uu=-\ou$ is a subsolution.
\enddemo

\subheading{Notes on \S 5} The discussion of \S 
\variations.A
regarding regularity of $Du$ is new in this framework, but 
such results were first  obtained by R. Jensen [\jenb] in 
a more complicated presentation.  The assertions of \S 
\variations.A
regarding twice differentiable solutions (classical or  
$W^{2,p}$, etc.) 
rely on various versions of the maximum principle as it 
evolved through
the  works of Aleksandrov [\aleb, \alec], Y. Bakelman 
[\baka], C.
Pucci [\pucci],  J. M. Bony [\bony], and P. L. Lions 
[\pllbony].  

Section \variations.C recalls a classical strategy and 
corresponds to remarks used in M. G.
Crandall and P. L. Lions [\clb],
 H. Ishii [\iseik], H. Ishii and P. L. Lions [\lii], 
$\ldots$.
More sophisticated uses of estimates from comparison 
arguments 
occur in the study of regularity questions---see \S 
\perspectives\ for 
references.

Section \variations.D is concerned with the growth of 
solutions at
infinity and its influence on comparison-uniqueness 
results.  The
relevance of the class $\roman{UC}(\RN)$ for general 
results of this  sort
was progressively understood in a series of papers by the 
authors
and we refer to M. G. Crandall and P. L. Lions [\clii, 
Part II] for a few
examples showing how natural this class is and how it 
interacts with
structure conditions on the nonlinearity.  Of course, if one
restricts the nonlinearity further, other  asymptotic 
behaviors can
be allowed.  See H. Ishii [\ish, \ishse], M. G. Crandall 
and P. L.
Lions [\cplld,  \clinv],  and
  M. G. Crandall,  R. Newcomb, and Y. Tomita [\cnt] for more
information in this direction;  in particular, general 
functions 
$F(x,r,p,X)$ can be used in place of $F(p,X)-f(x)$ 
(although 
subtleties concerning existence arise).

\resultno=0 \equationno=0 
\spre\rlimits 
\heading \sli. Limit operations with viscosity 
solutions\endheading 
\def\secno{\the\sectionno} 
Suppose we have a sequence $u_n$,  $n=1,2,\ldots$, of 
subsolutions of  an
equation $F=0$ on the locally compact set $\co$.    It turns
out that the following ``limit"   
$$\hU(z)={\limsup_{n\rtra\infty}}^*
u_n(z)=\lim_{j\rtra\infty}\sup\{u_n(x):n\ge j,\ 
x\in\co,\ti{and} |z-x|\le \tfrac 1j\}\eq$$ \re\eqeb 
\noindent in which,
roughly speaking, the ``lim\,sup" operation and the $^*$
operations are performed simultaneously,  will also be a
solution of $F\le 0$ in $\co$.  Indeed, since $\hU(x)<r$
only if there  are $\gep, j>0$ such that $u_n(z)<r-\gep$
for $n\ge j$ and $z\in \co$ with $|z-x|\le 1/j$,  it is
clear that $\{x\in\co:\hU(x)<r\}$ is open in $\co$ and thus
$\hU\in \roman{USC}(\co)$.    We have
 \rpre\prea \myProclaim Lemma \li.  Let $u_n\in 
\roman{USC}(\co)$ for
$n=1,2,\ldots$ and $\hU$ be given by {\rm \eqeb)}, 
$z\in\co$ and
$\hU(z)<\infty$.  If $(p,X)\in\jtpco \hU(z)$, then there
exist sequences  
$$n_j\rtra\infty, \qquad \hx_j\in\co,  \qquad
(p_j,X_j)\in\jtpco u_{n_j}(\hx_j)\eq$$
 \re\eqec 
\noindent such that $$(\hx_j,u_{n_j}(\hx_j),p_j,X_j)\rtra
(z,\hU(z),p,X).\eq$$ \re\eqed 
\noindent In particular, if each $u_n$
is a solution of $F\le 0$ and $\hU<\infty$ on $\co$, then
$\hU$ is a solution of $F\le 0$ on $\co$. 

\demo{Proof} By definition, there  are sequences   
$$n_j\rtra\infty, \
x_j\in\co\text{ such that } x_j\rtra z\text{ and }
u_{n_j}(x_j)\rtra\hU(z);$$  it is also clear that if 
$z_j\rtra
x$ in $\co$, then $\limsup u_{n_j}(z_j)\le\hU(x)$.  The
result now follows at once from Proposition \prcc.\enddemo

\rpre\reea 
\rem{Remark\  {\rm \li}}  If the
assumptions of the lemma are altered to assume instead that
$u_n\in \roman{LSC}(\co)$ and $(p,X)\in\jtmco v(z)$, then 
the
assertions are changed by replacing $\hU$ by  $$\underline
U(z)={\liminf_{n\rtra\infty}}\!\sas
u_n(z)=\lim_{j\rtra\infty}\inf\{u_n(x):n\ge j,\ 
x\in\co,\ti{and} |z-x|\le \tfrac1j\},\eq$$ \re\eqebx
 putting $\jtmco$ in place of $\jtpco$ in \eqec) and
supersolutions in place of subsolutions.
\endrem
 
\rpre\reeb 
\rem{Remark\ {\rm \li}}  In
fact, the above proof shows more---suppose $u_n$ is a
solution of a \good\   equation $F_n\le 0$ that
varies with $n$.  Then the conclusion is that $\hU$ is a
solution of $G\le 0$ where  
$$G(x,r,p,X)={\liminf_{n\rtra\infty}}\!\sas 
F_n(x,r,p,X);\eq$$ 
\re\eqee  
\noindent note that $F_n$ need not be continuous and that if
$F_n=F$ is independent of $n$ but discontinuous, then
$G=F_*$ is the lower semicontinuous envelope of $F$. 
Analogous remarks hold for supersolutions. 
\endrem

\rpre\reec 
\rem{Remark\ {\rm \li}} The above
results are related to uniform convergence as explained
next.  Let  $u_n$ be a sequence of functions on $\co$ and
$$\hU(x)={\limsup_{n\rtra\infty}}^*u_n(x),\qquad \uU(x)
={\liminf_{n\rtra\infty}}\!\sas u_n(x).\eq$$ \re\eqef  
\noindent Suppose
$\hU(x)=\uU(x)$ on $\co$, let $U(x)$ denote this common 
value
and assume that $-\infty<U(x)<\infty$ on $\co$.  Then $U$ is
continuous (since it is both upper and lower semicontinuous)
and $\lim_{n\rtra\infty}u_n(x)=U(x)$ uniformly on compact
sets.  Indeed, if this were not the case and uniform
convergence failed on some compact set $K$, there would be
an $\gep>0$ and sequences $n_j\rtra\infty$, $x_j\in K$ such
that $u_{n_j}(x_j)-U(x_j)>\gep$ or
$u_{n_j}(x_j)-U(x_j)<-\gep$.  Assuming $x_j\rtra x$ and 
using
the continuity of $U$, we would conclude that 
$|U(x)-U(x)|\ge \gep$, a contradiction. In order to prove
that $\hU=\uU$ one notes that $\uU\le\hU$ by definition and
typically uses comparison results to prove the other
inequality.  The next result provides  a simple example; a 
typical case in which the hypotheses of this result are 
easily verified is mentioned after the short proof.
\endrem

\rpre\thea \myProclaim Theorem \li.  Let $\gO$ be a bounded
open set in $\RN$, $H\in C(\RN)$, and $f\in C(\goc)$.  
Consider the problem  $$u+
H(Du)-\gep\gD u=f(x)\quad\text{in}\  \gO,\ \ \
u=0\ti{on}\partial\gO,\leqno(\roman{DP})_\gep$$ and assume 
that
$(\roman{DP})_\gep$ has a subsolution $\underline u\in 
C(\goc)$
and  a supersolution $\overline u\in C(\goc)$ independent
of $\gep\in[0,1]$ that vanish on $\partial\gO$.  Then
$(\roman{DP})_\gep$ has a solution $u_\gep$ for  
$\gep\in[0,1]$ and
$\lim_{\gep\downa 0}u_\gep(x)=u_0(x)$ uniformly for 
$x\in\goc$.

\demo{Proof} We know that (DP)$_\gep$ has a solution 
$u_\gep$ for
$\gep\in[0,1]$ by Perron's method and the assumptions. 
Moreover, $\hU(x)={\limsup_{\gep\downa 0}}^*u_\gep(x)$ and
$\uU(x)={\liminf_{\gep\downa0}}\sas u_\gep(x)$ are a 
subsolution
and a supersolution of (DP)$_0$ by the above and then
$\hU\le \uU$ by comparison while $\hU\ge \uU$ by
definition.  Thus $u_0=\hU=\uU$ and the convergence is
uniform.  

Harking back to the discussion of Example \exp, the reader 
may 
easily verify the hypotheses of Theorem \thea\ when 
$\bgo$ is of class $C^2$,  $H$ satisfies  
$\lim_{|p|\rtra\infty}H(p)=\infty$
and $H(0)-f\le 0$.  To do so,  set $\ub\equiv 0$ and 
$\ovp(x)=(\gl d(x))\wedge C$ where  $\gl$ is first 
chosen sufficiently large and then $C$ is chosen 
in a suitable manner.\enddemo

\subheading{Notes on \S 6}
The fact that viscosity solutions pass to the limit under 
uniform
convergence is an immediate consequence of  the 
definition.  This
consistency and stability property can also be seen as the 
analogue
(in the space of  continuous functions) of Minty's device 
for
monotone operators as developed by L. C. Evans [\evco, 
\evac].

The idea of using only ``half-relaxed" limits like $\bar 
U$ in
\eqeb) arises naturally when attempting to  pass to limits 
with
maxima and minima and has been used extensively in the 
calculus of
variations and homogenization theory.   The use of this 
concept in
the area of viscosity solutions was introduced by G. 
Barles and B.
Perthame [\bapa, \bapb] and  H. Ishii [\isbvp].  This 
passage to
the limit and the  associated  a  posteriori uniform 
convergence
through comparison of  semicontinuous functions has become 
one of
the main features of the theory.  Viscosity solutions 
were first used to prove 
theorems analogous to Theorem \tha\  (but more complicated) 
with applications 
by L. C. Evans and H. Ishii [\evish]; they were motivated 
by 
questions arising in probability (see also \S 
\perspectives).  
Proofs in the spirit we have presented 
are much simpler.

\medskip \resultno=0 \equationno=0 

\spre\boundary \heading \sli. General and
generalized boundary conditions\endheading 
\def\secno{\the\sectionno}  
\subheading{\num{\secno.A.}  Boundary conditions in the 
viscosity
sense} 
In this section we consider more
general boundary value problems of the form  
$$\cases (\roman E)\ \
\ \ F(x,u,Du,D^2u)=0\text{\quad in \  $\gO$} 
\cr(\roman{BC})\ \ 
B(x,u,Du,D^2u)=0\text{\quad on  
$\bgo$}\endcases\leqno(\roman{BVP})$$  
where (E)
denotes  the equation and (BC)  the boundary condition.  
Here $F$ is \good, as always, and $B$ is a
given function on $\bgo\times\R\times\RN\times\m$ that is
also to be \good.  While it is convenient to
allow $B$ to depend on $D^2u$ at this stage, as the reader
will see,  when we turn to existence and uniqueness theorems
$B=B(x,r,p)$  will be of first order.  For example, the
Dirichlet condition $u=f(x)$ on $\bgo$  arises from the
choice $B(x,r,p,X)=B(x,r,p)=r-f(x)$ of $B$ while the Neumann
condition $u_n=f(x)$  (here and later, $n(x)$ denotes the
outward unit normal to $x\in\bgo$) arises if
$B(x,r,p)=\ip{n(x),p} -f(x)$.  The Neumann condition is
generalized by the  ``oblique derivative" problem, in which
$B$ has the form $B(x,r,p)=\ip{\gn(x),p}-f(x)$ where
$\gn(x)$ is a  vector field on $\bgo$ satisfying
$\ip{\gn(x),n(x)}>0$ on $\bgo$.   

We have already given a meaning to statements like (BC) in
Definition \dea, which is, however, clearly inappropriate
here since it would only involve the behavior of $u$ on 
$\bgo$.  The simplest  viscosity definition of (BC) we
might try is
\rpre\debca
\dfn{Definition \li}  A function
$u\in\uscc$ is a subsolution of \bc\  in the strong
(viscosity) sense if  
$$B(x,u(x),p,X)\le
0\iq{for}x\in\bgo, (p,X)\in\jtpcgo u(x)\eq $$ 
\re\eqbcc 
\noindent and $v\in\lscc$ is a supersolution of \bc\ in the
strong (viscosity) sense if  
$$B(x,v(x),p,X)\ge
0\iq{for}x\in\bgo,\text{ and } (p,X)\in\jtmcgo v(x).\eq$$ 
\re\eqbcd 
\noindent Finally, $u\in C(\goc)$ satisfies  \bc\ in
the strong (viscosity) sense if it is  both a  subsolution
and a supersolution in the strong sense.
\enddfn

Interpreting (E) in the usual way and (BC) in the strong
sense, we obtain notions of subsolutions (which will lie in
$\uscc)$, supersolutions (which will lie in $\lscc)$ and
solutions (which lie in $C(\goc))$.   
   
As indicated by the
appearance of the modifier ``the strong
sense," it will be
necessary to relax this interpretation to obtain full
generality.  Before explaining this, however, we pause to
relate ``the strong sense" to the ``classical sense."

  It is important to recognize early on that if $u\in
C^2(\goc)$ (i.e., $u$ extends to a twice continuously
differentiable function in a neighborhood of $\goc$) and
$$
B(x,u(x),Du(x),D^2u(x))\le 0\quad\text{for }x\in\bgo,$$
then it is
not necessarily true that $u$ is a subsolution  of (BC). 
Indeed, in Remark \rebx\  we determined $\jtpco u(x)$ in the
situation where  $\co$ is an N-submanifold of $\RN$ with
boundary, $x\in\bgo$, and $\bgo$ is twice differentiable at
$x$.  In   particular, it follows from Remark \rebx\ that
if $S\le 0$ and $X_\gl$ is given by \eqzfx) with $u$ in 
place of $\gf$, 
\ then  $(Du(x)-\gl n(x),X_\gl
+\gm n(x)\otimes n(x))\in \jtpco u(x)$ for $\gl=\gm=0$ and
for $\gl>0$, $\gm\in\R$.  Thus, assuming that $B$ is
lower semicontinuous, $u$ will not be a subsolution
unless   $$B(x,u(x),Du(x)-\gl  n(x),X_\gl+\gm n(x)\otimes
 n(x))\le 0\iq{for}\gl\ge 0,\ \gm\in\R\eq$$ \re\eqbcdx 
\noindent (to
see that \eqbcdx) holds with $\gl=0$ and $\gm$ arbitrary,  
just let $\gl\downa0$ with $\gm$ arbitrary but fixed).  

Notice, however, that if $B(x,r,p)$ is first order and
$n\in\RN$, then
$B(x,u(x),Du(x))\le 0$ implies
$B(x,u(x),Du(x)-\gl n)\le 0$ for $\gl\ge 0$ provided
$$\gl\rtra B(x,r,p-\gl n)\text{ is nonincreasing in 
}\gl\ge 0.\eq$$
\re\eqbcdxx Moreover, by Remark \rebx, if  $u$ is
differentiable at $x$ and $(p,X)\in\jtpco u(x)$, then
$Du(x)-p$ is a ``generalized normal" to $\co$ at $x$.  Since
analogous remarks evidently hold for supersolutions, we
have found a  case in which classical satisfaction of (BC)
implies (BC) in the strong sense.  More precisely, we have 
   
proved \rpre\prbca 
\myProclaim  Proposition \li.  If $B\in
C(\bgo\times\R\times\RN)$ satisfies {\rm \eqbcdxx)} for 
every
$x\in\bgo,\ r\in\R,\ p\in \RN$ and $n$ is a generalized 
normal to $\co$
at $x$, $u\in C^1(\goc)$  and $B(x,u(x),Du(x))\le 0$ {\rm 
(}respectively,
$\ge 0)$ for $x\in\bgo$, then  $u$ is a subsolution
{\rm (}respectively, supersolution{\rm )} of {\rm (BC)} in 
the strong sense.

This result will be needed in the construction of
subsolutions and supersolutions below.

As mentioned above, it will be necessary to relax the
interpretation of the boundary condition ``in the strong
sense;''  one  way to motivate this is to consider what
happens when we try to pass to the limit in the following
situation, which is an analogue of Remark \reeb\ ``with
boundary:''
 Suppose $F_n\in C(\goc\times\R\times\RN\times\m)$, $B_n\in
C(\bgo\times\R\times\RN\times\m)$ for $n=1,2,\ldots$ and
$u_n$ be a subsolution of  $F_n=0$ in $\gO$, $B_n=0$ on
$\bgo$ in the strong sense.   Let $F_n\rtra F$, $B_n\rtra B$
uniformly on compact sets as $n\rtra \infty$ and
$u=\limsup_{n\rtra \infty}^* u_n$ be bounded.   What problem
does $u$ solve?  Of course, we know from Remark \reeb\ that
$u$ is a solution of $F\le 0$ in $\gO$.   However, when we
attempt to check the boundary condition $B(x,u(x),p,X)\le 0$
for $x\in\bgo$ and $(p,X)\in\jtpcgo u(x)$, the proof fails. 
We know, from Proposition \prcc, that (passing to a
subsequence if necessary), there exists $(p_n,X_n)\in\jtpcgo
u_n(x_n)$ such that $(x_n,u_n(x_n),p_n,X_n)\rtra
(x,u(x),p,X)$, but we do not know $x_n\in\bgo$ (and the
example below shows this cannot be achieved  in general). 
If $x_n\in\gO$, we have $F_n(x_n,u_n(x_n),p_n,X_n)\le 0$
(rather than $B_n(x_n,u_n(x_n),p_n,X_n)\le 0$).  However, in
this way we can conclude  that either 
$$
B(x,u(x),p,X)\le 0\quad\text{or}\quad
F(x,u(x),p,X)\le 0;$$
in other words, $u$ is a
subsolution for the boundary condition 
$$F(x,u,Du,D^2u)\wedge
B(x,u,Du,D^2u)= 0$$
in the strong sense (where
$a\wedge
b=\min\{a,b\}$ and $a\vee b=\max\{a,b\})$.  The following
example, in the context of a Neumann problem, shows that
this is the best that can be expected. 

\rpre\exbca 
\ex{Example\  \li} 
The     linear problem  $$-\gep u'' +u'+u=x+1
\eq$$ \re\eqbca on $(0,1)\subset\R$ with the
homogeneous Neumann conditions $$u'(0)=u'(1)=0\eq$$ 
\re\eqbce
\noindent has a unique classical solution $u_\gep$ for 
$\gep>0$ and it
follows from Proposition \prbca\  that $u_\gep$ satisfies
\eqbce) in the strong sense when we chose the function $B$
that defines the boundary condition \eqbce) to be given by
$B(1,r,p)=p$  and $B(0,r,p)=-p$;  note the monotonicity in
the direction of the exterior normals to $(0,1)$.   We will
see that the limit $\lim_{\gep\downa0}u_\gep=u$ exists
uniformly on $[0,1]$ and the limit $u$  does not satisfy 
\eqbce)  in the strong sense. 
\endex

It is elementary to compute 
$$u_\gep(x)=x+
{\e{\glm}-1\over\glp\paren{\e{\glp}-\e{\glm}}}\e{\glp
x}+ {1-\e{\glp}\over\glm\paren{\e{\glp}-\e{\glm}}}\e{\glm
x},$$ where $\gl_\pm=(1/2\gep)\paren{1\pm
(1+4\gep)^{1/2}}$.   Noting that $\glp\rtra\infty$ and
$\glm\rtra-1$ as $\gep\downa0$, one sees that  
$$u_\gep(x)\rtra
u(x)\equiv x+\e{-x}\iq{uniformly on} [0,1]\text{ as 
}\gep\downa0.$$ 
The function $u(x)=x+\e{-x}$ satisfies 
$$u'+u=x+1\text{ in }(0,1)\quad\text{and}\quad u'(0)=0$$  
in the classical
sense.  However, $u'(1)=1-1/e>0$ and therefore $u$ is not a
subsolution of $B(1,u'(1))=0$ in the strong sense.  Note,
however, that $(u'+u-(x+1))\wedge u'\le 0$ does hold at
$x=1$, in agreement with the above discussion.  

These considerations suggest the appropriate  definitions
regarding solutions of (BVP).   \rpre\debcb  

\dfn{Definition \li}   Let $\gO$ be an open subset of 
$\RN$, $F\in
C(\goc\times\R\times\RN\times\m)$, and $B\in
C(\bgo\times\R\times\RN\times\m)$ be  \good. 
Then $u$ is  a (viscosity) subsolution of (BVP) if
$u\in\uscc$ and  
$$\cases F(x,u(x),p,X)\le 0\quad\text{for }x\in\gO,
(p,X)\in\jtpcgo u(x),\\
F(x,u(x),p,X)\wedge
B(x,u(x),p,X)\le0\quad\text{ for }x\in\bgo, (p,X)\in\jtpcgo
u(x).\endcases\eq$$ 
\re\eqcf 
Similarly, $u$ is a supersolution
of (BVP) if $u\in\lscc$ and 
$$\cases F(x,u(x),p,X)\ge
0\quad\text{for }x\in\gO,(p,X)\in\jtmcgo u(x),\\
F(x,u(x),p,X)\vee B(x,u(x),p,X)\ge0\quad\text{for }x\in\bgo,
(p,X)\in\jtmcgo u(x).\endcases\eq$$ \re\eqcg 
\noindent Finally, $u$ is a
solution of (BVP) if it is both a subsolution and a
supersolution.
\enddfn

Another way to express these definitions is to introduce the
functions 
$$G_-(x,r,p,X)=\cases F(x,r,p,X)&\text{if $x\in\gO$,}\cr
F(x,r,p,X)\wedge B(x,r,p,X)&\text{if 
$x\in\bgo$}\endcases\eq$$ 
\re\eqbch
\noindent and 
$$G_+(x,r,p,X)=\cases F(x,r,p,X)&\text{if $x\in\gO$,}\cr
F(x,r,p,X)\vee B(x,r,p,X)&\text{if 
$x\in\bgo$,}\endcases\eq$$ 
\re\eqbci
\noindent and note that ``$u$ is a subsolution of (BVP)" 
means
exactly that $u$ is a solution of $G_-\le 0$ in $\goc$,
etc.  Observe that $G_-$ is lower semicontinuous and $G_+$
is upper semicontinuous.  In fact, $G_-$ is precisely the
lower semicontinuous envelope of the function 
$$G(x,r,p,X)=\cases F(x,r,p,X)&\text{if $x\in\gO$,}\cr 
B(x,r,p,X)&\text{if $x\in\bgo$,}\endcases$$ 
while $G_+$ is the
upper semicontinuous envelope of $G$.\smallskip

\noindent
\subheading{\num{\secno.B} Existence and
uniqueness for the Neumann problem} 
We are concerned here with the following choice 
$B=B(x,r,p)$ of the boundary condition $B$,
 $$B(x,r,p)=\ip{n(x),p}+f(x,r)\eq$$ \re\eqbcj where $f\in
C(\bgo\times\R)$.  With appropriate restrictions on $\gO$,
$F$, and $f$, we  prove that (BVP) has a unique solution by
a program analogous to that developed in \S\S \maximum\ 
and \perron.  Regarding $\gO$ we will always assume
that  
$$
\ctrblock
$\goc$ is a compact $C^1$ $N$-submanifold  with
boundary of $\RN$.\endblock \eq$$ 
\re\eqbcjx   
\noindent In addition, we
impose the uniform exterior sphere condition  
$$
\ctrblock
$\exists
r>0$ such that $b(x+rn(x),r)\cap
\gO=\emptyset$ for $x\in\bgo$\endblock\eq$$ 
\re\eqbck 
\noindent where $b(x,r)$
denotes the closed ball of radius $r$ centered at $x$.  For
$F$ we will require \eqbk), restated here for convenience,
$$ 
\eqalign{
&\gamma(r-s)\le F(x,r,p,X)-F(x,s,p,X)\cr
&\qquad\qquad\qquad\qquad\text{for } r\ge s,\
(x,p,X)\in\goc\times\RN\times\m\cr}\eq$$  \re\eqbckx 
\noindent as well as 
uniform continuity in $(p,X)$ near $\bgo$; more precisely,
we require that there be a neighborhood $V$ of $\bgo$
relative to $\goc$ such that 
$$\cases|F(x,r,p,X)-F(x,r,q,Y)|\le\go(|p-q|+\vno{X-Y})\\
\text{for }x\in V,\ p,q\in\RN,\ X,Y\in\m\endcases\eq
$$
\re\eqbcl 
\noindent for some $\go\:[0,\infty)\rtra[0,\infty]$
satisfying $\go(0+)=0$.  In addition, we strengthen \eqbl)
to the condition 
\pre\eqbcm
$$\left\{\ctrblock
$F(y,r,p,Y)-F(x,r,p,X)\le
\go(\ga\xmys+|x-y|(|p|+1))$ \quad whenever $x,y\in\goc,\
r\in\R,\ p\in\RN,\  X,Y\in\m$, and \eqbi)\ holds
\endblock\right.\eq$$  
where we may use the same $\go$
as in \eqbcl).  Observe that \eqbcm) implies \eqbl).  For
$f$ we assume that  $$f(x,r)\text{ is nondecreasing in }
r\text{ for } x\in\bgo.\eq$$ \re\eqbcn 
\noindent We have: \rpre\thbca
\myProclaim Theorem \li.  Let {\rm \eqbcj)--\eqbcn)} hold. 
 If $u$
is a subsolution of {\rm (BVP)} and $v$ is a supersolution 
of
{\rm (BVP)}, then $u\le v$ in $\goc$.  Moreover, {\rm 
(BVP)} has a
unique solution. 

The proof consists of verifying the comparison assertion,
then producing a subsolution and a supersolution, and
invoking Perron's method.  We begin with the comparison
proof.  

\demo{Proof of comparison} 
The strategy involves two types of approximations.    As a
first step, we produce approximations $u_\gep, v_\gep$ such
that $u\sep\rtra u,\ v\sep\rtra v$ uniformly as 
$\gep\downa 0$ and
$u\sep, v\sep$ are a subsolution and a supersolution  of
(BVP) with $B(x,r,p)$ replaced by $B(x,r,p)+\gep$ and
$B(x,r,p)-\gep$ respectively.   (Please note that {\it B
itself will retain its meaning below.})
If we then can prove $u\sep\le v\sep$
we find $u\le v$ in the limit $\gep\downa 0$.  After this, 
we
arrange to deal only with the inequalities involving   $F$
(and not $B$)  via another perturbation built into the
``test functions" we will use.
\enddemo
Let us proceed with the first step.  We need the following
elementary lemma \rpre\lebca \myProclaim Lemma \li.  Let
{\rm \eqbcjx)}  hold and $\gn\in C(\bgo,\RN)$ satisfy
$\ip{n(x),
\gn(x)}>0$
on $\bgo$.  Then there is a function
$\gf\in C^2(\goc)$ such that
$$
\ip{\gn(x),D\gf(x)}\ge1\quad\text{for }x\in\bgo\text{ and }
\gf\ge0\ \text{on }\goc.\eq
$$
\re\eqbcp

The proof of the  lemma is deferred until the theorem is
established.  Choosing $\gn=n$ in Lemma \lebca, let $\gf$ be
provided by the lemma; we may assume that the support of
$\gf$ lies in $V$.  Let $u$ be a solution of $G_-\le 0$ in
$\goc$ (that is, a subsolution of (BVP)) and $v$ a
solution of $G_+\ge 0$, and put
$$u\sep(x)=u(x)-\gep\gf(x)-C\quad\text{and}\quad 
v\sep(x)=v(x)+\gep\gf(x)+C\eq$$
\re\eqbcq 
\noindent where $C=C(\gep)$ will be chosen later on.  Using
\eqbckx) and \eqbcl),   for $x\in V$ and $(p,X)\in\jtpcgo
u_\gep(x)$ we have 
$$F(x,u\sep(x),p,X)\le F(x,u(x),p+\gep
D\gf(x),X+\gep D^2\gf(x))-\gamma C+\go(\gep M),$$ 
\noindent where $M=\max_{\goc}(|D\gf(x)|+\vno{D^2\gf(x)})$ 
and, if
$x\in\bgo$,
$$
\aligned
B(x,u\sep(x),p)&=B(x,u\sep(x),p+\gep
D\gf(x))-\gep\ip{D\gf(x),n(x)}\\
& \le B(x,u(x),p+\gep
D\gf(x))-\gep.\endaligned$$ 
Using next that $\jtpcgo u\sep(x)+(\gep
D\gf(x),\gep D^2\gf(x))=\jtpcgo u(x)$ (see Remarks \rebx),
we conclude that if $C=\go(\gep M)/\gamma$, then $u\sep$ 
is a
subsolution of (BVP) with $B+\gep$ in place of $B$. 
Similarly, $v\sep$ is a supersolution of (BVP) with $B$
replaced by $B-\gep$.

We are reduced to proving that if $u$ is a subsolution and
$v$ is a supersolution of (BVP) with $B$ replaced by
$B+\gep$ and $B-\gep$ respectively, then $u\le v$.   Assume,
to the contrary, that $\max_{\goc}(u-v)>0$.  We know that
$\max_{\bgo}(u-v)\ge\max_{\goc}(u-v)$ by \S \variations.B,  
so there must exist $z\in\bgo$ such that
$u(z)-v(z)=\max_{\goc}(u-v)>0$. Put 
$$\gF(x,y)=u(x)-v(y)-\ot\ga\nos{x-y}+
f(z,u(z))\ip{n(z),x-y}-|x-z|^4\quad\text{on 
}\goc\times\goc,$$
 and let $(x\sga,y\sga)$ be a maximum point of
$\gF$.  By Proposition \pra,  since $u(x)-v(x)-|x-z|^4$ has
$x=z$ as a unique maximum point,  $$x\sga\rtra z,\qquad
\ga\nos{x\sga-y\sga}\rtra 0,\qquad u(x\sga)\rtra 
u(z),\qquad v(x\sga)\rtra
v(z)\eq$$  
\re\eqbcs 
\noindent as $\ga\rtra\infty$.  For simplicity, we
now write $\ph$ for $(x\sga,y\sga)$ and put 
$$\psi(x,y)=\ot\ga\nos{x-y}-f(z,u(z))\ip{n(z),x-y}+
|x-z|^4.$$ 
We have not yet invoked the exterior sphere condition
\eqbck) and do so now.  Restating \eqbck) as 
$$\nos{y-x-rn(x)}>r^2\quad\text{for }x\in\bgo\text{ and 
}y\in\gO,$$ we
see it is equivalent to
$$\ip{n(x),y-x}<\oo{2r}\nos{y-x}\quad\text{for 
}x\in\bgo\text{ and }y\in\gO.\eq$$
\re\eqbct 
\noindent Using \eqbct) we compute that if $\hx\in\bgo$ then
$$\aligned
B(\hx,u(\hx),D_x\psi\ph)&=B(\hx,u(\hx),\ga(\hx-\hy)\\
&\hphantom{=}-f(z,u(z))n(z)+
4|\hx-z|^2(\hx-z))\\
& \ge-{\ga\over
2r}|\hx-\hy|^2-f(z,u(z))\ip{n(\hx),n(z)}\\
&\hphantom{=}+O(|\hx-z|^3)+f(\hx,u(\hx))\endaligned$$
and in view of \eqbcs) this implies
$$B(\hx,u(\hx),D_x\psi\ph)\ge
o(1)\quad\text{as }\ga\rtra\infty\text{ if }\hx\in\bgo.$$ 
A similar
computation (this time involving \eqbcn)) shows that 
$$B(\hy,v(\hy),-D_y\psi\ph)\le
o(1)\quad\text{as }\ga\rtra\infty\text{ if }\hy\in\bgo.$$ 
Recalling that
$u$ is a subsolution with respect to $B+\gep$, etc., and the
definition of subsolutions and supersolutions, we conclude 
$$\cases F(\hx,u(\hx),D_x\psi\ph,X)\le
0\ \,\iq{for}(D_x\psi\ph,X)\in\overline\jtpcgo u(\hx),\\
F(\hy,v(\hy),-D_y\psi\ph,Y)\ge
0\iq{for}(-D_y\psi\ph,Y)\in\overline\jtmcgo 
v(\hy),\endcases\eq
$$
\re\eqbcux 
\noindent provided $\ga$ is large.

The remaining arguments are very close to those given in the
proof of Theorem \thb.   We apply Theorem \tha\  with $k=2$,
$\co_1=\co_2=\goc$, $u_1=u$, $u_2=-v$, and $\gep=1/\ga$ to
find $X,Y\in\m$ such that  $$(D_x\psi\ph,X)\in\overline J
^{2,+}_{\overline\Omega},\qquad
(-D_y\psi\ph,Y)\in\overline J^{2,-}_{\overline\Omega} 
v(\hy)$$ and $$-(\ga+\vno
A)\pmatrix I&0\cr0&I\endpmatrix\le \pmatrix 
X&0\cr0&-Y\endpmatrix\le A+\oo\ga
A^2,\eq$$ \re\eqbcz 
\noindent where $A= D^2\psi\ph$.  Observing that
$$A=\ga\pmatrix I&-I\cr-I&I\endpmatrix+O(\nos{\hx-z}),\ \
A^2=2\ga^2\pmatrix I&-I\cr-I&I\endpmatrix
+O(\ga\nos{\hx-z}+|\hx-z|^4),$$ we see from \eqbcz) that 
$$-3\ga\pmatrix I&0\cr0&I\endpmatrix+o(1)\le 
\pmatrix X&0\cr0&-Y\endpmatrix\le 
3\ga\pmatrix I&-I\cr-I&I\endpmatrix+
o(1)\iq{as}\ga\rtra\infty.$$ We want
to use this information in conjunction with \eqbcux); to
this end, we formulate the next lemma (which will be used
later as well). \rpre\lebcb \myProclaim Lemma \li.  Let
{\rm \eqbcl)} and {\rm \eqbcm)} hold.  Let $\gd>0$, 
$x,y\in V$,
$r\in\R$, $p\in\RN$, $X,Y\in\m$, and 
$$-(3\ga+\gd)\pmatrix I&0\cr0&I\endpmatrix\le \pmatrix
X&0\cr0&-Y\endpmatrix\le 3\ga\pmatrix
I&-I\cr-I&I\endpmatrix+\gd\pmatrix I&0\cr0&I\endpmatrix$$ hold.  Then
$$F(y,r,p,Y)-F(x,r,p,X)\le \go((\ga+\tfrac 23
\gd)\nos{x-y}+|x-y|(|p|+1))+2\go(\gd).$$

\demo{Proof} By assumption, $$-(3\ga+2\gd)\pmatrix 
I&0\cr0&I\endpmatrix\le
\pmatrix X-\gd I&0\cr0&-Y-\gd I\endpmatrix\le 
3\ga\pmatrix I&-I\cr-I&I\endpmatrix.$$ Using \eqbcm)  we 
therefore have
$$F(y,r,p,Y+\gd I)-F(x,r,p,X-\gd I)\le \go((\ga+{\tfrac 23}
\gd)\nos{x-y}+|x-y|(|p|+1)).$$ Next, \eqbcl) implies
$$F(y,r,p,Y+\gd I)-F(x,r,p,X-\gd I)\ge F(y,r,p,Y)-
F(x,r,p,X)  -2\go(\gd)$$ and the result follows.

We complete the proof of comparison.  From the lemma and the
developments above,  
$$
\aligned
0&\le
F(\hy,v(\hy),-D_y\psi\ph,Y)- F(\hx,u(\hx),D_x\psi\ph,X)\\
&\le  F(\hy,u(\hx),-D_y\psi\ph,Y)-
F(\hx,u(\hx),-D_y\psi\ph,X)\\
&\hphantom{\leq}-\gamma(u(\hx)-v(\hy))
+\go(4|\hx-z|^3)\\
&\le
\go(\ga\nos{\hx-\hy}+|\hx-\hy|(|D_y\psi\ph|+1))-\gamma
(u(\hx)-v(\hy))+o(1)\endaligned$$
as $\ga\rtra\infty$.  From this we obtain a contradiction as
in the proof of Theorem \thb. \enddemo

\demo{Proof of existence}  In order to establish
existence, we may use Perron's Method from Theorem \thp\  as
supplemented by Remark \reb\ to deduce that it suffices to
produce a  subsolution and a supersolution.  Using
Proposition \prbca, we see it suffices to produce a
classical subsolution and supersolution.  
    Let $\gf\in
C^2(\goc)$ be as in Lemma \lebca.  Put 
$$\uu(x)=-A\gf(x)-B,\qquad \ou(x)=A\gf(x)+B$$ where $A,\ 
B$ are
constants to be determined.  Observing that  
$$B(x,\ou(x),D\ou(x))\ge A\ip{n(x),D\gf(x)}+f(x,0)\ge
A+f(x,0)$$  for $x\in\bgo$ and
$$F(x,\ou(x),D\ou(x),D^2\ou(x))\ge
F(x,0,AD\gf(x),AD^2\gf(x))+\gamma B$$  for $x\in\goc$  we 
see
that if we first choose $A$ and then $B$ by
$$A=\sup_{\bgo}|f(x,0)|,\qquad B=\oo \gamma\sup_{\goc,|p|+
\vno X\le
AM}| F(x,0,p,X|,$$ where
$M=\sup_{\goc}(|D\gf(x)+\vno{D^2\gf(x)})$, then $\ou$ is a
supersolution.   Likewise, $\uu$ will be a subsolution. It
remains to prove Lemma \lebca.
\enddemo

\demo{Proof of Lemma {\rm \lebca}}  In view of \eqbcjx),
there is a $C^1$ function $\psi$ on $\RN$ such that
$\psi\ge0$ on $\RN\backslash\gO$, $\psi <0$ on $\gO$, and
$|D\psi|>0$ on $\bgo$.  Multiplying by a constant if
necessary, we may also assume that
$|D\psi(x)|>\ip{\gn(x),n(x)}^{-1}$ for $x\in\bgo$.  By
standard approximation theorems, we see that there is a
sequence of $C^2$ functions $\psi_k$ such that
$\psi_k\rtra\psi$ in $C^1(\goc)$ as $k\rtra\infty$.  In
particular, we have $$\ip{\gn(x),D\psi_k(x)}\rtra
\ip{\gn(x),D\psi(x)}\iq{uniformly on}\bgo\text{ as }k\rtra
\infty.$$ 
Since $n(x)=D\psi(x)/|D\psi(x)|$ and
$|D\psi(x)|\ip{\gn(x),n(x)}>1$ for $x\in\bgo$, if $k$ and
then $C$ are sufficiently large, $\gf=\psi_k+C$ has the
desired properties.\enddemo  

\subheading{\num{\secno.C.}
The generalized Dirichlet problem} 
We turn to the Dirichlet problem
$$\cases F(x,u,Du,D^2u)=0&\text{in $\gO$,}\\ u-f(x)=0&
 \text{on $\bgo$}\endcases\eq$$\re\eqbctx where $f\in 
C(\bgo)$.

The simplest example of this problem with a unique solution
that does not satisfy boundary conditions in the strong
sense is the following: 
  $F(x,u,Du,D^2u)
=u-h(x)$ and $f(x)\equiv 0$.  If $h$ is
continuous on $\goc$, then $u\equiv h$ is evidently the only
solution of this problem and it does not satisfy $u=0$ on
$\bgo$ in the strong sense unless $h=0$ on $\bgo$.

We will provide a comparison theorem, Theorem \thbcb
\ below, with this example as a special case.  This
comparison theorem, however will be different in character 
from those
we have obtained so far in that it will assert comparison
for a {\it continuous\/} subsolution and supersolution.   As
we know, the ability to compare semicontinuous subsolutions
and supersolutions allows us to prove existence.  Thus the
following counterexample to existence indicates that the 
restriction to continuous functions (or some other
assumption) is necessary in Theorem \thbcb. 
 
\rpre\exbcb 
\ex{Example\ \li}   Let
$N=2$ and $\gO=\{(x,y): -1<x<1, \ 0<y<1\}$.  Let $f\in
C(\bgo)$ satisfy $0\le f\le 1$, $f(x,0)=0$, and $f(x,1)=1$ 
for
$|x|<1$.  We claim that the problem 
$$\cases u+xu_y=0&\text{in
$\gO$,}\cr u=f&\text{on $\bgo$}\endcases\eq$$ 
\re\eqbcu 
\noindent does not have a
solution.  To see this, suppose the  contrary and let
$u(x,y)$ be a solution.  We claim that then for fixed
$|x|<1$, $v(y)=u(x,y)$ is a solution of 
$$v+xv'=0\text{ in }(0,1),\qquad v(0)=0,\qquad 
v(1)=1.\eq$$ \re\eqbcv 
\noindent Assume
for the moment that this is the case.  Then we remark that 
Theorem \thbcb\ below shows that $v$ is uniquely determined
by \eqbcv), $v\equiv 0$ is a solution of \eqbcv) if $x>0$
and $v(y)=e^{(1-y)/x}$ is a solution of \eqbcv) if $x<0$
(these last assertions the reader can check by
computations).   Thus we see that $u(x,y)=0$ if $x>0$ and
$u(x,y)=e^{(1-y)/x}$ if $x<0$; however,  $u$ cannot then be
continuous at $(0,1)$ and we conclude that \eqbcu) does not
have a solution.
\endex

 It remains to
verify that if $u$ is a solution of \eqbcu), then $v$ is a
solution of \eqbcv).  One way to check this is as follows:
Fix $\ox\in(-1,1)$ and set $v(y)=u(\overline x,y)$.  Let 
$\gf\in
C^2([0,1])$ and assume that $v-\gf$ has a unique maximum at
some point $\oy\in[0,1]$.  For $\ga>0$ let $\ph$ be a
maximum point of the function $u\p 
-\gf(y)-\ga\nos{x-\overline x}$. 
It follows from Proposition \pra\ that then $\hx\rtra\ox$,
$\hy\rtra\oy$ as $\ga\rtra\infty$; moreover, we have
$u\ph+\hx\gf'(\hy)\le 0$.  From this information we deduce
that $v(\oy)+\ox\gf'(\oy)\le 0$ and conclude that $v$ is a
subsolution of \eqbcv) with $x=\ox$.  In a similar way, one
shows that $v$ is a supersolution.  Finally we remark that 
$\ou\equiv
1$, $\uu\equiv 0$ are, respectively, a supersolution and a
subsolution of \eqbcu).  Thus comparison of semicontinuous
semisolutions of \eqbcu) must fail.  (If the reader is
concerned about regularity questions and the ``corners,"
note that the corners can be smoothed up without changing
the essential points of the discussion.)

We turn to the comparison result for continuous solutions.
\rpre\thbcb \myProclaim Theorem \li.    Let {\rm \eqbcjx),
\eqbckx), \eqbcl)}, and {\rm \eqbcm)} hold.  If $u,\ v\in 
C(\goc)$,
$u$ is a subsolution of {\rm \eqbctx)} and $v$ is a 
supersolution
of {\rm \eqbctx)}, then $u\le v$.

\demo{Proof} As before, we argue by contradiction and so 
suppose that
$\max_{\goc}(u-v)>$$0$.  We may assume that
$\max_{\goc}(u-v)=u(z)-v(z)>0$ for some $z\in\bgo$.  We
divide our considerations into two cases.

First, we consider the case when $v(z)<f(z)$.  For $\ga>1$
and $0<\gep<1$, we define the function 
$$\gF\p=u(x)-v(y)-|\ga(x-y)+\gep
n(z)|^2-\gep|y-z|^2\iq{on}\goc\times\goc.$$ Let $\ph$ be a
maximum point of $\gF$.  We assume $\ga$ is so large that
$z-(\gep/\ga) n(z)\in \gO$.  The inequality
$\gF\ph\ge\gF(z-(\gep/\ga) n(z),z)$ reads
$$|\ga(\hx-\hy)+\gep n(z)|^2+\gep|\hy-z|^2\le
u(\hx)-v(\hy)-u\left(z-{\gep\over \ga}n(z)\right)+v(z)$$ 
from which it
follows  (using the continuity of $u$!) that  if $\gep$ is
fixed,  $\hx,\hy\rtra z$ and $\ga(\hx-\hy)+\gep n(z)\rtra 0$
as $\ga\rtra\infty$.   Indeed, it is clear that $\ga(\hx 
-\hy)$
remains bounded as $\ga\rtra\infty$, so assuming that
$\ga(\hx-\hy)\rtra w$, $\hx,\hy\rtra \tilde z$ (along a
subsequence) as $\ga\rtra\infty$, one finds  $$|w+\gep
n(z)|^2+\gep|\tilde z-z|^2\le u(\tilde z)-v(\tilde
z)-u(z)+v(z)\le 0,$$ where the last inequality follows from
the definition of $z$, whence the claim is immediate.   
From this we see that  $\hx=\hy -(\gep n(z)+o(1))/\ga$ as
$\ga\rtra\infty$ and hence that $\hx\in\gO$ if $\ga$ is 
large
enough.  Since $v(z)<f(z)$, we have (using the continuity of
$v)$ $v(\hy)<f(\hy)$  if $\ga$ is large enough.  Thus, if
$\ga$ is large enough, we have  
$$F(\hx,u(\hx),p,X)\le
0\iq{for} (p,X)\in\jtpcgo u(\hx)$$  and
$$F(\hy,v(\hy),q,Y)\ge 0\iq{for} (q,Y)\in\jtmcgo v(\hy).$$ 
Observe that if we set $\gf(x,y)=|\ga(x-y)+\gep
n(z)|^2+\gep\nos{y-z}$, then
$$
D_x\gf(x,y)=2\ga(\ga(x-y)+\gep n(z)),
$$
$$
-D_y\gf(x,y)=2\ga(\ga(x-y)+\gep n(z))-2\gep (y-z),
$$ 
and 
$$
D^2\gf(x,y)=2\ga^2\pmatrix \ I&-I\cr-I&\ I\endpmatrix
+2\gep\pmatrix 0&0\cr0&I\endpmatrix.$$
 Now, using Theorem \tha\ together with Lemma \lebcb,
calculating as usual, sending $\ga\rtra\infty$ and then
$\gep\downa0$, we obtain a contradiction.

It remains to treat the case when $v(z)\ge f(z)$.  Since
$u(z)>v(z)$ this entails $u(z)>f(z)$.  Replacing the above
$\gF$ by  $$\gF\p=u(x)-v(y)-|\ga(x-y)-\gep
n(z)|^2-\gep|x-z|^2,$$ we argue as above and obtain a
contradiction.\enddemo

\subheading{\num{\secno.C\<$'$\<.} The state
constraints problem}
The problem  
\pre\eqbcaa 
$$\cases F(x,u,Du,D^2u)\le 0&\text{in\ \ $\Omega,$}
\cr F(x,u,Du,D^2u)\ge 0&\text{in\ \ $\goc$}\endcases\eq$$  
corresponds to an important
problem in optimal control, the so-called state 
constraints problem provided $F$ has
the  form of the left member of \eqxh).  This is, indeed, 
the extreme form of
\eqbctx) where $f(x)\equiv-\infty$.  To make this clear, we
give the definitions precisely.  A function $u\in 
\roman{USC}(\goc)$
(respectively, $v\in \roman{LSC}(\goc))$ is called a 
subsolution
(supersolution) of \eqbcaa) if  $$F(x,u(x),p,X)\le 0\iq{for}
x\in\goc\text{ and }(p,X)\in\overline J^{2,+} u(x)$$ 
(respectively,
$$F(x,v(x),p,X)\ge 0\iq{for} x\in\gO\text{ and 
}(p,X)\in\overline J
_{\overline\Omega}^{2,-}
v(x)).$$ An application of Theorem \thbcb\ immediately
yields \rpre\thbcd \myProclaim Theorem \li.  
Let {\rm \eqbcjx), \eqbckx), \eqbcl)}, and {\rm \eqbcm)} 
hold.  If $u,
v\in C(\goc)$, $u$ is a subsolution of {\rm \eqbcaa)} and 
$v$ is a
supersolution of {\rm \eqbcaa)}, then $u\le v$ on $\goc$.

Indeed, $v$ is a supersolution and $u$ is a subsolution of
\eqbctx) with $f(x)\geq u(x)$ on $\bgo$.  

In fact, a subsolution $v$ of \eqbcaa) is a supersolution
of \eqbctx) for any $f\in C(\bgo)$ and hence a ``universal
bound" on all solutions.  This implies that solutions of
\eqbcaa) do not exist in general.  For example, if $F$ is a
linear  uniformly elliptic operator, \eqbctx) is uniquely
solvable for arbitrary $f\in C(\bgo)$ and then there can be
no universal bound that is continuous in  $\goc$.  On the
other hand, if $F$ is suitably degenerate in the normal
direction at $\bgo$, one can sometimes find  such a 
supersolution.  

\subheading{\num{\secno.D.}   A remark {\rm \bc} in the 
classical sense} 
The reader will have noticed by now that we gave some 
examples of solutions of  (BC) that were smooth but did
not satisfy the boundary condition in the classical sense. 
In these examples, the equation (E) is rather degenerate. 
In fact, ``degeneracy" of a sort is necessary for this to
happen,  as we now establish.    \rpre\prbcq \myProclaim
Proposition \li.  Let $B\in C(\bgo\times\R\times\RN)$,  $F$
be \good, and  $u\in C^2(\goc)$ be a
subsolution {\rm (}respectively, supersolution{\rm )} of 
{\rm (BVP)}.  If 
$$
\cases
\limsup_{\gm\rtra\infty}F(x,r,p,X-\gm
n\otimes n)>0\\
\text{for }(x,r,p,X)\in
\bgo\times\R\times\RN\times\m,n\in\RN\backslash\{0\}%
\endcases\eq$$
\re\eqbcw 
\noindent {\rm (}respectively,
$$
\cases\liminf_{\gm\rtra\infty}F(x,r,p,X+\gm
n\otimes n)<0\\
\text{for }(x,r,p,X)\in
\bgo\times\R\times\RN\times\m,n\in\RN\backslash\{0\}),%
\endcases\eq$$
\re\eqbcwx
\noindent then $u$ is a classical subsolution {\rm 
(}respectively,
supersolution{\rm )} of {\rm \bc}.

\demo{Proof} It will suffice to treat the subsolution 
case.  We rely
again on Remarks \rebx (iii).  Let $y\in \RN\backslash\goc$
and $z\in\bgo$ be a nearest point to $y$ in $\goc$.  It is
clear that the set of such nearest points (as $y$ varies)
is dense in $\bgo$.  Then $\goc\subset\{ x\in\RN:|x-y|\ge
|y-z|\}=\co$ where the last equality is the definition of
$\co$; now this inclusion evidently implies $\jtpco
{\text{Zero}}\subset \jtpcgo {\text{Zero}}(z)$.  On the 
other hand, we computed $
\jtpco {\text{Zero}}(z)$ in Remark \rebx---with the 
notations $n= (y-z)/|y-z|$ and $r=|y-z|$,
 this computation
shows that  
$$\left(-\gl n,{\gl\over r}I+\gm n\otimes 
n\right)\in\jtpco {\text{Zero}}(z)$$
for $\gl>0$  and $\gm\in \R$.  Thus the assumption
that $u$ is a supersolution implies that  
$$B(z,u(z),Du(z)-\gl n)\wedge
F\left(z,u(z),Du(z)-\gl n,D^2u(z)+{\gl\over r}I+\gm 
n\otimes n
\right)\le 0 $$
for $\gl>0$ and $\gm\in\R$. Taking the limit superior as 
$\gm\rtra-\infty$ 
and   using the fact that assumption \eqbcw)  implies 
$B(z,u(z),Du(z)-\gl n)\le 0$ for 
$\gl>0$ and then    
  letting $\gl\downa0$ we find  $B(z,u(z),Du(z))\le 0$. 
Since $z$ was an arbitrary
``nearest point" and these are dense in $\bgo$,  the proof
is complete.  The reader will notice that it suffices to
choose $n$ in \eqbcw) and \eqbcwx) from 
$N_{\goc}(x)$.\enddemo\medskip

\subheading{\num{\secno.E.}   Fully nonlinear boundary 
conditions} 
We conclude this section by describing an extremely general
existence and uniqueness result  for a fully nonlinear 
first-order 
boundary operator of the form  $B=B(x,p),\ B\in
C(\bgo\times\RN)$.  We will require that $B$ satisfies 
\pre\eqbcaax 
$$ 
\ctrblock
$ B$ is  uniformly continuous
in $p$ uniformly in $x\in\bgo$\endblock\eq$$ and  
\pre\eqbcab
$$|B(x,p)-B(y,p)|\le \go(|x-y|(1+|p|)) \eq$$ for some
$\go\:[0,\infty)\rtra[0,\infty)$ satisfying $\go(0+)=0$
and   for some $\gn>0$  \pre\eqbcac $$B(x,p+\gl n(x))\ge
B(x,p)+\gn \gl \iq{for} x\in\bgo,\ \gl\ge 0,
p\in\RN.\eq$$
 We will also need to strengthen a bit the regularity of
$\gO$ by assuming that  \pre\eqbcad $$\bgo\text{ is of 
class }
C^{1,1}.\eq$$ Observe that \eqbcad) is stronger  than
\eqbcjx) and \eqbck) since $\gO$ is bounded.  The loss of
generality,  however, is not really more  restrictive for
applications, since $\gO$ typically is either smooth  or has
``corners,'' and corners are not allowed by  \eqbcjx) in any
case.   (Corners typically require some  ad hoc analysis.)

We have

\rpre\thbcc  \myProclaim  Theorem \li.  Let {\rm 
\eqbcaax)--\eqbcad)  }
and {\rm \eqbckx)--\eqbcm)} hold.  If $u$ is a
subsolution of  {\rm (BVP)} and $v$ is a supersolution of 
{\rm (BVP)},
then $u\le v$ in $\goc$.  Moreover, {\rm (BVP)} has a unique
solution.

\rpre\exbcc 
\ex{Example\  \li}   The
Neumann boundary conditions  of \eqbcj) are special cases 
of those considered above; note that  \eqbcab) holds since
$n$ is Lipschitz continuous  on $\bgo$ because of \eqbcad). 
Similarly,  this setting includes the  more general linear
oblique derivative boundary condition corresponding to 
\pre\eqbcae $$B(x,p)= \ip{\gamma(x), p}+f(x)\eq$$ where 
$\gamma$
is Lipschitz and satisfies  \pre\eqbcaf
$$\ip{\gamma(x),n(x)}\ge\gn>0\iq{for}x\in\bgo;\eq$$ note 
that
$f$ is merely required to be continuous on $\bgo$. 
\endex

Nonlinear examples arising in optimal control and
differential games  are typically written in the form 
\pre\eqbcag $$B(x,p)= \sup_{\ga}(\ip{\gamma_\ga(x),
p}+f_\ga(x))\eq$$ and  \pre\eqbcah $$B(x,p)=
\inf_{\gb}\sup_{\ga}(\ip{\gamma_{\ga,\gb}(x),
p}+f_{\ga,\gb}(x))\eq$$ where $\gamma_\ga,\ f_\ga$,
$\gamma_{\ga,\gb},\ f_{\ga,\gb}$ are functions that  satisfy
the conditions (including continuity) laid on $\gamma,  f$
above uniformly in the parameters $\ga, \gb$.

As a  last example, we offer the  famous capillarity
boundary condition  that  corresponds to  \pre\eqbcabi $$
B(x,p)=\ip{n(x),p}-\gk(x)(1+|p|^2)^{1/2}\eq$$ where $\gk\in
C(\bgo)$ and $\sup_{\bgo}|\gk|<1$.  

We will not prove Theorem \thbcc\ since its proof, while
based on  the strategies presented above, requires some
rather complex  adaptations  of them.   We remark that the
form of the boundary operator $B(x,p)$ could be  generalized
to  $B(x,r,p)$ above, provided that $B$ is  nondecreasing
with respect to $r\in\R$ and \eqbcaax)--\eqbcac) hold
uniformly  in bounded $r$.

 \subheading{Notes on \S7} Boundary
conditions in the viscosity sense first appeared in P. L. 
Lions
[\pllneumann] and B. Perthame and  R. Sanders [\persan] 
for Neumann
or oblique derivative problems for first order equations, 
making
this problem  well posed!   A general formulation was 
given later
by H. Ishii [\isbvp].  Existence and uniqueness was shown  
in P. L.
Lions [\pllneumann] (for the first-order case), H. Ishii 
and P. L.
Lions [\lii] (second-order case), and
 P. Dupuis and H. Ishii [\dupish, \dupisha]
(second-order and nonsmooth domains).  However, it is 
 worth emphasizing that the presentation we give using 
Theorem \tha\
\ simplifies the arguments substantially (even in this 
special case).

 Concerning the Dirichlet problem, most of the works 
required (and
often still require) continuous solutions 
 up to the boundary and prescribed data on the entire 
boundary.  It
was first noted in P. L. Lions [\plla], G. Barles 
 [\barlesa, \barlesb] and H. Ishii [\isp] that this can be 
 achieved for special classes of equations by imposing 
compatibility
conditions on 
 the boundary data or, in general, by assuming the 
existence of
appropriate super and subsolutions.  It was also 
 noted in M. G. Crandall and R. Newcomb [\cn] and P. E. 
Souganidis
[\takb] that some part of the boundary may be 
 irrelevant (in the case of a first-order equation).  

 Progress on the understanding of Dirichlet boundary 
conditions was
stimulated by the ``state constraints"
 problems studied first by M. Soner [\sona, \sonb] and 
later by I.
Capuzzo-Dolcetta and P. L. Lions [\capl].  
 This led to the ``true" viscosity formulation of Dirichlet
conditions as considered in H. Ishii [\isbvp] and
 G. Barles and B. Perthame [\bapa, \bapb].  
 This formulation, in some sense, automatically selects 
the relevant
part of 
  the boundary for degenerate problems and yields 
uniqueness for
continuous solutions.  However, the existence of a 
  continuous solution fails in general (Example \exbcb), 
and the
situation is not entirely clear except for first-order
  optimal control problems [\bapa, \bapb, \isbvp].

  We also mention at this stage that some uniqueness 
results for
semicontinuous solutions have begun 
  to emerge (G. Barles and B. Perthame [\bapa, \bapb] and 
E. N.
Barron and R. Jensen [\barjena, \barjenb])
  for optimal control problems 
  and the study of some second-order nondegenerate
state constraints problems (J. M. Lasry and P. L. Lions 
[\lalib]).

  Finally, the possibility of solving general equations 
with general
fully nonlinear oblique derivative type boundary 
  conditions---a rather startling fact---was illustrated 
in G. Barles
and P. L. Lions [\barlespll]
   for first-order equations.  The full second-order 
   result of Theorem \thbcc\ is a generalization by 
   G. Barles [\barlesg]  of a similar result of 
   H. Ishii [\ishiob].

\medskip \resultno=0 \equationno=0 \spre\parabolic
\heading \sli. Parabolic
problems\endheading 
\def\secno{\the\sectionno}

In this section we indicate how to extend the results of the
preceeding sections to problems involving the parabolic
equation $$u_t+F(t,x,u,Du,D^2u)=0\leqno{\roman{(PE)}}$$ 
where now $u$
is to be a function of $(t,x)$ and $Du,D^2u$ mean
$D_xu(t,x)$ and $D^2_xu(t,x)$.   We do this by discussing
comparison for the Cauchy-Dirichlet problem on a bounded
domain; it will then be clear how to modify other proofs as
well.  Let $\co$ be a locally compact subset of $\RN$,
$T>0$, and 
 $\co_T=(0,T)\times\co$.  We denote by $\dppco,\ \dpmco$ the
``parabolic" variants of the  semijets $\jtpco, \jtmco$;
for example,  if $u\:\cy\rtra\R$ then  $\dppco u$ is defined
by $(a,p,X)\in \R\times\ \RN\times\m$ lies in  $\dppco u(s,
z)$ if $(s,z)\in\cy$ and  $$\eqalign{u(t,x)\le & u(s,z) +
a(t- s) + \ip{p,x- z}+\tfrac 12 \ip{X(x-z),x-z}\cr
&+o(|t-s|+|x-z|^2)\quad\text{as 
}\cy\ni(t,x)\rtra(s,z);}\eq$$ 
\re\eqpa 
\noindent similarly, $\dpmco u=-\dppco (-u)$.  The
corresponding definitions of $\dppcoc,$ $\dpmcoc$ are then
clear.   

 A {\it subsolution} of (PE) on $\cy$ is a function
$u\in \roman{USC}(\cy)$ such that  $$a+
F(t,x,u(t,x),p,X)\le 0\ti{for}
(t,x)\in\cy \ti{and} (a,p,X)\in\dppco u(t,x);\eq$$ \re\eqpb
\noindent likewise, a {\it supersolution} is a function 
$v\in 
\roman{LSC}(\cy)$
such that  $$a+F(t,x,v(t,x),p,X)\ge 0\ti{for} (t,x)\in\cy
\ti{and} (a,p,X)\in\dpmco v(t,x);\eq$$ \re\eqpc 
\noindent and a {\it
solution} is a function that is simultaneously a
subsolution and a supersolution.   

\rpre\repa
\rem{Remark\ {\rm \li}}  Suppose we set
$x_0=t$ and $\tilde x=(x_0,x)$.  Then (PE) is an equation of
the form $\hat F(\tilde x,u,D_{\tilde x}u,D_{\tilde
x}^2u)=u_{x_0}+\cdots$ that is \good\ if and
only if $F(t,x,u,Du,D^2u)$ is \good\  when $t$
is held fixed.   The definition of subsolutions, etc., of
(PE), which takes into account that only the first
derivative with respect to $t$ appears, does not coincide
with the definition of subsolutions, etc., of $\hat F=0$
under this correspondence.   However, it is not hard to see
that the two notions are equivalent.
\endrem

We will illustrate  the additional considerations that
arise in dealing with the Cauchy-Dirichlet problem for (PE)
as opposed to the pure Dirichlet problem in the elliptic
case.  The problem of interest has the form 
 $$\left\{\eqalign{&\text{(E)}\ \ \ \ u_t+
F(t,x,u,Du,D^2u)=0\ti{in} 
 (0,T)\times\gO,\cr &\text{(BC)}
 \ \ \ u(t,x)=0\ti{for} 0\le t< T\ti{and}x\in\bgo, \cr 
&\text{(IC)}
 \ \ \  u(0,x)=\psi(x)\ti{for}x\in\goc,\cr}\right.\eq$$  
\re\eqpr
\noindent where $\gO\subset\RN$ is open and  $T>0$ and 
$\psi\in
C(\goc)$ are given. By a subsolution of \eqpr) on
$[0,T)\times\goc$ we mean a function $u\in 
\roman{USC}([0,T)\times
\goc)$ such that $u$ is a subsolution of (E), $u(t,x)\le 0$
for  $ 0\le t< T\ti{and}x\in\bgo$ and $u(0,x)\le\psi(x)$ for
$x\in\goc$--the appropriate notions of supersolutions and
solutions are then obvious.

\rpre\thpa   \myProclaim Theorem \li.  Let $\gO\subset\RN$ 
be
open and bounded.   Let $F\in
C([0,T]\times\goc\times\R\times\RN\times\m)$ be continuous,
\good, and satisfy {\rm \eqbl)} for each fixed
$t\in[0,T)$, with the same function $\go$.  If $u$ is a
subsolution of {\rm \eqpr)} and   $v$ is  a supersolution 
of {\rm \eqpr)},
then $u\le v$ on $[0,T)\times\gO$.

 To continue, we require the parabolic analogue of Theorem
\tha.  It takes the following form

\rpre\thpb  \myProclaim Theorem \li.  Let $u_i\in
\roman{USC}((0,T)\times\co_i)$ for $i=1,\ldots,k$ where 
$\coi$ is a
locally compact subset of $\R^{N_i}$. Let  $\gf$ be defined
on an open neighborhood of
$(0,T)\times\co_1\times\cdots\times\co_k$ and such that
$(t,x_1,\cdots,x_k)\rtra\gf(t,x_1,\cdots,x_k)$ is once
continuously differentiable in $t$ and twice continuously
differentiable in
$(x_1,\cdots,x_k)\in\co_1\times\cdots\times\co_k$.  
Suppose that 
$\hat t\in(0,T),\ \hx_i\in\coi$ for $i=1,\ldots,k$ and  
$$
\aligned
w(t,x_1,\dotsc,x_k)&\equiv u_1(t,x_1)+
\dotsb+u_k(t,x_k)-\gf(t,x_1,\dotsc,x_k    )\\
&\le
w(\hat t,\hx_1,\dotsc,\hx_k)\endaligned$$ for $0< t< T$  
and $x_i\in\co$. 
Assume, moreover,  that there is an $ r>0$ such that for
every $M>0$ there is a $C$ such that for $i=1,\ldots,k$ 
$$\eqalign{&b_i\le C \ti{whenever}
(b_i,q_i,X_i)\in\dppco u_i(t,x_i),\cr& |x_i-\hx_i|+|t-\hat 
t|\le 
r\text{ and\ }|u_i(t,x_i)|+|q_i|+\vno{X_i}\le
M.\cr}\eq$$  \re\eqpdx 
\noindent Then for each $\gep>0$ there
are $ X_i\in {\cal S}(N_i)$ such that
$$\left\{\eqalign{&\text{{\rm (i)}}\ \ \
(b_i,D_{x_i}\gf(\hat 
t,\hx_1,\dotsc,\hx_k),X_i)\in\dppc_{\co_i}
u_i(\hat t,\hx_i)  \ti{for}i=1,\ldots,k,\cr& \text{{\rm 
(ii)}}\ \
-\left(\oo\gep +\vno A\right)I\le  \diag\le A+\gep A^2,\cr
&\text{{\rm (iii)}}
\ b_1+\dotsb+b_k=\gf_t(\hat 
t,\hx_1,\dotsc,\hx_k),\cr}\right.\eq$$
\re\eqprx
\noindent where $A=(D^2_x\gf)(\hat t,\hx_1,\dotsc,\hx_k)$.

Observe  that the condition \eqpdx) is guaranteed by having
each $u_i$ be a subsolution of a parabolic equation.   

\demo{Proof of Theorem {\rm \thpa}}   We first observe
that for $ \gep>0$,  $\tilde u=u-\gep/(T-t)$ is also a
subsolution of \eqpr) and satisfies (PE) with a strict
inequality; in fact, $$\tilde u_t+F(t,x,\tilde u,D\tilde
u,D^2\tilde u)\le -{\gep\over (T-t)^2}.$$ Since $u\le v$
follows from $\tilde u\le v$ in the limit $\gep\downa0$, it
will simply suffice to prove the comparison under the
additional assumptions  $$\left\{\eqalign{&\text{(i)}\
u_t+F(t,x,u,Du,D^2u)\le -\varepsilon/T^2<0\ti{and}\cr 
&\text{(ii)}\
\lim_{t\ua T}u(t,x)=-\infty\ti{uniformly
on}\goc.\cr}\right.\eq$$ \re\eqqa

We will assume  $$(s,z)\in(0,T)\times\gO\quad\text{and}\quad
u(s,z)-v(s,z)=\gd>0\eq$$  \re\eqqb 
\noindent and then contradict this
assumption. We may assume that $u, -v$ are bounded  above. 
 Let  $(\hat t,\hx,\hy)$  be a maximum point of
$u(t,x)-v(t,y)-(\ga/2)\nos{x-y}$ over
$[0,T)\times\goc\times\goc$ where $\ga>0$;  such a maximum
exists in view of the assumed bound above on $u,\ -v$, the
compactness of $\goc$, and \eqqa)(ii).  The purpose of the
term  $(\ga/2)\nos{x-y}$ is as in the elliptic case.   Set
$$M\sa=u(\hat t,\hx)-v(\hat 
t,\hy)-\ot\ga\nos{\hx-\hy}.\eq$$   By
\eqqb), $M\sa\ge\gd$.  If $\hat t=0$, we  have $$0<\gd\le
M\sa\le\sup_{\goc\times\goc}\left(\psi(x)-\psi(y)-\ot\ga%
\nos{x-y}
\right);$$
however, the right-hand side above tends to zero as
$\ga\rtra\infty$ by Lemma \lea, so $\hat t>0$ if $\ga$ is 
large. 
Likewise, $\hx,\hy\in\gO$ if $\ga$ is large by $u\le v$ on
$[0,T)\times\bgo$.  Thus we may apply Theorem \thpb\  at
$(\hat t,\hx,\hy)$ to learn  that there are numbers $a, b$ 
and
$X, Y\in\m$ such that
 $$(a,\ga(\hx-\hy),X)\in \dppcoc u(\hat t,\hx),\qquad
(b,\ga(\hx-\hy),Y)\in \dpmcoc v(\hat t,\hy)$$ such that 
$$a-b=0\quad\text{and}\quad-3\ga\pmatrix I&0\cr 
0&I\endpmatrix
\le\pmatrix \ X&\ 0\cr\ 0&-Y\endpmatrix\le 
3\ga\pmatrix \ I&-I\cr-I&\ I\endpmatrix.\eq$$ \re\eqps The 
relations
$$\eqalign{a+F(\hat t,\hx,u(\hat 
t,\hx),\ga(\hx-\hy),X)&\le -c,\cr 
 b+F(\hat t,\hy,v(\hat t,\hy),\ga(\hx-\hy),Y)&\ge 0,\cr}$$
and \eqps) imply $$\align
c&\le 
F(\hat t,\hy,v(\chatt,\hy),\ga(\hx-\hy),Y)-
F(\chatt,\hx,u(\chatt,\hx),\ga(\hx-\hy),X)\\
&\le\go(\ga|\hx-\hy|+|\hx-\hy|)\endalign$$ which
leads to a contradiction as in the proof of Theorem \thb.

Let us mention a couple of other adaptations of results
above to parabolic problems.  Section \variations.B  may be
regarded as establishing continuity of solutions with
respect to boundary data and the equation itself.  In the
parabolic case, there is the initial data, the boundary 
data,
and the equation to consider.  In the above context,  we may
consider a solution of $u_t+F(t,x,Du,D^2u)\le 0$ in
$(0,T)\times\gO$ and a solution of $v_t+G(t,x,Dv,D^2v)\ge0$
in $(0,T)\times\gO$ with the continuity and boundedness
properties assumed in Theorem \thpa.   Suppose $g(t)\ge
(G(t,x,p,X)-F(t,x,p,X))^+$ where $g$ is continuous and 
$(u(t,x)-v(t,x))^+\le K_1$ for $(t,x)\in(0,T)\times\partial
\gO$, $(u(0,x)- v(0,x))^+\le K_2$ for $x\in\goc$.  Then the
function $w(x,t)=v(x,t)+\max(K_1,K_2)+\int_0^t\,g(s)\,ds$ 
is a
solution of $w_t+F(t,x,Dw,D^2w)\ge 0$.  In this way, if we
have comparison, we conclude that   $$u(t,x)\le
v(t,x)+\max\left(\sup_{(0,T)\times\partial\gO}(u-v)^+
,\sup_{\{0\}\times\goc}(u-v)^+\right)+
\int_0^tg(s)\,ds.$$ In many cases
we may put 
$$g(t)=\sup_{\RN\times\RN\times\m}(G(t,\cdot)-F(t,\cdot))^+
,$$
and a simple example is $G(t,x,p,X)=F(t,x,p,X)+f(t,x)$.
\enddemo

\subheading{Notes on \S 8} The presentation follows M. G. 
Crandall and H.
Ishii [\craish].  Let us also mention that it is well 
recognized
that most results concerning stationary equations have
straightforward parabolic analogues; indeed, the parabolic 
situation
is often better,  since $$u_t+F(x,u,Du,D^2u)=0$$ more or 
less
corresponds to $\gl u+F(x,u,Du,D^2u)=0$ with large $\gl>0$.
Moreover, the special linear dependence on $u_t$ in the 
parabolic
case allows one to let $F$ depend on $t$ in a merely 
measurable
manner; see, e.g., H. Ishii [\ishm], B. Perthame and P. L. 
Lions \  
[\pllper], N. Barron and R. Jensen [\barjenm], 
and  D. Nunziante [\nuna, \nunb].  Another special
feature of the parabolic case is that, owing to the special
structure, one can often allow rather singular initial 
data.  For
example,  infinite values may be allowed and 
semicontinuity may
suffice, depending on the situation  (see, e.g.,
M. G. Crandall, P.
L. Lions, and P. E. Souganidis [\cls] and E. N. Barron and 
R. Jensen
[\barjena, \barjenb]).

\medskip
\resultno=0 \equationno=0 \spre\singular \heading \sli. 
Singular equations:  
an example from geometry\endheading 
\def\secno{\the\sectionno}

Let $p\otimes q=\{p_iq_j\},$ the matrix with entries
$p_iq_j$.   It can be shown that if $\psi$ is a smooth
function and $D\psi$ does not vanish on the level set
$\gG=\{\psi=c\}$ and $u$ is a classical solution of the
Cauchy problem
 $$u_t-\trace\lf(\lf(I-\tenpn{Du}\rt)D^2u\rt)=0, \qquad 
u(0,x)=\psi(x)\eq$$
\re\eqsa 
\noindent on some strip $(0,T)\times\RN$,  then 
$\gG_t=\{u(t,\cdot)=c\}$ represents the result of  evolving 
$\gG$ according to its mean curvature to the time $t$, 
whence
there is geometrical interest in \eqsa).   Indeed, in less
regular situations where we  have viscosity solutions,
$\gG_t$ has been proposed as a definition of the  result of
evolving $\gG$ in this way.   The nonlinearity involved, 
which we hereafter denote $F(p,X)=-\trace((I-(p\otimes 
p)/|p|^2) X)$,
is degenerate elliptic on the set $p\not=0$ and undefined at
$p=0$.  Thus the results of the preceeding sections do not
apply immediately.  However, this is easy to remedy using
the special form of the equation.    The extensions of $F$
to $(0,X)$ given by
 $$\uf(p,X)=\cases F(p,X)&\text{if $p\not=0$,}\cr -2\vno X &
\text{if $p=0$,}\endcases \ \ \ \ \ \  \of(p,X)=\cases 
F(p,X)&
\text{if $p\not=0$,}\cr 2\vno X &\text{if 
$p=0$}\endcases\eq$$ 
\re\eqsb 
are lower semicontinuous and upper semicontinuous
respectively.  We define $u$ to be a subsolution
(respectively, supersolution) of $u_t+F(Du,D^2u)=0$  if it
is a subsolution of $u_t+\uf(Du,D^2u)=0$ (respectively, of
$u_t+\of(Du,D^2u)=0$) and a solution  if it is both a
subsolution and a  supersolution.   The reason that this
will succeed is roughly that, in the analysis, we will need
only to insert $X=0$ when we have to deal with $p=0$, so
what one does with $p=0$ is not important so long as it is
consistent.    To illustrate matters in a slightly simpler
setting, let us consider instead  the stationary problem  $$
u+F(Du,D^2u)-f(x)=0\quad\text{in }\RN,\eq$$ \re\eqsc 
 where we use the
corresponding definitions of subsolutions, etc. Theorem
\thda\ remains valid for the current $F$. \rpre\thsa
\myProclaim Theorem \li.   Let $f\in \roman{UC}(\RN)$.  
Then {\rm (\eqsc)}
has a unique solution $u\in \roman{UC}(\RN).$

We sketch the proof, which proceeds according to the
outline given in the proof of Theorem \thda\  with  slight
twists.  

\demo{Proof of comparison}  The comparison proof
is a slight modification of that of \S 5.D.  We begin
assuming that $u, v$ are a subsolution and a supersolution
of \eqsc) and  $u(x)-v(y)\le L(1+|x|+|y|)$ (which is \eqdp),
and then proceed as before, ending up with 
\pre\eqsd
$$
\eqalign{u(\hx) - v(\hy) \leq f(\hx) - f(\hy)
+ \of(\hat p - D\gb_R(\hy),- \hZ - D^2\gb_R(\hy)) \cr
-\uf(\hat p + D\gb_R(\hx),\hZ + D^2\gb_R(\hx))\cr}
 \eq$$
in 
place of \eqds).  This still implies a bound on
$u(x)-v(y)-2K|x-y|$ as before.  Then proceeding still
further with the proof, we  adapt slightly 
 and consider a maximum point $\ph$ of  
$$\gF(x,y)=u(x)-v(y)-(\ga|x-y|^4+\gep(\nos x+\nos y)),$$
which will exist by virtue of the bound already obtained. We
assume, without loss of generality for what follows, that
$\gF\ph\ge0$ so that for some $C$
$$\ga|\hx-\hy|^4+\gep(\nos{\hx}+\nos{\hy})\le 
u(\hx)-v(\hy)\le
2K|\hx-\hy|+C.\eq$$ \re\eqse Using Remark \reba, we then
have the existence of $X,Y\in\m$ such that
$$(\hp+2\gep\hx,X+2\gep I)\in\jtpc u(\hx),\qquad
(\hp-2\gep\hy,Y-2\gep I)\in\jtmc v(\hy)\eq$$  \re\eqsf 
where 
$$
\gathered
\hp=4\ga|\hx-\hy|^2(\hx-\hy),\\
 \vno X,\vno Y\le
C_1\ga|\hx-\hy|^2\quad\text{and}\quad X\le 
Y.\endgathered\eq$$ 
\re\eqsg 
\noindent Now \eqse)
implies that $\gep(\nos {\hx}+\nos{\hy})$ and
$\ga|\hx-\hy|^3$, and hence $\hp$,  are bounded 
independently
of  $\gep\le 1$ for fixed $\ga\ge1$ while $\ga|\hx-\hy|^4$ 
is also 
bounded independently of $\ga$.   Hence, $\gep\hx,
\gep\hy\rtra0$ as $\gep\downa0$. We have the  following 
analogue
of \eqdz)
$$
u(\hx)-v(\hy)\le (\of(\hp-2\gep\hy,Y-2\gep I)
-\uf(\hp+2\gep\hx,X+2\gep
I))+f(\hy)-f(\hx)$$ 
that leads, by
use of the estimates above, to 
$$\limsup_{\gep\downa0}(u(\hx)-v(\hy))\le \of (p,Y_0)-\uf
(p,X_0)+\gk(C/(\ga)^{1/4})$$ where $C$ is some constant
 and $(p, X_0, Y_0)$ is a limit point 
of $(\hp,  X,  Y)$ as $\gep\downa 0$.  If $p\not = 0$,  we 
are
done since 
 $$\of (p,Y_0)-\uf (p,X_0)= F(p,Y_0)-F(p,X_0)\le 0$$
 because $X_0\le Y_0$.  If $p=0$, we use the information
$4\ga|\hx-\hy|^2(\hx-\hy)\rtra p=0$  and \eqsg)  (recall 
$\ga$
is fixed)  to conclude that $X_0=Y_0=0$, and then $ \of
(p,Y_0)=\uf (p,X_0)=0$, and we are still done.
\enddemo

\demo{Proof that solutions lie in $\uc$} This is
exactly as in \S\variations.D.
\enddemo

\demo{Proof of existence}  A supersolution and
subsolution are produced exactly as in \S\variations.D. 
Perron's method still applies here, since $\uf$ is
lower semicontinuous and, we may use Lemma \lepb \ as is
with $F$ replaced by $\uf$ and the variant of Lemma \lepb
\ given in Remark \reb\ with $G_-=\uf$  and 
$G^+=\of$.
\enddemo

\subheading{Notes on \S 9} 
As shown in \S \equations, a number of equations arising 
from 
geometrical considerations present singularities at $p=0$. 
 The
fact  that this can easily be circumvented was shown 
independently
by  L. C. Evans and J. Spruck [\evsp] and Y. Chen, Y. 
Giga, and S.
Goto [\giga]. We also mention the work of H. M. Soner 
[\sond] on the
equation \eqsa) and  the papers by G. Barles [\barlesc] 
and S. Osher
and J. Sethian [\osse], which showed  how various 
geometrical
questions about ``moving fronts" could be  reduced to 
equations
that can be handled by viscosity theory. A general class 
of singular 
equations is treated in Y. Giga, S. Goto, H. Ishii, and M. 
H. Sato 
[\gigi] and M. H. Sato [\satoa]; these works establish 
existence 
and uniqueness as well as convexity properties of solutions.

\medskip

\spre\perspectives \heading \sli. Applications and 
perspectives \endheading 
\def\secno{\the\sectionno}  

In this section we list some applications of the theory of
viscosity solutions and indicate some of the promising
directions for development of the theory  in the next few
years.  We give some important references but they are not
exhaustive.

To begin this rather long list of applications, we recall
that  perhaps the main motivation for developing the theory
was its relevance for  the theories  of  Optimal Control and
Differential Games. Indeed, as is well known, in the theory
of optimal control of  ordinary differential equations or
stochastic differential equations (with complete
observations) or in the theory of zero sum, two player
deterministic or stochastic differential games, the Dynamic
Programming Principle (DPP for short)  states that the
associated value functions  should be characterized as the
solutions of associated partial differential equations. 
These equations are  called  Bellman  or  or
Hamilton-Jacobi-Bellman (HJB for short) equations in
control theory or Isaacs equations in differential games. 
The DPP was, however, heuristic and  proofs of it required
more regularity of the value functions than they usually
enjoy.
  The flexibility of the theory of viscosity solutions  has
completely filled this regularity  gap:  roughly speaking,
value functions are viscosity solutions and are uniquely
determined by this fact (via the uniqueness of viscosity
solutions). See, e.g.,   [\evsoug, \plla,  \pllc,
\pllv, \taka]. This basic theoretical fact allows a
spectacular simplification of the theory of  deterministic
differential games [\evsoug, \taka] and also provided 
the possibility of  creating sound mathematical foundations
for stochastic differential games [\flsougb].  See also   
[\ishiirep, \newc] in addition to [\bara, \barish].

The generalization of the definition of viscosity solutions 
to systems in diagonal form is rather straightforward and 
has applications to optimal control and differential games.
For these topics see, e.g., [\capev, \leneng, \ishps, 
\ishiikoike, \lenyama, \yam].   In the case of 
systems, a combination of viscosity solutions and weak 
solutions
 based on distribution theory may define a natural notion 
 of weak solutions [\gigi].

More generally, as usual, a better understanding of
existence-uniqueness issues for classes of equations leads
to a better understanding of more specific issues.  Typical
examples here are   perturbation questions, asymptotic
problems, and a more detailed solution of some specific
applications to Engineering or Finance problems
[\barjenfa, \zara, \zarb].  Also, the part of  the
theory concerned with  boundary conditions  has led to a
rather complete theory for problems with state-constraints
(at least for deterministic problems) that are enforced by
cost or boundary mechanisms [\sona, \sonb, \capl,
\isbvp,   \bapa, \bapb, \barjena,
\barjenb].  Let us also mention that an interesting link
between viscosity solutions and the other main argument of
Control Theory, 
 the Pontryagin principle, has been shown [\barjenpon,
\barbua].  Of course, last but not least, these results
have led to numerical approaches to Control or Differential
Games problems via the resolution of the HJB (or Isaacs)
equation.

Indeed, the viscosity solutions theory is intimately
connected with numerical analysis and scientific computing. 
First of all, it provides efficient tools to perform
convergence analyses (e.g., [\rclap, \barlessoug, 
\pllsoa, \takd]).  It also indicates how to build
discretization methods or schemes for other general boundary
conditions and in particular for classical boundary
conditions when working with rather degenerate equations
[\falcone, \roto].  

Another consequence of efficient existence, uniqueness,
approximation, and convergence results is the possibility of
establishing or discovering various qualitative properties
of solutions (formulae, representations, singularities, 
geometrical properties, characterizations and properties of
semigroups, $\dots$---[\bardiev, \cansona, \cansonb,
\cplld, \clinv, \evsoug, \ishiirep, \pla, \pllro]).   
Of course, one of
the most important qualitative properties is the regularity
of solutions.  Viscosity solutions, because of their
flexibility and their pointwise  definition, have led to
regularity results  that are spectacular  either in their
generality (regularizing effects, Lipschitz regularity, or
semicontinuity---[\pllreg,  \plla, \cls,
\barlesf,  \barlese, \barlesgeo]) or by their originality
($C^{1,\ga},\ C^{2,\ga}$ or pointwise  $L^{p}$ 
estimates---[\cafa, 
\lii, \trua, \trub, \truc, \jenr]).

The uniqueness and convergence parts of the theory have made
possible partial differential equations approaches to
various asymptotic problems like large deviations
[\evish,  \flsoa, \bardi, \pera],
geometrical optics [\evsougb, \barlesevsoug,
\barlesbronsoug, \bardiper,   \evsonsoug],  or
homogenization problems [\evpert,  \pll, \pllpap]  by
arguments that are both powerful and simple.

More specific applications concern the interplay between the
behavior of solutions at infinity   and structure conditions
[\ish, \cplld, \clinv, \pllsov] or the treatment of
integrodifferential operators [\lena, \sayaa].
 Other applications  concern some particular
classes of equations arising in Engineering like some models
of the propagation of fronts in Combustion Theory
[\barlesc, \osse], or the so-called shapes from shading
models in Vision Theory [\roto].

Finally, a large part of the theory (but not yet  all of it)
has been \lq\lq raised" to infinite-dimensional equations
both for first-order and second-order equations
[\clii, \pllcid--\plle]. In addition to \lq\lq
standard" extensions to infinite-dimensional spaces,
specific applications like the optimal control of partial
differential equations or even stochastic partial
differential equations---this last topic being motivated in
particular by the optimal control of the so-called Zakai's
operation, a well-known formulation of optimal control
problems of finite-dimensional diffusion processing with
partial observations---require some new developments of the
viscosity solutions theory in order to accommodate unbounded
terms in the equations
[\clii, Parts IV and V; \plld, \barbu, \barbuc, \sonc,
\cansoni, \tata].

The above rather vague and general comments on applications
already contain  many hints concerning promising directions
the theory of viscosity solutions  may take in the near
future.  In particular, the infinite-dimensional part of
the theory will most probably explode in view of the
unbounded  possible avenues of investigation.  Much more
progress is also to be expected for degenerate second-order
equations and boundary conditions.   In particular,  
progress
is to be made on existence questions for Dirichlet (and
state constraints)  boundary conditions---this might have a
considerable impact on various applications like models in
Finance.  Similarly, we expect progress on uniqueness and
regularity questions for uniformly elliptic second-order
equations.

It is also reasonable to hope that the insight gained by
viscosity solutions will help to devise efficient high-order
schemes for numerical approximations and prove their
convergence.

More specific developments should (and will) concern
questions of behavior at infinity, geometrical optics
problems, the use and the theory of discontinuous solutions,
and the investigation of \lq\lq second-order"
integrodifferential operators \dots.

Of course, the reader should not restrict his imagination
to the borders we drew above.

\equationno=0
\def\secno{A} 
\heading Appendix.  The proof of 
Theorem \tha \endheading

In this section we sketch the proof of Theorem \tha, which
we reproduce here for convenience.   \myProclaim Theorem
\tha.  Let $\co_i$ be a locally compact subset of $R^{N_i}$
for  $i=1,\ldots,k$, $$\co=\co_1\times\cdots\times\co_k,$$
$u_i\in \roman{USC}(\co_i),$ and $\gf$  be twice 
continuously
differentiable in a neighborhood of $\co$.    Set
$$w(x)=u_1(x_1)+\cdots+u_k(x_k)\quad\text{for 
}x=(x_1,\cdots,x_
k)\in\co,$$ and suppose $\hx=(\hx_1,\dotsc,\hx_k)\in\co$ 
is a
local maximum of $w-\gf$ relative to $\co$. Then for each
$\gep>0$ there exists $X_i\in\md{N_i}$  such that
$$(D_{x_i}\gf(\hx),X_i)\in\jtpcoic u_i(\hx_i) 
\quad\text{for }i=1,\ldots,k$$ and the block diagonal 
matrix with
entries $X_i$ satisfies \pre\eqapax 
$$-\left(\oo\gep +\vno
A\right)I\le \diag\le A+\gep A^2\eq$$ where
$A=D^2\gf(\hx)\in\m$, $N=N_1+\cdots+N_k$. 
\re\eqbhn

The sketch we give below will be somewhat  abbreviated, but
the main points  are all displayed and the reader should  be
able to provide any omitted details  from the indications
given.

As defined and used in the main text, $\roman{USC}(\co)$ 
consisted
of the upper semicontinuous  functions mapping $\co$ into
$\R$; however, at this point it is convenient to  allow the
value $``-\infty$," so we hereafter set  $$
\roman{USC}(\co)=\set{\ti{upper semicontinuous
functions}u\:\co\rtra\R\cup\set{-\infty}}.$$ The
convenience is well illustrated by the following reduction.
Below, conditions written  for the index $i$ are understood
to hold for $i=1,\ldots,k$.

\subheading{First Reduction}  We may as well assume that
$\co_i=\R^{N_i}$. Indeed, if not,  we first restrict $u_i$ 
to a 
compact neighborhood $K_i$ of $\hx_i$ in $\co_i$ and then 
extend  
the restriction to
$\R^{N_i}$ by $u_i(x_i)=-\infty$ if $x_i\not\in K_i$ 
(abusing notation 
and still writing $u_i$).  The compactness of $K_i$ 
guarantees that $u_i\in
\roman{USC}(\R^{N_i})$; one
then  checks that $\jtpcoic u_i(\hx_i)=\jtpc u_i(\hx_i)$
(provided $-\infty<u_i(\hx_i))$.  It is  clear that $\hx$
is still a local maximum of $w-\gf$ relative to
$\R^N=\R^{N_1+\cdots+N_k}$.  
  This device of putting functions
equal to $-\infty$ on unimportant sets will be used several 
times below to localize considerations to any compact
neighborhood of a point under consideration.

\subheading{Second Reduction}   We may as well assume
that $\hx=0$, $u_i(0)=0$, $\gf(x)=\oo2\ip{Ax,x}$ for some
$A\in\m$ is a pure quadratic, and  $0$ is a global maximum
of  $w-\gf$.  Indeed, a translation puts $\hx$ at the 
origin and then by replacing $\gf(x),\ u_i(x_i)$ by
$\gf(x)-(\gf(0)+\ip{D\gf(0),x})$ and 
$u_i(x_i)-(u_i(0)+\ip{D_{x_i}\gf(0),x_i})$, we reduce to the
situation   $\hx = D\gf(0)=0$, and $\gf(0)=u_i(0)=0$.  Then,
since $\gf(x)=\oo 2\ip{Ax,x}+o(|x|^2)$ where $A=D^2\gf(0)$
and $w(x)-\gf(x)\le w(0)-\gf(0)=0$ for small $x$, if
$\eta>0$ we will have $w(x)-\oo 2\ip{(A+\eta I)x,x}<0$ for
small $x\not=0$.  Globality of the (strict) maximum at 0 may
be achieved by localizing via the  first reduction.  If the
result holds in this case, we may use the bounds asserted
in Theorem \tha\ to pass to the limit  as $\eta  \downa 0$ 
to
obtain the  full  result.

Thus we have to prove the following

\myProclaim Theorem \tha\cprime.  Let $u_i\in 
\roman{USC}(\R^{N_i})$,  
$u_i(0)=0$, for $i=1,2,\ldots,k$,  $N=N_1+\dotsb+N_k$,
$A\in\m$, and    \pre\eqapa
$$w(x)=u_1(x_1)+\dotsb+u_k(x_k)\le
\tfrac12\ip{Ax,x}\ti{for}x=(x_1,\cdots,x_ k)\in\R^N.\eq$$ 
Then
for each $\gep>0$, there exists $X_i\in\md{N_i}$  such that
$$(0,X_i)\in\jtpcoic u_i(\hx_i)\quad  \text{for 
}i=1,\ldots,k,$$ and
the block diagonal matrix with entries $X_i$ satisfies
\pre\eqapw $$-\left(\oo\gep +\vno A\right)I\le \diag\le
A+\gep A^2.\eq$$

Here are the main steps in the proof:

 \step{1. 
Introduction  of the sup convolutions} If $\gep>0$, the
Cauchy-Schwarz inequality yields   $$\ip{Ax,x}
\le\ip{(A+\gep A^2)\xi,\xi}+ \left(\oo\gep+\vno
A\right)|x-\xi|^2\quad\text{for },x,\xi\in\RN.\eq$$ 
Putting \pre\eqapb
$$\gl=\oo\gep+\vno A\eq$$ and using \eqapa) and \eqapb) we
find \pre\eqapc
$$\left(u_1(x_1)-\ot\gl|x_1-\xi_1|^2\right)+\cdots+
\left(u_k(x_k)-\ot\gl|x_k-\xi_k|^2\right)\le \oo
2\ip{(A+\gep A^2)\xi,\xi}\eq$$ 
or $w(x)-(\gl/2)|x-\xi|^2\le
(1/2)\ip{(A+\gep A^2)\xi,\xi}$.  If necessary, we modify 
the 
definition of the $u_i$ off a neighborhood of 0 to assure
that they are bounded above. Put
 \pre\eqapd $$
\hw(\xi)=\sup_{x\in\RN}\left(w(x)-\ot\gl|x-\xi|^2\right),\ 
\hu_i(\xi_i)=\sup_{x_i\in\R^{N_i}}\left(u_i(x_i)-\ot%
\gl|x_i-\xi_i|^2\right)
\eq$$
so that $\hw(\xi)=\hu_1(\xi_1)+\dotsb+\hu_k(\xi_k)$. 
Observe that, choosing $x_i=0$, $\hu_i(0)\ge u_i(0)=0$ while
$\hu_1(0)+\dotsb+\hu_k(0)\le 0$, so $\hu_i(0)=0$.

\step{2. Two results on semiconvex functions: Aleksandrov's
Theorem and Jensen's Lemma}  The functions $\hw$ and $\hu_i$
are semiconvex; more precisely,  since the supremum of
convex functions is convex, $\hw(\xi)+(\gl/2)|\xi|^2$ and
$\hu_i(\xi)+(\gl/2)|\xi_i|^2$ are  convex.  We then call
$\gl$ a semiconvexity constant for $\hw$ and $\hu_i$.   This
is the reason we will employ two nontrivial facts about
semiconvex functions. The first assertion is a classical
result of Aleksandrov: 
\rpre\thapa \myProclaim Theorem \li. 
Let $\gf\:\RN\rtra \R$ be semiconvex.  Then  $\gf$ is twice 
differentiable almost everywhere on $\RN$.

\demo{Proof} \kern-.39pt Since \kern-.39pt we \kern-.39pt 
assume \kern-.39pt
that \kern-.39pt the \kern-.39pt members \kern-.39pt of 
\kern-.39pt our
\kern-.39pt audience \kern-.39pt who \kern-.39pt are 
\kern-.39pt not
\kern-.39pt functional analytically oriented have already 
seen a
proof of this
result, we will give one for those who are so oriented. In
the  course of proof we assume certain facts about convex
and Lipschitz functions that are known to  most readers. 
References will be given in the notes.

We may as well assume that $\gf$ is convex; a convex
function on $\RN$ is locally  Lipschitz and hence once
differentiable almost everywhere.
 We let 
$$F_1=\{x\in\RN\: D\gf(x)\text{ exists}\};$$ the letter $F$
being  used to indicate that $F_1$ has full measure (i.e.,
its complement is a  null set). The subdifferential $\sgf$
of $\gf$ is the set-valued function given by 
$$\hy\in\sgf(\hx) \quad\text{if }\gf(x)\ge
\gf(\hx)+\ip{\hy,x-\hx}\quad\text{for }x\in\RN.$$ It is 
clear that
$\sgf(x)=\set{D\gf(x)}$ whenever $\gf$ is differentiable at
$x$.  In particular, $\sgf$ is single valued almost
everywhere.  By convex analysis, $J=(I+\sgf)^{-1}$ is a 
single-valued mapping of $\RN$ into itself and is
nonexpansive (has Lipschitz  constant 1).   (Given
$x\in\RN$, the  solution $z$ of $z+\sgf(z)\ni x$ is the
value of $J$ at $x$;  $z$ is the minimum of $q\rtra
\gf(q)+\oo2|x-q|^2$.)  Since $\gf$ is defined on all of
$\RN$,    $\sgf(z)$ is nonempty for every $z$.  We claim
that for almost all $x\in\RN$ $D\gf$ is differentiable  at
$x$, i.e., there exists $A\in\m$ such that
$$D\gf(y)=D\gf(x)+A(y-x)+o(|y-x|)\iq{for}y\in F_1$$ for
almost all $x\in\RN$.  

To establish this, let  $$F_2=\{J(x): J\ti{is
differentiable at} x \ti{and} DJ(x)\text{ is 
nonsingular}\}.$$ 
Since $J$ is Lipschitz continuous and onto $\RN$ (because
$\sgf(z)$ is nonempty for  every $z$), $F_2$ has full
measure.  We use here the facts that
$\set{x: DJ(x)\ti{exists}}$ is a set of full measure, $J$ 
maps
null sets to null sets and $J(\{x:DJ(x)$ exists and is 
singular\}) is 
null.  It is
to be shown that  $\dgf$ is differentiable on  
$F_3=F_1\cap F_2$, which is a set of
full measure.  By the definition of $J$, 
$$\dgf(J(x))=x-J(x)\iq{for}J(x)\in F_3\subset F_1.$$ Thus
assume $J(x)+\dey\in F_1$ where $\dey$ is small. By
assumption, 
 $\dgf(J(x)+\dey)$ exists. Since $J$ is Lipschitz and
$DJ(x)$ is nonsingular, if $\dey$ is sufficiently small, 
there is a $\dex$ solving $J(x+\dex)=J(x)+\dey$, and  we may
choose $\dex$ to satisfy $|\dex|\le K|\dey|$  for some
constant $K$. Since $J$ is a contraction, we also have
$|\dey|\le |\dex|$ and  thus $|\dex|$ and $|\dey|$ are
comparable.  Then  
$$\eqalign{\dgf(J(x)+\dey)&=\dgf(J(x+\dex))=x+\dex-J(x+
\dex)\cr
&=\dgf(Jx)+(\dex- DJ(x)\dex)+o(\dex)\cr
&=\dgf(Jx)+(I-DJ(x))\dex+o(\dey).}$$ It remains to see that
$\dex=(DJ(x))^{-1}\dey+o(\dey)$; however, this follows from
the fact that  $|\dex|$ and $|\dey|$ are comparable and the
relation $J(x)+DJ(x)\dex+o(\dex)=J(x)+\dey.$ We conclude
that $D(D\gf)(J(x))$ exists and is $DJ(x)^{-1}-I$.

It remains to establish that  
$$
\gf(J(x)+\dey)=\gf(J(x))+\dgf(J(x))\dey+\tfrac
12\ip{(DJ(x)^{-1}-I)\dey,\dey}+o(|\dey|^2)$$
for $J(x)\in F_3$.  However, if we let
$\psi(\dey)=\gf(J(x)+\dey)$ and 
$\tpsi(\dey)=\gf(J(x))+\dgf(Jx)\dey+
\oo2\ip{(DJ(x)^{-1}-I)\dey,\dey}$,
we have  $\psi(0)=\tpsi(0)$ and for almost all small
$\dey$
$$
\aligned
D\psi(\dey)&=\dgf(J(x)+\dey)=\dgf(J(x))+
(DJ(x))^{-1}\dey-\dey+o(\dey)\\
&=D\tpsi(\dey)+o(\dey);\endaligned$$
thus $\Psi(\dey)=\psi(\dey)-\tpsi(\dey)$ is locally 
Lipschitz continuous  and
satisfies $\Psi(0)=0$ and $D\Psi(\dey)=o(\dey)$ for almost 
all small $\dey$.  It is
clear that then $\Psi(\dey)=o(|\dey|^2)$, whence the  
result.

The next result we will need concerning semiconvex
functions, which we call Jensen's lemma,  follows.  In the
statement, $B(x,r)$ is the closed ball of radius $r$
centered at $x$ and $B_r$ is the ball centered at the
origin.\enddemo 

\rpre\lej \myProclaim  Lemma  \li.  Let $\gf\:\RN\rtra
\R$ be semiconvex and $\hx$ be a strict local maximum point
of $\gf$.   For $p\in\RN$, set $\gf_p(x)=\gf(x)+\ip{p,x}$. 
Then for $r,\ \gd>0$,  $$K=\set{\!x\in 
B(\hx,r)\!:\!\ti{there
exists}p\in B_\gd \ti{for which} \gf_p \ti{has a local 
maximum
at}x}$$	 has positive measure.

\demo{Proof} We assume that $r$ is so small that $\gf$ has 
$\hx$ as a
unique maximum point in $B(\hx,r)$ and assume for the 
moment that 
$\gf$ is $C^2$. It follows from this
that if $\gd$ is sufficiently small and $p\in B_\gd$,  then
every maximum of $\gf_p$ with respect to $B(\hx,r)$ lies in
the interior of $B(\hx,r)$.   Since $D\gf+p=0$ holds at
maximum points of $\gf_p$, $D\gf(K)\supset B_\gd$. Let
$\gl\ge0$ and $\gf(x)+(\gl/2)|x|^2$ be convex; we then 
have $-\gl I\le D^2\gf$; moreover,
on $K$, $D^2\gf\le 0$ and then $$-\gl I\le D^2\gf(x)\le
0\iq{for}x\in K.$$ In particular, $|\det D^2\gf(x)|\le
\gl^N$ for  $x\in K$. Thus 
$$ \text{meas}(B_\gd)\le\text{meas}(D\gf(K))\le \int_K
|\det D^2\gf(x)|dx\le \text{meas}(K)|\gl|^N$$ 
(see the notes) and we have
a lower bound on the measure of $K$ depending only on $\gl$.

In the general case, in which $\gf$ need not be smooth, we
approximate it via  mollification with smooth functions
$\gf_\gep$ that have the same semiconvexity  constant $\gl$
and  that converge uniformly to $\gf$ on $B(\hx,r)$.  The
corresponding  sets $K_\gep$ obey the above estimates for
small $\gep$ and $$K\supset
\bigcap_{n=1}^\infty\bigcup_{m=n}^\infty K_{1/m}$$ is
evident.  The result now follows.\enddemo

\step{3. A  consequence of Step 2 and magic properties of
sup convolution}
Lemma \leat\ below applies to $\hw$ of Step 1 and we shall
see it provides us with matrices $X_i\in\md{N_i}$ such that
$(0,X_i)\in\jtpc \hu_i(0)$ and \eqapw) holds. 

   We use the notation  
    \pre\eqw
     $$J^2 f(z)=\jtp f(z)\cap \jtm f(z)\eq$$
     from which one defines $\overline J^2$ analogously to
$\jtpc, \ \jtmc$.  Note that
     $(p,X)\in J^2 f(x)$ amounts to
     $$f(y)=
f(x)+\ip{p,y-x}+\tfrac 12\ip{X(y-x),y-x}+
o(|x-y|^2)\quad\text{as }y\rtra x;$$
     i.e., $f$ is twice differentiable at $x$ and $p=Df(x),\
X=D^2f(x)$.
      
      \rpre\leat
      \myProclaim Lemma \li .  If $f\in C(\RN), \  B\in\m 
$,  
$f(\xi)+ (\gl/2) |\xi|^2$ is
      convex and
$\max_{\RN}(f(\xi)-\oo2\ip{B\xi,\xi})=f(0)$, then there is
an $X\in\m$
      such that $(0,X)\in\overline J^2 f(0)$ and $-\gl I\le
X\le B$.
       
       \demo{Proof}  Clearly  
$f(\xi)-\oo2\ip{B\xi,\xi}-|\xi|^4$ has
a strict maximum at $\xi=0$.  By Aleksandrov's Theorem  
and the convexity assumption,
$f$ is twice differentiable a.e. and
       then by Jensen's Lemma, for every $\gd>0$ there
exists $q_\gd\in\RN$
       with $|q_\gd|\le\gd$ such that $f(\xi)
+\ip{q_\gd,\xi}-\oo2\ip{B\xi,\xi}-|\xi|^4$ has a
       maximum at a point $\xi_\gd$ with $|\xi_\gd|\le \gd$ 
and $f(\xi)$ is twice
       differentiable at $\xi_\gd$.   From  calculus,
$|q_\gd|,|\xi_\gd|\le\gd$, and the
       convexity of $f(\xi)+(\gl/2) |\xi|^2$,
       $$ Df(\xi_\gd)=O(\gd),\qquad  -\gl I\le 
D^2f(\xi_\gd)\le
B+ O(\gd^2).$$
       Observing that $(Df(\xi_\gd),D^2f(\xi_\gd))\in J^2
f(x_\gd)$, we  conclude the
       proof upon passing to a subsequential limit (selected
so that $D^2f(\xi_\gd)$ is
       convergent) as $\gd\downa0$.
	
	Applying Lemma \leat\  to $\hw$ from \eqapd) with
	$B=A+\gep A^2$ and noting that the values of $J^2 \hw$ are
formed as one expects
	given the  representation
	$\hw=\hu_1+\dotsb+\hu_k$, we obtain $(0,X_i)\in \overline
	J^2\hu_i(0)$, $i=1,\ldots,k$ such that \eqapw) holds.  The
final step in the proof of
	Theorem \tha\cprime\ is provided by applying the next 
lemma to
each $u_i$ to conclude
	that $(0,X_i)\in \jtpc u_i(0)$.
	
	Here are magical properties of the sup convolution.\enddemo

	\rpre\leb
	\myProclaim Lemma  \li.  Let $\gl>0$,  $v\in 
\roman{USC}(\R^M)$
	bounded above,  and
$$\hv(\xi)=\sup_{x\in\R^M}(v(x)-(\gl/2)|x-\xi|^2).$$
If $\eta\!, q\in\R^M$\!, $Y\in\md M$\!, and 
$(q,Y)\in\jtp\hv(\eta)$,
then 
$$(q,Y)\in \jtp
	v(\eta+q/\gl)\quad\text{and}\quad
\hv(\eta)+(1/2\gl)|q|^2=v(\eta+q/\gl).
$$  
In consequence,
	if $(0,Y)\in\jtpc \hv(0)$, then $(0,Y)\in\jtpc v(0)$.

	 \demo{Proof}  We assume $(q,Y)\in\jtp\hv(\eta)$, and let 
$y\in\RN$
be a point such that $
\hv(\eta)=v(y)-(\gl/2)|y-\eta|^2$.   Then for any $x,\
\xi\in\RN$
	 \pre\eqzw
	 $$\eqalign{v(x)-\ot
\gl|\xi-x|^2& \le\hv(\xi)\cr
& \le\hv(\eta)+\ip{q,\xi-\eta}\cr
&\hphantom{=}+\oo
2\ip{Y(\xi-\eta),\xi-\eta}+o(|\xi-\eta|^2)\cr
&=v(y)-\ot\gl|y-\eta|^2
	 +\ip{q,\xi-\eta}\cr
&\hphantom{=}+\oo2\ip{Y(\xi-\eta)\xi-\eta}+
o(|\xi-\eta|^2)\cr
&=v(y)-\ot\gl|y-\eta|^2
	 +\ip{q,\xi-\eta}+O(|\xi-\eta|^2).\cr}\eq$$
	 Putting $\xi=x-y+\eta$ in the appropriate one of these 
relations yields
	 \pre\eqzax
	 $$
v(x)\le
v(y)+\ip{q,x-y}+\oo2\ip{Y(x-y),x-y}+o(|x-y|^2),\eq$$
	 while substituting $x=y$ and
$\xi=\eta+\ga(\gl(\eta-y)+q)$   yields
	 \pre\eqzbx
	 $$0\le\ga|\gl(\eta -y)+q|^2+O(\ga^2).\eq$$
	 The former simply says that $(q,Y)\in\jtp v(y)$ while the
latter with small  $\ga<0$
	 implies that $\gl(\eta-y)+q=0$ or $y=\eta+q/\gl$ as
claimed.  The relation
	 $\hv(\eta)+(1/2\gl)|q|^2=v(\eta+q/\gl)$ follows at once. 
Assume now that
	 $(q_n,Y_n)\in \jtp\hv(\xi_n)$ and
$(\xi_n,\hat v(\xi),q_n,Y_n)\rtra(0,0,0,Y)$; by the
	 foregoing,   $(q_n,Y_n)\in\jtp v(\xi_n+q_n/\gl)$ and
	 $v(\xi_n+q_n/\gl)=\hv(\xi_n)+(1/2\gl)|q_n|^2$.  By the
definitions,
	 $v(0)\le \hv(0)$ and from this,  the upper semicontinuity
and the foregoing, we have
	 $$v(0)\ge\limsup_{n\rtra\infty}v(\xi_n+{q_n\over \gl})=
	 \limsup_{n\rtra\infty}(\hv(\xi_n)+
\oo{2\gl}|q_n|^2)=\hv(0)\ge
v(0),$$
	 which provides the final piece of information we needed to
conclude that
	 $(0,Y)\in\jtpc v(0).$\enddemo

\rpre\reapa 
\rem{Remark {\rm \li}} We close with a final
and more elegant reformulation of what was proved above:  if
the $u_i$ satisfy the conditions of  Theorem \tha,
$w(x)=u_1(x_1)+\cdots+u_k(x_k)$ and
$((p_1,\ldots,p_k),A)\in\jtpoc w(\hx)$, then for 
 each $\gep>0$, there exists $X_i\in
\cs(N_i)$ such that $(p_i,X_i)\in\jtpcoic u_i(\hx_i)$ and 
\eqapax) holds. 
\endrem

\subheading{Notes on the appendix}
Above, we provided (for the first time) a self-contained 
proof of
Theorem \tha.  Except for  the two auxiliary results on 
semiconvex
functions, the main tool is the so-called sup convolution. 
  This
approximation procedure (more often in the guise of inf 
convolution)
is well known in functional analysis and, in particular, 
in  convex
analysis and the theory of maximal monotone operators 
(see, for
example, the text [\brezis] of  H. Brezis).  It was 
noticed in
[\lali] that it may provide an efficient regularization 
procedure 
for (even degenerate) elliptic equations; some of its 
properties are
given there.  See also [\jls].   Its ``magical properties" 
can be
seen as related to the Lax formula for the solution of  
$${\partial
w\over\partial t}-\oo2|\nabla w|^2=0\quad\text{for 
}x\in\RN,\ t\ge 0,\
w|_{t=0}=v\ti{on}\RN,$$ which is  $$
w(x,t)=\sup_{y}\set{v(y)-\oo{2t}|x-y|^2}.$$ Indeed, the 
coincidence
of this solution formula and solutions produced by the  
method of
charteristics leads to the properties used.  Of course, 
this is a
heuristic  connection, since  characteristic methods 
require too
much regularity to be rigorous here. 

The inf convolution  can also be seen as a nonlinear 
analogue of the
standard mollification when  replacing the ``linear 
structure of
$L^2$ and its duality" by the ``nonlinear structure of 
$L^\infty$ 
or $C$."  One can also interpret this analogy in terms of 
the
so-called exotic algebra ($\R$,  max, +).  

Theorem \thapa\ is a classical result of A. D. Aleksandrov 
[\alea]. 
The proof given here  is in fact a slightly stronger  form 
of this
result ($\dgf$ is  differentiable a.e.) and our proof 
follows F.
Mignot [\mignot] (although we have made his proof more 
complicated for 
pedagogical reasons).  In the proof we used the fact that 
Lipschitz
functions are differentiable a.e., which is called  
Rademacher's Theorem.
A proof may be found in L. C. Evans and R. Gariepy 
[\evgar] and F. Mignot 
[\mignot].  We also used that Lipschitz functions map 
null sets to null sets and the solvability of
$J(x+\dey)=J(x)+\dex$ when $DJ(x)$ exists and is 
nonsingular. 
These  are proved in, respectively,  Lemma 7.25 and the 
proof of
Theorem 7.24 in  W. Rudin [\rudin]; see also [\mignot].
We also used the fact that
$J(\{x:DJ(x)\ti{exists and is singular\!\!}\!\!\})$ is 
null.   This
follows  from the general formula $$\int_{\RN} \#(A\cap
J^{-1}(y))\,dy = \int_A |\det DJ(x)|\,dx$$   where \# is
counting measure.  This formula  holds for Lipschitz 
continuous
functions and  measurable sets $A\subset\RN$ and is a  
special case
of the Area Formula for Lipschitzian maps between 
Euclidean spaces. 
The  Area Formula may be found in L. C. Evans and R. Gariepy
[\evgar] and H. Federer [\fed]. 

Jensen's Lemma (Lemma \lej) is a variation on an aspect of 
the theme known as ``Aleksandrov's maximum  principle" 
(see [\aleb, \alec, \baka, \bony, \pllbony,  \pucci]).

\enddocument